\documentclass[mnsc,nonblindrev]{informs3} 

\OneAndAHalfSpacedXI

\usepackage{natbib}
 \bibpunct[, ]{(}{)}{,}{a}{}{,}%
 \def\bibfont{\small}%

\newtheorem{observation}{Observation}

\TheoremsNumberedThrough     
\ECRepeatTheorems

\EquationsNumberedThrough    

\MANUSCRIPTNO{}
                 
\usepackage{enumitem}
\usepackage{booktabs}
\usepackage{graphicx}
                 
\usepackage{bm}        

\usepackage[ruled,vlined]{algorithm2e}
\include{pythonlisting}
\usepackage{pdflscape}

\usepackage{multirow}

\newcommand{\field}[1]{\ensuremath{\mathbb{#1}}}
\newcommand{\sets}[1]{\ensuremath{\mathcal{#1}}}

\newcommand{\reals}{\ensuremath{\field{R}}} 
\newcommand{\integers}{\ensuremath{\field{Z}}} 
\newcommand{\1}{\ensuremath{{\rm \mathbf e}}} 
\newcommand{\I}[1]{\ensuremath{\mathbb{I}\left(#1\right)}} 

\newcommand{\subjectto}{\text{\rm subject to}} 

\DeclareMathOperator{\st}{s.t.}

\newcommand{\minimize}{\ensuremath{\mathop{\mathrm{minimize}}\limits}}
\newcommand{\maximize}{\ensuremath{\mathop{\mathrm{maximize}}\limits}}

\newcommand{\newpv}[1]{{\color{black}{#1}}}

\renewcommand{\arraystretch}{1.5}

\DeclareRobustCommand{\cpluspluslogo}{\hbox{C\hspace{-0.5ex}
                       \protect\raisebox{0.5ex}
                       {\protect\scalebox{0.67}{++}}}}

\def\ROCPP{RO\cpluspluslogo}

\begin{document}


\RUNAUTHOR{Vayanos, Georghiou, Yu}

\RUNTITLE{Robust Optimization with Decision-Dependent Information Discovery}

\TITLE{Robust Optimization with \\ Decision-Dependent Information Discovery}

\ARTICLEAUTHORS{%
\AUTHOR{Phebe Vayanos}
\AFF{University of Southern California, Center for Artificial Intelligence in Society, \EMAIL{phebe.vayanos@usc.edu}} 
\AUTHOR{Angelos Georghiou}
\AFF{University of Cyprus, Department of Business and Public Administration, \EMAIL{georghiou.angelos@ucy.ac.cy}
\AUTHOR{Han Yu}
\AFF{University of Southern California, Center for Artificial Intelligence in Society, \EMAIL{hyu376@usc.edu}}
}
} 

\ABSTRACT{%

Robust optimization (RO) is a popular paradigm for modeling and solving two- and multi-stage decision-making problems affected by uncertainty. \newpv{In many real-world applications, such as R\&D project selection, production planning, or preference elicitation for product or policy recommendations, the \emph{time of information discovery} is decision-dependent and the uncertain parameters only become observable after an often costly investment.} Yet, \newpv{most of the literature on robust optimization} assumes that the uncertain parameters can be observed \emph{for free} and that the sequence in which they are revealed is \emph{independent} of the decision-maker's actions. To fill this gap \newpv{in the practicability of RO}, we consider two- and multi-stage robust optimization problems in which part of the decision variables control the time of information discovery. Thus, \emph{information} available at any given time is \emph{decision-dependent} and can be discovered (at least in part) by making strategic exploratory investments in previous stages. We propose a novel \emph{dynamic} formulation of the problem and prove its correctness. We leverage our model to provide a solution method inspired from the $K$-adaptability approximation, whereby $K$ candidate strategies for each decision stage are chosen \emph{here-and-now} and, at the beginning of each period, the best of these strategies is selected \emph{after} the uncertain parameters that were chosen to be observed are revealed. We reformulate the problem as a finite \newpv{mixed-integer (resp.\ bilinear)} program \newpv{if none (resp.\ some of the) decision variables are real-valued. This finite program is} solvable with off-the-shelf solvers. We generalize our approach to the minimization of piecewise linear convex functions. We demonstrate the effectiveness of our method \newpv{in terms of interpretability, optimality, and speed} on synthetic instances of the Pandora box problem, \newpv{the preference elicitation problem with real-valued recommendations,} the best box problem, and the R\&D project portfolio optimization problem. Finally, we evaluate it on an instance of the active preference elicitation problem used to recommend kidney allocation policies to policy-makers at the United Network for Organ Sharing based on real data from the U.S.\ Kidney Allocation System.

}


\maketitle

%



\section{Introduction}
\label{sec:introduction}


\subsection{Background \& Motivation}
\label{subsec:motivation}

Over the last two decades, robust optimization has emerged as a popular approach for decision-making under uncertainty in both \emph{single-} and \emph{multi-stage} settings, see e.g.,~\cite{BenTal_Book, RobustLP_UncertainData,RSols_UncertainLinear_Programs_BT,BenTalNemirovski_RCO,Robust_LO_generalNorms,Price_Robustness, AdjustableRSols_uncertainLP,Bertsimas_polynomial_policies,ZHS16:fourier_motzkin,ConstraintSampling_VKR,Goyal_OnPowerLimitations,XB16:copositive}. 
In multi-stage models, the uncertain parameters are revealed sequentially as time progresses and the decisions are allowed to depend on all the information made available in the past. Mathematically, decisions are modeled as functions of the history of observations, thus capturing the \emph{adaptive} and \emph{non-anticipative} nature of the decision process. 

Most models and solution approaches in multi-stage robust optimization are tailored to problems where the uncertain parameters are \emph{exogenous}, being independent of the decision-maker's actions. In particular, they assume that uncertainties can be observed \emph{for free} and that the \emph{sequence} in which they are revealed \emph{cannot be influenced} by the decision-maker. Yet, these assumptions fail to hold in many real-world applications where the \emph{time of information discovery} is decision-dependent and the uncertain parameters only become observable after an often costly investment. Mathematically, some binary \emph{measurement} (or \emph{observation}) decisions control the time of information discovery and the non-anticipativity requirements depend upon these decisions, severely complicating solution.


\subsubsection{Motivating Applications.}


We now detail several applications areas where the time of revelation of the uncertain parameters is decision-dependent.

\paragraph{{R\&D Project Portfolio Optimization.}} Research and development firms typically maintain long pipelines of candidate projects whose returns are uncertain, see \cite{Solak_RandDportfolios}. For each project, the firm can decide whether and when to start it and the amount of resources to be allocated to it. The return of each project will only be revealed once the project is completed. Thus, project start times and resource allocation decisions impact the time of information discovery in this problem.

\paragraph{{Clinical Trial Planning.}} Pharmaceutical companies typically maintain long R\&D pipelines of candidate drugs, see e.g., \cite{Colvin_ClinicalTrials}. Before any drug can reach the marketplace it needs to pass a number of costly clinical trials whose outcome (success/failure) is uncertain and will only be revealed after the trial is completed. Thus, the decisions to proceed with a trial control the time of information discovery in this problem. 

\paragraph{{Offshore Oilfield Exploitation.}} Offshore oilfields consist of several reservoirs of oil whose volume and initial deliverability (maximum initial extraction rate) are uncertain, see e.g., \cite{Jonsbraten_thesis,GoelGrossmanGasFields}, and \cite{DDI_VKB}. While seismic surveys can help estimate these parameters, current technology is not sufficiently advanced to obtain accurate estimates. In fact, the volume and deliverability of each reservoir only become precisely known if a very expensive oil platform is built at the site and the drilling process is initiated. Thus, the decisions to build a platform and drill into a reservoir control the time of information discovery in this problem.

\paragraph{{Production Planning.}} Manufacturing companies can typically produce a large number of different items. For each type of item, they can decide whether and how much to produce to satisfy their demand given that certain items are substitutable, see e.g.~\cite{Jonsbraten_DecDepRandElmts}. The production cost of each item type is unknown and will only be revealed if the company chooses to produce the item. Thus, the decisions to produce a particular type of item control the time of information discovery in this problem.

\paragraph{{Active Preference Elicitation.}} Preference elicitation refers to the problem of developing a decision support system capable of generating recommendations to a user, thus assisting in decision making. In active preference elicitation, one can ask users a (typically limited) number of questions from a potentially large set before making a recommendation, see e.g., \cite{Vayanos_ActivePreferences}. The answers to the questions are initially unknown and will only be revealed if the particular question is asked. Thus, the choices of questions to ask control the time of information discovery in this problem. 


\subsection{Literature Review}
\label{subsec:literature_review}

\paragraph{Decision-Dependent Information Discovery.} Our paper relates to research on optimization problems affected by uncertain parameters whose time of revelation is decision-dependent and which originates in the literature on stochastic programming. The vast majority of these works assumes that the uncertain parameters are discretely distributed. In such cases, the decision process can be modeled by means of a finite scenario tree whose branching structure depends on the binary measurement decisions that determine the time of information discovery. This research began with the works of~\cite{Jonsbraten_DecDepRandElmts} and \cite{Jonsbraten_thesis}. \cite{Jonsbraten_DecDepRandElmts} consider the case where all measurement decisions are made in the first stage and propose a solution approach based on an implicit enumeration algorithm. \cite{Jonsbraten_thesis} generalizes this enumeration-based framework to the case where measurement decisions are made over time. More recently, \cite{GoelGrossmanGasFields} showed that stochastic programs with discretely distributed uncertain parameters whose time of revelation is decision-dependent can be formulated as deterministic mixed-binary programs whose size is exponential in the number of endogenous uncertain parameters. To help deal with the ``curse of dimensionality,'' they propose to precommit all measurement decisions, i.e., to approximate them by here-and-now decisions, and to solve the multi-stage problem using either a decomposition technique or a folding horizon approach. Later, \cite{GoelGrossman_ClassStochastic_DDU,GoelGrossman_NovelBB_GasFields}, and \cite{Colvin_Pharmaceutical} propose optimization-based solution techniques that truly account for the adaptive nature of the measurement decisions and that rely on branch-and-bound and branch-and-cut approaches, respectively. Accordingly, \cite{Colvin_Pharmaceutical} and \cite{GuptaGrossman_SolutionStrategies} have proposed iterative solution schemes based on relaxations of the non-anticipativity constraints for the measurement variables. Our paper most closely relates to the work of~\cite{DDI_VKB}, wherein the authors investigate two- and multi-stage stochastic and robust programs with decision-dependent information discovery that involve continuously distributed uncertain parameters. They propose a decision-rule based approximation approach that relies on a prepartitioning of the support of the uncertain parameters. Since this approach applies in our context, we will benchmark against it in our experiments.

\paragraph{Robust Optimization with Decision-Dependent Uncertainty Sets.} Our work also relates to the literature on robust optimization with uncertainty sets parameterized by the decisions. Such problems capture the ability of the decision-maker to influence the set of possible realizations of the uncertain parameters and have been investigated by~\cite{Spacey2012,Nohadani2016,Nohadani2017,ZHANG2017249}, and \cite{BV_dynamic_pricing}. These models do not apply in our context since they do not capture the ability of the decision-maker to influence the \emph{information} available. In particular, the problems investigated by~\cite{Spacey2012,Nohadani2016}, and \cite{Nohadani2017} are all single-stage, while \emph{problems with decision-dependent information discovery are inherently sequential in nature}.


\paragraph{Robust Optimization with Binary Adaptive Variables.} Two-stage, and to a lesser extent also multi-stage, robust binary optimization problems have received considerable attention in the recent years. One stream of works proposes to restrict the functional form of the recourse decisions to functions of benign complexity, see~\cite{Bertsimas2017} and \cite{Georghiou2015,bertsimas2018binary}. A second stream of work relies on partitioning the uncertainty set into finite sets and applying constant decision rules on each partition, see~\cite{DDI_VKB,Dunning2016,Postek2016MultistageSet,BV_dynamic_pricing}. The last stream of work investigates the so-called $K$-adaptability counterpart of two-stage problems, see \cite{Caramanis_FiniteAdaptability,Hanasusanto2015,subramanyam2017k,CHASSEIN2019308}, and~\cite{Rahmattalabi_SuicidePreventionNIPS}. In this approach, $K$ candidate policies are chosen here-and-now and the best of these policies is selected after the uncertain parameters are revealed. Most of these papers assume that the uncertain parameters are \emph{exogenous} in the sense that they are \emph{in}dependent of the decision-maker's actions. Our paper most closely relates to the works of \cite{Caramanis_FiniteAdaptability} and \cite{Hanasusanto2015}. 
It generalizes and subsumes the approach from~\cite{Hanasusanto2015} to problems with \emph{decision-dependent information discovery}, to \emph{multi-stage} problems, and to problems with \emph{piecewise linear convex objective}.

\paragraph{Stochastic Probing.} Our paper also fits in a line of work on stochastic probing in the computer science literature, see~\cite{Gupta2016_SP,Gupta2017_SP} and~\cite{Singla_PhDThesis}. Here, the problem consists of a set of elements with uncertain value whose distribution is known but whose realization becomes observable only after the element is probed. However, probing is costly (incurs a cost or consumes budget) and irrevocable and the goal is to choose the set of elements to probe and the order in which to probe them to maximize profit (e.g., the value of the item with the highest value that has been probed). Concrete examples include the best box problem and the Pandora box problem, see e.g.,~\cite{Singla_PhDThesis}. The techniques presented in this stream of work do not apply to the case where the distributions are unknown, to general optimization problems with decision-dependent information discovery, nor to problems with general, potentially uncertain, constraints.



\paragraph{Worst-Case Regret Optimization.} Finally, our work relates to two-stage worst-case absolute regret minimization problems, see e.g., \cite{Assavapokee2008a,Assavapokee2008b,Zhang2011,Jiang2013,NG2013483,Chen2014,NING2018425}, and~\cite{poursoltani19Regret}. To the best of our knowledge, our paper is the first to investigate worst-case regret minimization problems in the presence of uncertain parameters whose time of revelation is decision-dependent.





\subsection{Proposed Approach and Contributions}
\label{subsec:contributions}

We now summarize our approach and main contributions in this paper:
\begin{enumerate}[label=\textit{(\alph*)}]
    \item We consider general two- and multi-stage robust optimization problems with decision-dependent information discovery. These encompass as special cases the R\&D project portfolio optimization problem, the Pandora box problem (which can be used to model job candidate selection and house hunting, among others), the active preference elicitation problem, and many more. To the best of our knowledge, only one other paper in the literature studies such problems in the robust optimization setting. We propose novel ``min-max-min-max-...-min-max'' reformulations of these problems and prove correctness of our formulations. These reformulations unlock new approximate (and potentially also exact) solution approaches for addressing problems with decision-dependent information discovery.
    \item We leverage our new reformulations to propose a solution approach based on the $K$-adaptability approximation, wherein $K$ candidate strategies are chosen here-and-now and the best of these strategies is selected after the uncertain parameters that were chosen to be observed are revealed. This approximation allows us to control the trade-off between complexity and solution quality by tuning a single design parameter, $K$.  We propose practicable reformulations of the $K$-adaptability counterpart of problems with decision-dependent information discovery in the form of moderately sized finite programs solvable with off-the shelf solvers. These programs can be written equivalently as mixed-binary linear programs if all decision-variables are binary. Our reformulations subsume those from the literature that apply only to two-stage problems with exogenous uncertain parameters.
    \item We generalize the $K$-adaptability approximation scheme to multi-stage problems and to problems with piecewise linear convex objective function. The piecewise linear convex objective enables us, among others, to address worst-case absolute regret minimization problems. These generalizations and associated algorithm that we provide apply also to problems with exogenous uncertain parameters.
    \item We perform a wide array of experiments on the R\&D project portfolio selection problem, \newpv{the preference elicitation problem with real-valued recommendations,} the best box selection problem, Pandora's box problem, and the preference elicitation problem. We show that our proposed approach outperforms the state-of-the-art in the literature in terms of \newpv{interpretability, optimality, and speed. Indeed, our approach reduces the number of subsets in the recourse strategy by a factor of 3, improves the quality of the returned solution by a factor of 1.9, and results in an 8.5$\times$ speed-up}.
    We perform a case study \newpv{showcasing the benefits of our approach} on real data from the U.S.\ Kidney Allocation System (KAS) to recommend policies that meet the needs of policy-makers at the Organ Procurement and Transplantation Network (OPTN) and the United Network for Organ Sharing (UNOS), the lead agency in charge of allocating organs for transplantation in the United States.\footnote{See \url{https://www.srtr.org}, \url{https://optn.transplant.hrsa.gov}, and \url{https://unos.org}.}
\end{enumerate}

\subsection{Organization of the Paper and Notation}

The paper is organized as follows. Sections~\ref{sec:problem_formulation_exo} and~\ref{sec:problem_formulation_endo} introduce two-stage robust optimization problems with exogenous uncertainty and with decision-dependent information discovery (DDID), respectively. In particular, Section~\ref{sec:problem_formulation_endo} introduces our novel formulation. Section~\ref{sec:kadaptability_overview} proposes reformulations of the $K$-adaptability counterparts of problems with DDID as finite programs solvable with off-the-shelf solvers. Section~\ref{sec:pwl_convex_objective} generalizes the $K$-adaptability approximation to problems with piecewise linear convex objective and proposes an efficient solution procedure. Section~\ref{sec:multistage} generalizes the $K$-adaptability approximation to multi-stage problems. Section~\ref{sec:computational_studies} presents computational results on synthetic instances of the two-stage R\&D project portfolio optimization problem, the two-stage best box selection problem, and the multi-stage Pandora's box problem. Finally, Section~\ref{sec:preference_elicitation_KAS} formulates the preference elicitation problem for learning the preferences of policy-makers at the OPTN/UNOS as a two-stage robust problem with decision-dependent information discovery, and presents numerical results on real data from the U.S.\ Kidney Allocation System. The proofs of all statements can be found in the Electronic Companion to the paper. Proposed extensions to our methods, algorithms, and speed-up strategies are also deferred to the Electronic Companion.


\paragraph{Notation.} Throughout this paper, vectors (matrices) are denoted by boldface lowercase (uppercase) letters. The $k$th element of a vector ${\bm x} \in \reals^n$ ($k \leq n$) is denoted by ${\bm x}_k$. Scalars are denoted by lowercase letters, e.g., $\alpha$ or $u$. For a matrix ${\bm H} \in \mathbb R^{n \times m}$, we let $[{\bm H}]_k \in \mathbb R^m$ denote the $k$th row of ${\bm H}$, written as a column vector. We let $\mathcal L_{n}^k$ denote the space of all functions from $\reals^n$ to $\reals^k$. Accordingly, we denote by $\mathcal B_{n}^k$ the spaces of all functions from $\reals^n$ to $\{0,1\}^k$. Given two vectors of equal length, ${\bm x}$, ${\bm y} \in \mathbb R^n$, we let ${\bm x} \circ {\bm y}$ denote the Hadamard product of the vectors, i.e., their element-wise product. With a slight abuse of notation, we may use the maximum and minimum operators even when the optimum may not be attained; in such cases, the operators should be understood as suprema and infima, respectively. We use the convention that a decision is feasible for a minimization problem if and only if it attains an objective that is $<+\infty$. Finally, for a logical expression $E$, we define the indicator function $\I{E}$ as $\I{E}:=1$ if $E$ is true and 0 otherwise.


\section{Two-Stage RO with Exogenous Uncertainty}
\label{sec:problem_formulation_exo}


To motivate our formulation from Section~\ref{sec:problem_formulation_endo}, we introduce two equivalent models of two-stage robust optimization with \emph{exogenous uncertainty} from the literature and discuss their relative merits. 

In two-stage robust optimization with \emph{exogenous} uncertainty, first-stage (or here-and-now) decisions ${\bm x} \in \sets X \subseteq \reals^{N_x}$ are made today, \emph{before} any of the uncertain parameters are observed. Subsequently, all of the uncertain parameters ${\bm \xi} \in \Xi \subseteq \reals^{N_\xi}$ are revealed. Finally, once the realization of ${\bm \xi}$ has become available, second-stage (or wait-and-see) decisions ${\bm y} \in \sets Y \subseteq \reals^{N_y}$ are selected. We assume that the uncertainty set $\Xi$ is a non-empty bounded polyhedron expressible as $\Xi := \{ {\bm \xi} \in \reals^{N_\xi} \; : \; {\bm A} {\bm \xi} \leq {\bm b} \}$ for some matrix ${\bm A} \in \reals^{R \times N_\xi}$ and vector ${\bm b} \in \reals^R$. As the decisions ${\bm y}$ are selected after the uncertain parameters are revealed, they are allowed to \emph{adapt} or \emph{adjust} to the realization of ${\bm \xi}$. In the literature, there are two formulations of generic two-stage robust problem with exogenous uncertainty: they differ in the way in which the ability of ${\bm y}$ to adapt to~${\bm \xi}$ is modeled.


\paragraph{Decision Rule Formulation.}

In the first model, one optimizes today over both the here-and-now decisions ${\bm x}$ and over recourse actions ${\bm y}$ to be taken in each realization of ${\bm \xi}$. The decision ${\bm y}$ is modeled as a function (or \emph{decision rule}) of ${\bm \xi}$ that is selected today, along with ${\bm x}$. Under this paradigm, a two-stage linear robust problem with exogenous uncertainty is expressible as:
\begin{equation}\renewcommand{\arraystretch}{1.5}
    \begin{array}{cl}
         \minimize & \quad \displaystyle \max_{{\bm \xi} \in \Xi} \;\;\;  {\bm \xi}^\top {\bm C} \; {\bm x} + {\bm \xi}^\top {\bm Q} \; {\bm y}({\bm \xi}) \\
         \subjectto & \quad {\bm x} \in \sets X, \; {\bm y} \in \mathcal L_{N_\xi}^{N_y} \\
         & \quad \!\! \left. \begin{array}{l} 
         {\bm y}({\bm \xi}) \in \sets Y  \\
         {\bm T} {\bm x} + {\bm W}{\bm y}({\bm \xi}) \leq {\bm H}{\bm \xi} 
         \end{array} \quad \right\} \quad \forall {\bm \xi} \in \Xi,
    \end{array}
\label{eq:exo_1}
\end{equation}
where ${\bm C} \in \reals^{N_\xi \times N_x}$, ${\bm Q} \in \reals^{N_\xi \times N_y}$, ${\bm T} \in \reals^{L \times N_x}$, ${\bm W} \in \reals^{L \times N_y}$, and ${\bm H} \in \reals^{L \times N_\xi}$. We assume that the objective function and right hand-sides are \emph{linear} in ${\bm \xi}$. We can account for \emph{affine} dependencies on~${\bm \xi}$ by introducing an auxiliary uncertain parameter ${\bm \xi}_{N_\xi+1}$ restricted to equal unity.


\paragraph{Min-Max-Min Formulation.}

In the second model, only ${\bm x}$ is selected today and the recourse decisions ${\bm y}$ are optimized explicitly, in a \emph{dynamic} fashion, \emph{after} nature is done making a decision. Under this model, a two-stage robust problem with exogenous uncertainty is expressible as:
\begin{equation}\renewcommand{\arraystretch}{1.5}
    \begin{array}{cl}
         \minimize & \quad \displaystyle \max_{{\bm \xi} \in \Xi} \;\;\;  \left[ {\bm \xi}^\top {\bm C} \; {\bm x} + \min_{ {\bm y} \in \sets Y } \; \; \left\{ {\bm \xi}^\top {\bm Q} \; {\bm y} \; : \; {\bm T} {\bm x} + {\bm W}{\bm y} \leq {\bm H}{\bm \xi}  \right\} \right] \\
         \subjectto & \quad {\bm x} \in \sets X. 
    \end{array}
\label{eq:exo_2}
\end{equation}

Problems~\eqref{eq:exo_1} and~\eqref{eq:exo_2} are equivalent, see e.g.,~\cite{shapiro_interchangeability}. 
However, each of them has proved successful in different contexts. Problem~\eqref{eq:exo_1} has been the building block of most of the literature on the decision rule approximation, see Section~\ref{sec:introduction}. Problem~\eqref{eq:exo_2} has enabled the advent and tremendous success of the $K$-adaptability approximation approach to two-stage robust problems with binary recourse, see~\cite{Caramanis_FiniteAdaptability, Hanasusanto2015}. It has also facilitated the development of algorithms and efficient solution schemes, see e.g.,~\cite{zeng2013solving,ayoub2016decomposition}, and \cite{Shtern2018}.


\section{Two-Stage RO with Decision-Dependent Information Discovery}
\label{sec:problem_formulation_endo}

In this section, we describe two-stage robust optimization problems with decision-dependent information discovery (DDID) and propose an entirely new modeling framework for studying such problems. This framework underpins our ability to generalize the popular $K$-adaptability approximation approach from the literature to problems affected by uncertain parameters whose time of revelation is decision-dependent, see Sections~\ref{sec:kadaptability_objective} and~\ref{sec:kadaptability_constraint}.



\subsection{Problem Description}
\label{sec:exogenous_description}

In two-stage robust optimization with DDID, the uncertain parameters~${\bm \xi}$ do not necessarily become observed (for free) between the first and second decision-stages. Instead, some (typically costly) first stage decisions control the \emph{time of information discovery} in the problem: they decide whether (and which of) the uncertain parameters will be revealed \emph{before} the wait-and-see decisions ${\bm y}$ are selected. If the decision-maker chooses to not observe some of the uncertain parameters, then those parameters will still be uncertain at the time when the decision~${\bm y}$ is selected, and~${\bm y}$ will only be allowed to depend on the portion of the uncertain parameters that have been revealed. On the other hand, if the decision-maker chooses to observe all of the uncertain parameters, then there will be no uncertainty in the problem at the time when~${\bm y}$ is selected, and~${\bm y}$ will be allowed to depend on all uncertain parameters. 

In order to allow for endogenous uncertainty, we introduce a here-and-now binary measurement (or observation) decision vector ${\bm w} \in \{0,1\}^{N_\xi}$ of the same dimension as ${\bm \xi}$ whose $i$th element ${\bm w}_i$ is 1 if and only if we choose to observe ${\bm \xi}_i$ between the first and second decision stages. In the presence of such endogenous uncertain parameters, the recourse decisions ${\bm y}$ are selected after the \emph{portion} of uncertain parameters that was \emph{chosen} to be observed is revealed. In particular, ${\bm y}$ must be constant in (i.e., robust to) those uncertain parameters that remain unobserved at the second decision-stage. The requirement that ${\bm y}$ only depend on the uncertain parameters that have been revealed at the time it is chosen is termed \emph{non-anticipativity}. In the presence of uncertain parameters whose time of revelation is decision-dependent, this requirement translates to \emph{decision-dependent non-anticipativity constraints}.


\subsection{Decision Rule Formulation}
\label{sec:exogenous_literature}

In the literature and to the best of our knowledge, two-stage robust optimization problems with DDID have been formulated (in a manner paralleling Problem~\eqref{eq:exo_1}) by letting the recourse decisions~${\bm y}$ be functions of~${\bm \xi}$ and requiring that those functions be constant in~${\bm \xi}_i$ if ${\bm w}_i=0$, see~\cite{DDI_VKB}. Under this (decision rule based) modeling paradigm, generic two-stage robust optimization problems with decision-dependent information discovery take the form 
\begin{equation}\renewcommand{\arraystretch}{1.5}
    \begin{array}{cl}
         \minimize & \quad \displaystyle \max_{{\bm \xi} \in \Xi} \;\;\;  {\bm \xi}^\top {\bm C} \; {\bm x} + {\bm \xi}^\top {\bm D} \; {\bm w}  + {\bm \xi}^\top {\bm Q} \; {\bm y}({\bm \xi}) \\
         \subjectto & \quad {\bm x} \in \sets X, \; {\bm w} \in \sets W, \; {\bm y} \in \mathcal L_{N_\xi}^{N_y} \\
         & \quad \!\! \left. \begin{array}{l} 
         {\bm y}({\bm \xi}) \in \sets Y  \\
         {\bm T} {\bm x} + {\bm V} {\bm w} + {\bm W}{\bm y}({\bm \xi}) \leq {\bm H}{\bm \xi} 
         \end{array} \quad \right\} \quad \forall {\bm \xi} \in \Xi \\
         & \quad {\bm y}({\bm \xi}) = {\bm y}({\bm \xi}') \quad \forall {\bm \xi}, \; {\bm \xi}' \in \Xi \; : \; {\bm w} \circ {\bm \xi} = {\bm w} \circ {\bm \xi}',
    \end{array}
\label{eq:endo_1}
\end{equation}
where $\sets W \subseteq \{0,1\}^{N_\xi}$, ${\bm D} \in \reals^{N_\xi \times N_\xi}$, ${\bm V} \in \reals^{L \times N_\xi}$, and the remaining data elements are as in Problem~\eqref{eq:exo_1}. The set $\sets W$ can encode requirements on the measurement decisions. For example, it can enforce that a given uncertain parameter~${\bm \xi}_i$ may only be observed if another uncertain parameter~${\bm \xi}_{i'}$ has been observed using ${\bm w}_i \leq {\bm w}_{i'}$. Accordingly, it can postulate that the total number of uncertain parameters that are observed does not exceed a certain budget~$Q$ using $\sum_{i=1}^{N_\xi} {\bm w}_i \leq Q$. \newpv{If only some (or all) of the uncertain parameters have a time of information discovery that is exogenous, our models and solution approaches can be used by restricting the observation decisions ${\bm w}_i$ to equal~1 (resp.~0) for each \emph{exogenous} uncertain parameter~$i$ that is (resp.\ is not) observed between the first and second decision stages. These restrictions can be conveniently added as constraints to the set~$\mathcal W$.} The last constraint in the problem is a decision-dependent non-anticipativity constraint: it ensures that the function ${\bm y}$ is constant in the uncertain parameters that remain unobserved at the second stage. Indeed, the identity ${\bm w} \circ {\bm \xi} = {\bm w} \circ {\bm \xi}'$ evaluates to true only if the elements of ${\bm \xi}$ and ${\bm \xi}'$ that were observed are indistinguishable, in which case the decisions taken in scenarios ${\bm \xi}$ and ${\bm \xi}'$ must be equal. We omit joint (first stage) constraints on ${\bm x}$ and ${\bm w}$ to minimize notational overhead but emphasize that our approach remains applicable in their presence.

Note that Problem~\eqref{eq:endo_1} generalizes Problem~\eqref{eq:exo_1}. Indeed, if we set ${\bm w}=\1$, ${\bm D}={\bm 0}$, and ${\bm V}={\bm 0}$ in Problem~\eqref{eq:endo_1}, we recover Problem~\eqref{eq:exo_1}. In addition, it generalizes the single-stage robust problem: if we set ${\bm w}={\bm 0}$ in Problem~\eqref{eq:endo_1}, all uncertain parameters are revealed \emph{after} the second stage so that the second stage decisions are forced to be static (i.e., constant in ${\bm \xi}$).

To the best of our knowledge, the only approach in the literature for (approximately) solving problems of type~\eqref{eq:endo_1} is presented in~\cite{DDI_VKB} and relies on a decision rule approximation. The authors propose to approximate the binary (resp.\ continuous) wait-and-see decisions by functions that are piecewise constant (resp.\ piecewise linear) on a pre-selected partition of the uncertainty set of the form 
$
\Xi_{\bm s} \; := \; \left\{ {\bm \xi}\in \Xi \; : \; {\bm c}_{{\bm s}_i-1}^i \leq {\bm \xi}_i < {\bm c}_{{\bm s}_i}^i, \; i=1,\ldots,k \right\},
$ 
where ${\bm s} \in \sets S := \times_{i=1}^{N_\xi} \{1,\ldots,{\bm r}_i\} \subseteq \integers^{N_\xi}$ and $
{\bm c}_1^i \; < \; {\bm c}_2^i \; < \; \cdots \; < \; {\bm c}_{{\bm r}_i-1}^i \quad \text{for } i=1,\ldots,N_\xi
$ 
represent ${\bm r}_i-1$ breakpoints along the ${\bm \xi}_i$ axis. Unfortunately, as the following example illustrates, this approach is highly sensitive to the choice of breakpoint configuration.

\begin{example}
Consider the following instance of Problem~\eqref{eq:endo_1}
\begin{equation}\renewcommand{\arraystretch}{1.5}
    \begin{array}{cl}
         \minimize & \quad 0 \\
         \subjectto & \quad {\bm w} \in \{0,1\}^2, \; {\bm y} \in \mathcal B_2^2 \\
         & \quad \!\! \left. \begin{array}{l} 
         {\bm \xi} - {\bm \epsilon} \; \leq \; {\bm y}({\bm \xi}) \; \leq \; \1 + {\bm \xi} - {\bm \epsilon} \\
         \end{array} \;\; \right\} \;\; \forall {\bm \xi} \in \Xi \\
         & \quad {\bm y}({\bm \xi}) = {\bm y}({\bm \xi}') \quad \forall {\bm \xi}, \; {\bm \xi}' \in \Xi \; : \; {\bm w} \circ {\bm \xi} = {\bm w} \circ {\bm \xi}',
    \end{array}
    \label{eq:vayanos_DDID_conservative}
\end{equation}
where $\Xi := [-1,1]^2$. The inequality constraints in the problem combined with the requirement that ${\bm y}({\bm \xi})$ be binary imply that we must have ${\bm y}_i({\bm \xi}) = 1$ (resp.\ 0) whenever ${\bm \xi}_i > {\bm \epsilon}_i$ (resp.\ ${\bm \xi}_i < {\bm \epsilon}_i$). Thus, from the decision-dependent non-anticipativity constraints, the only feasible choice for ${\bm w}$ is $\1$. It is easy to show that if ${\bm \epsilon}=1{\rm e}{-3} \1$ and if we uniformly partition each axis iteratively in 2, 3, 4, etc.\ subsets, then 1999 breakpoints along each direction will need to be introduced before reaching a feasible (and thus optimal) solution. The associated problem will involve over $8{\rm e}7$ binary decision variables and $16{\rm e}7$ constraints. In contrast, as will become clear later on, our proposed solution approach with approximation parameter $K=4$ will be optimal in this case.\Halmos
\label{ex:vayanos_DDID_conservative}
\end{example}

Example~\ref{ex:vayanos_DDID_conservative} is not surprising: the approach from~\cite{DDI_VKB} was motivated by stochastic programs which are less sensitive to the breakpoint configuration than robust problems. Thus, a more flexible approach is needed to address two-stage and multi-stage robust problems with DDID.


\subsection{Proposed Min-Max-Min-Max Formulation}
\label{sec:min_max_min_max_formulation}

\begin{figure}[t!]
    \centering
    \includegraphics[width=0.45\textwidth]{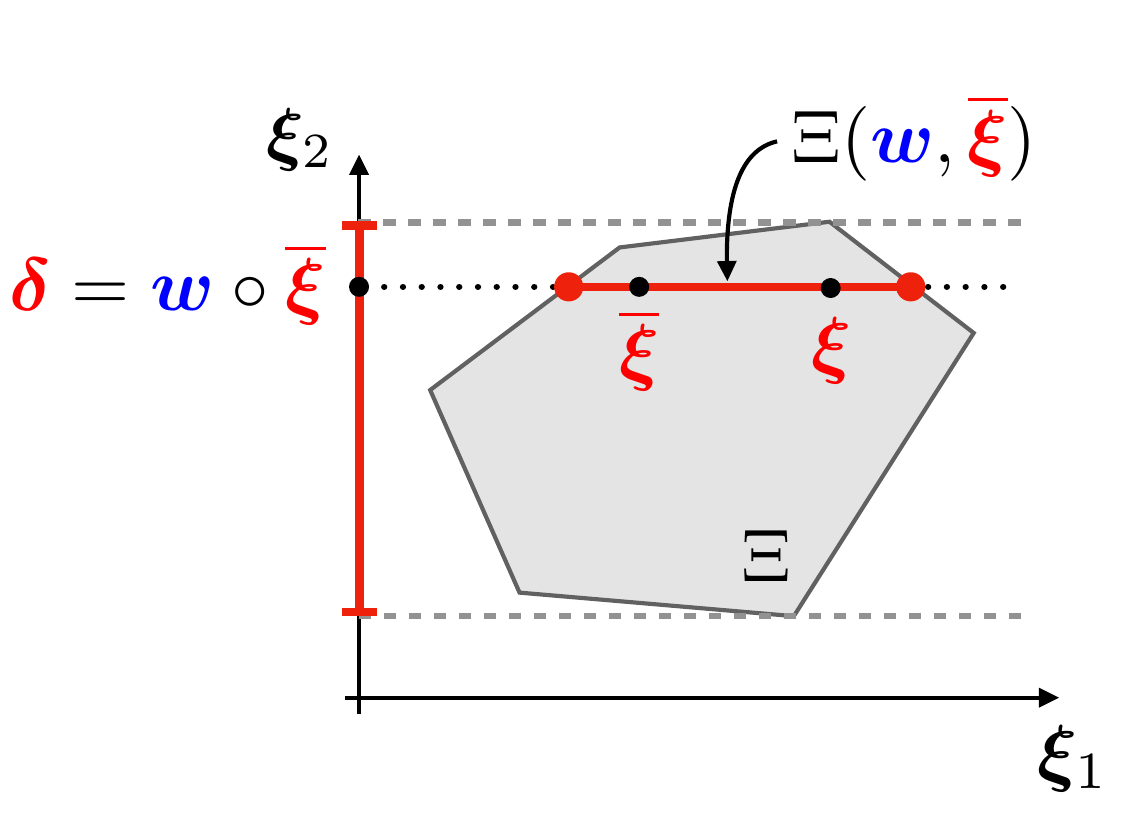} \vspace{0.4cm}
    \includegraphics[width=0.45\textwidth]{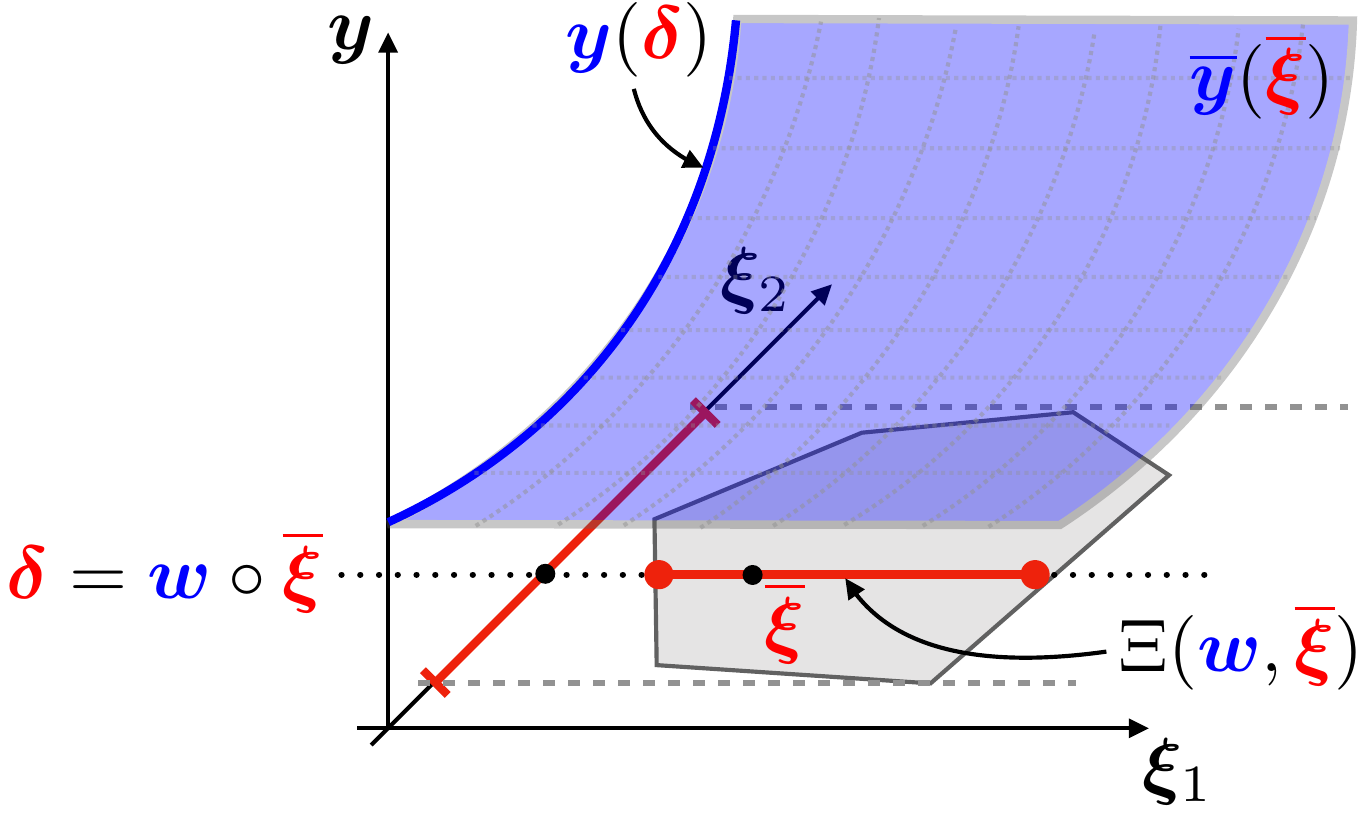}
    \caption{The figure on the left illustrates the role played by $\overline{\bm \xi}$ in the new formulation~\eqref{eq:endo_3} and the definition of the uncertainty sets $\Xi$ and $\Xi({\bm w},\overline{\bm \xi})$. Consider a setting where $\Xi \subseteq \reals^2$ (i.e., $N_\xi=2$) and suppose that ${\bm w} = (0,1)$ so that the decision-maker has chosen to only observe ${\bm \xi}_2$. In the figures, $\Xi$ is shown as the grey shaded area. Once $\overline{\bm \xi}$ is chosen by nature, the decision-maker can only infer that ${\bm \xi}$ will materialize in the set $\Xi({\bm w},\overline{\bm \xi})$ which collects all parameter realizations ${\bm \xi} \in \Xi$ that satisfy ${\bm \xi}_2 = \overline {\bm \xi}_2$, being compatible with our partial observation. The figure on the right illustrates the construction of an optimal non-anticipative decision $\overline {\bm y}$ from the an optimal solution ${\bm y}({\bm \delta})$ to $\displaystyle \min_{ {\bm y} \in \sets Y } \; \; \left\{ \max_{ {\bm \xi} \in \Xi({\bm w},{\bm \delta}) } \; \; {\bm \xi}^\top {\bm C} \; {\bm x} + {\bm \xi}^\top {\bm D} \; {\bm w} + {\bm \xi}^\top {\bm Q} \; {\bm y} \; : \; {\bm T} {\bm x} + {\bm V} {\bm w} + {\bm W}{\bm y} \leq {\bm H}{\bm \xi}  \; \; \; \; \forall {\bm \xi} \in \Xi({\bm w},{\bm \delta}) \right\}$, see Theorem~\ref{thm:endo_equiv}. We note that the policy $\overline{\bm y}$ constructed as in Theorem~\ref{thm:endo_equiv} is constant along the ${\bm \xi}_1$ direction since here ${\bm w}_1 = 0$. 
    }
    \label{fig:exo_2}
\end{figure}

Motivated by the success of formulation~\eqref{eq:exo_2} as the starting point to solve two-stage robust optimization problems with exogenous uncertainty, we derive an analogous \emph{dynamic} formulation for the case of endogenous uncertainties. In particular, we build a robust optimization problem in which the sequence of problems solved by each of the decision-maker and nature in turn is captured explicitly. The idea is as follows. Initially, the decision-maker selects ${\bm x} \in \sets X$ and ${\bm w} \in \sets W$. Subsequently, nature commits to a realization $\overline{\bm \xi}$ of the uncertain parameters from the set~$\Xi$. Then, the decision-maker selects a recourse action ${\bm y}$ that needs to be robust to those elements $\overline{\bm \xi}_i$ of the uncertain vector $\overline{\bm \xi}$ that they have not observed, i.e., for which ${\bm w}_i=0$. Indeed, the decision ${\bm y}$ may have to be taken under uncertainty if there is some $i$ such that ${\bm w}_i=0$, in which case not all of the uncertain parameters have been revealed when ${\bm y}$ is selected. Indeed, after ${\bm y}$ is selected, nature is free to choose any realization of ${\bm \xi} \in \Xi$ that is compatible with the original choice $\overline{\bm \xi}$ in the sense that ${\bm \xi}_i=\overline{\bm \xi}_i$ for all $i$ such that ${\bm w}_i=1$. This model captures the notion that, after ${\bm y}$ has been selected, nature is still free to choose the elements ${\bm \xi}_i$ that have not been observed, provided it does so in a way that is consistent with those parameters that \emph{have} been observed. Mathematically, given the measurement decisions ${\bm w}$ and the observation $\overline{\bm \xi}$, nature can select any element ${\bm \xi}$ from the set
$$
\Xi({\bm w},\overline{\bm \xi}) := \left\{ {\bm \xi} \in \Xi \; : \; {\bm w} \circ {\bm \xi} = {\bm w} \circ \overline{\bm \xi} \right\}.
$$
Note in particular that if ${\bm w}=\1$, then $\Xi({\bm w},\overline{\bm \xi}) = \{\overline{\bm \xi}\}$ and there is no uncertainty when ${\bm y}$ is chosen. Accordingly, if ${\bm w} ={\bm 0}$, then $\Xi({\bm w},\overline{\bm \xi}) = \Xi$ and ${\bm y}$ has no knowledge of any of the elements of ${\bm \xi}$. The realizations $\overline {\bm \xi}$, ${\bm \xi}$, and the sets $\Xi$ and $\Xi({\bm w},\overline{\bm \xi})$ are all illustrated on Figure~\ref{fig:exo_2}. 

Based on the above notation, we propose the following generic formulation of a two-stage robust optimization problem with decision-dependent information discovery:
\begin{equation}\renewcommand{\arraystretch}{1.5}
\tag{$\mathcal{P}$}
    \begin{array}{cl}
         \min & \;\; \displaystyle \max_{\overline{\bm \xi} \in \Xi} \;\;\min_{ {\bm y} \in \sets Y } \;  \left\{ \max_{ {\bm \xi} \in \Xi({\bm w},\overline{\bm \xi}) } \; \; {\bm \xi}^\top {\bm C} \; {\bm x} + {\bm \xi}^\top {\bm D} \; {\bm w} + {\bm \xi}^\top {\bm Q} \; {\bm y} \; : \; {\bm T} {\bm x} + {\bm V} {\bm w} + {\bm W}{\bm y} \leq {\bm H}{\bm \xi}  \; \; \; \; \forall {\bm \xi} \in \Xi({\bm w},\overline{\bm \xi}) \right\}  \\
         \st & \;\; {\bm x} \in \sets X, \; {\bm w} \in \sets W.
    \end{array}
\label{eq:endo_3}
\end{equation}
Note that, at the time when ${\bm y}$ is selected, some elements of ${\bm \xi}$ are still uncertain. The choice of ${\bm y}$ thus needs to be robust to the choice of those uncertain parameters that remain to be revealed. In particular, the constraints need to be satisfied for all choices of ${\bm \xi} \in \Xi({\bm w},\overline{\bm \xi})$. Accordingly, ${\bm y}$ is chosen so as to minimize the worst-case possible cost when ${\bm \xi}$ is valued in the set ${\bm \xi} \in \Xi({\bm w},\overline{\bm \xi})$.

Problems~\eqref{eq:endo_1} and~\eqref{eq:endo_3} are equivalent in a sense made precise in the following theorem.
%
\begin{theorem}
The optimal objective values of Problems~\eqref{eq:endo_1} and~\eqref{eq:endo_3} are equal.  Moreover, the following statements hold true:
\begin{enumerate}[label=(\roman*)]
    \item Let $({\bm x},{\bm w})$ be optimal in~\eqref{eq:endo_3} and, for each ${\bm \delta}$ such that ${\bm \delta} = {\bm w} \circ \overline{\bm \xi}$ for some $\overline{\bm \xi} \in \Xi$, define
$$
 {\bm y}'( {\bm \delta} ) \; \in \; \argmin_{ {\bm y} \in \sets Y } \; \; \left\{ \max_{ {\bm \xi} \in \Xi({\bm w},{\bm \delta}) } \; \; {\bm \xi}^\top {\bm C} \; {\bm x} + {\bm \xi}^\top {\bm D} \; {\bm w} + {\bm \xi}^\top {\bm Q} \; {\bm y} \; : \; {\bm T} {\bm x} + {\bm V} {\bm w} + {\bm W}{\bm y} \leq {\bm H}{\bm \xi}  \; \; \; \; \forall {\bm \xi} \in \Xi({\bm w},{\bm \delta}) \right\}.
$$
Also, for each ${\bm \xi}\in \Xi$, define ${\bm y}({\bm \xi}) := {\bm y}'({\bm w} \circ {\bm \xi})$. Then, $({\bm x}, {\bm w}, {\bm y}(\cdot))$ is optimal in Problem~\eqref{eq:endo_1}. 
    \item Let $( {\bm x}, {\bm w},{\bm y}(\cdot))$ be optimal in Problem~\eqref{eq:endo_1}. Then, $( {\bm x}, {\bm w})$ is optimal in Problem~\eqref{eq:endo_3}.
\end{enumerate}
\label{thm:endo_equiv}
\end{theorem}

The parameter ${\bm \delta}$ in item \textit{(i)} of the theorem above is introduced to ensure that the decision rule ${\bm y}(\cdot)$ defined on $\Xi$ is non-anticipative. Indeed, if for any given $({\bm x},{\bm w})$ and $\overline{\bm \xi}$, there are many optimal solutions to problem
$$
\min_{ {\bm y} \in \sets Y } \;  \left\{ \max_{ {\bm \xi} \in \Xi({\bm w},\overline{\bm \xi}) } \; \; {\bm \xi}^\top {\bm C} \; {\bm x} + {\bm \xi}^\top {\bm D} \; {\bm w} + {\bm \xi}^\top {\bm Q} \; {\bm y} \; : \; {\bm T} {\bm x} + {\bm V} {\bm w} + {\bm W}{\bm y} \leq {\bm H}{\bm \xi}  \; \; \; \; \forall {\bm \xi} \in \Xi({\bm w},\overline{\bm \xi}) \right\},
$$
the decision rule $\tilde{\bm y}(\cdot)$ defined on $\Xi$ through 
$$
 \tilde{\bm y}( \overline{\bm \xi} ) \; \in \; \argmin_{ {\bm y} \in \sets Y } \; \; \left\{ \max_{ {\bm \xi} \in \Xi({\bm w},\overline{\bm \xi}) } \; \; {\bm \xi}^\top {\bm C} \; {\bm x} + {\bm \xi}^\top {\bm D} \; {\bm w} + {\bm \xi}^\top {\bm Q} \; {\bm y} \; : \; {\bm T} {\bm x} + {\bm V} {\bm w} + {\bm W}{\bm y} \leq {\bm H}{\bm \xi}  \; \; \; \; \forall {\bm \xi} \in \Xi({\bm w},\overline{\bm \xi}) \right\},
$$
may not be constant in those parameters that remain unobserved. We note of course that other tie-breaking mechanisms could be used to build a non-anticipative solution. For example, we may select, among all optimal solutions, the one that is lexicographically first. 

The theorem above is the main result that enables us to generalize the $K$-adaptability approximation scheme to two-stage robust problems with decision-dependent information discovery and binary recourse. In Electronic Companion~\ref{sec:EC_optval_notattained}, we show that for any given choice of here-and-now decisions, the set of parameters ${\bm \xi}$ for which a particular wait-and-see decision is optimal may be non closed and non-convex and that the optimal value of the problem may not be attained. This result is expected from the analysis in~\cite{Hanasusanto2015}, since Problem~\eqref{eq:endo_3} generalizes Problem~\eqref{eq:exo_2}. Our example illustrates that this may be the case even if a portion of the uncertain parameters remain unobserved in the second stage.

Two-stage robust optimization problems with decision-dependent information discovery have a huge modeling power, see Sections~\ref{sec:introduction}, \ref{sec:computational_studies}, and~\ref{sec:preference_elicitation_KAS}. Yet, as illustrated by the preceding discussion, they pose several theoretical and practical challenges. As we will see in the following sections, whether we are or not able to reformulate the the $K$-adaptability counterpart of the problem exactly as a finite program solvable with off-the-shelf solvers depends on the absence or presence of uncertainty in the constraints. When in presence of constraint uncertainty, we can always compute an arbitrarily tight outer (lower bound) approximation, see Section~\ref{sec:kadaptability_constraint}.


\section{$K$-Adaptability for Problems with DDID}
\label{sec:kadaptability_overview}

Instead of solving Problem~\eqref{eq:endo_3} directly, we approximate it through its $K$-adaptability counterpart,
\begin{equation}\renewcommand{\arraystretch}{1.5}
\tag{$\mathcal{P}_K$}
    \begin{array}{cl}
         \min & \;\; \displaystyle \max_{\overline {\bm \xi} \in \Xi} \;\;\min_{ k \in \sets K } \;  \left\{ \max_{ {\bm \xi} \in \Xi({\bm w},\overline{\bm \xi}) } \; \; {\bm \xi}^\top {\bm C} \; {\bm x} + {\bm \xi}^\top {\bm D} \; {\bm w} + {\bm \xi}^\top {\bm Q} \; {\bm y}^k \; : \; {\bm T} {\bm x} + {\bm V} {\bm w} + {\bm W}{\bm y}^k \leq {\bm H}{\bm \xi}  \; \; \; \; \forall {\bm \xi} \in \Xi({\bm w},\overline{\bm \xi}) \right\}  \\
         \st & \;\; {\bm x} \in \sets X, \; {\bm w} \in \sets W, \; {\bm y}^k \in \sets Y, \; k \in \sets K,
    \end{array}
\label{eq:Kadapt}
\end{equation}
where $\sets K:= \{1,\ldots,K\}$. In this problem, $K$ candidate policies ${\bm y}^1,\ldots,{\bm y}^K$ are chosen here-and-now, that is before ${\bm w} \circ \overline{\bm \xi}$ (the portion of uncertain parameters that we chose to observe) is revealed. Once~${\bm w} \circ \overline{\bm \xi}$ becomes known, the best of those policies among all those that are robustly feasible (in view of uncertainty in the uncertain parameters that are still unknown) is implemented. If all policies are infeasible for some $\overline{\bm \xi} \in \Xi$, then we interpret the maximum and minimum in \eqref{eq:Kadapt} as supremum and infimum, that is, the $K$-adaptability problem evaluates to $+\infty$. Problem \eqref{eq:Kadapt} is a conservative approximation to program~\eqref{eq:endo_3}. Moreover, if $| \sets Y | < \infty$ and $K=| \sets Y |$, then the two problems are equivalent. In practice, we hope that a moderate number of candidate policies $K$ will be sufficient to obtain a (near) optimal solution to~\eqref{eq:endo_3}.

\paragraph{The Price of Usability.} We note that Problem~\eqref{eq:Kadapt} is interesting in its own right. Indeed, in problems where usability is important (e.g., if workers need to be trained to follow diverse contingency plans depending on the realization ${\bm w} \circ \overline {\bm \xi}$), Problem~\eqref{eq:Kadapt} may be an attractive alternative to Problem~\eqref{eq:endo_3}. In such settings, the loss in optimality incurred due to passing from Problem~\eqref{eq:endo_3} to Problem~\eqref{eq:Kadapt} can be thought of as the \emph{price of usability}. For example, consider an emergency response planning problem where, in the first stage, a small number of helicopters can be used to survey affected areas and, in the second stage, and in response to the observed state of the areas surveyed, deployment of emergency response teams is decided. In practice, to avoid having to train teams in a large number of plans (yielding significant operational challenges), only a moderate number of response plans may be allowed. The importance of interpretability/usability has been previously noted by e.g.,~\cite{Koc2015,McCarthyijcai2018,bertsimas2019price}, and \cite{Aghaei2019AAAI,Aghaei_StrongCT}.

\begin{remark}
If $\mathcal W=\{0,1\}^{N_\xi}$, ${\bm D}={\bm 0}$, and ${\bm V}={\bm 0}$, then ${\bm w}={\textbf{e}}$ is optimal in Problem~$(\mathcal P)$ and thus $\Xi({\bm w},\overline{\bm \xi})=\{ \overline{\bm \xi} \}$, implying that Problem~$(\mathcal P)$ reduces to Problem~(2) and Problem~\eqref{eq:Kadapt} reduces to the $K$-adaptability counterpart of Problem~\eqref{eq:exo_2}.
\label{rmk:Kadapt_exo}
\end{remark}


Relative to the problems studied by~\cite{Caramanis_FiniteAdaptability} and~\cite{Hanasusanto2015}, Problem~\eqref{eq:Kadapt} presents several challenges. First, the second stage problem in~\eqref{eq:Kadapt} is a robust (as opposed to deterministic) optimization problem. Second, the uncertainty sets involved in the maximization tasks of this robust problem are decision-dependent. While Problem~\eqref{eq:Kadapt} appears to be significantly more complicated than its exogenous counterpart, it can be converted to an equivalent min-max-min problem by \emph{lifting} the space of the uncertainty set as show in the following lemma that is instrumental in our analysis. 

\begin{lemma}
The $K$-adaptability problem with decision-dependent information discovery, Problem~\eqref{eq:Kadapt}, is equivalent to
\begin{equation}\renewcommand{\arraystretch}{1.5}
    \begin{array}{cl}
         \min & \;\; \displaystyle \max_{ \{ {\bm \xi}^k \}_{k \in \sets K} \in \Xi^K({\bm w}) } \;\;\min_{ k \in \sets K } \;  \left\{  ({\bm \xi}^k)^\top {\bm C} \; {\bm x} + ({\bm \xi}^k)^\top {\bm D} \; {\bm w} + ({\bm \xi}^k)^\top {\bm Q} \; {\bm y}^k \; : \; {\bm T} {\bm x} + {\bm V} {\bm w} + {\bm W}{\bm y}^k \leq {\bm H}{\bm \xi}^k \right\}  \\
         \st & \;\; {\bm x} \in \sets X, \; {\bm w} \in \sets W, \; {\bm y}^k \in \sets Y, \; k \in \sets K,
    \end{array}
\label{eq:Kadapt_cstr_min_max_min}
\end{equation}
where 
\begin{equation}
 \Xi^K({\bm w}) := \left\{
\{ {\bm \xi}^k \}_{k \in \sets K} \in \Xi^K \; : \; \exists \overline {\bm \xi} \in \Xi \text{ such that } {\bm \xi}^k \in \Xi( {\bm w},\overline {\bm \xi}) \text{ for all } k \in \sets K
\right\}.
\label{eq:Xi_K}
\end{equation}
\label{lem:Kadapt_cstr_min_max_min}
\end{lemma}
For any fixed ${\bm w} \in \sets W$, the subvector ${\bm \xi}^k$ in the definition of ${\bm \Xi}^K({\bm w})$ represents the uncertainty scenario that ``nature'' will choose if the decision-maker acts according to decisions ${\bm w}$ in the first stage and according to policy $k$ in the second stage. The set ${\bm \Xi}^K({\bm w})$ collects, for each $k\in \sets K$, all feasible choices that nature can take if the decision-maker acts according to ${\bm w}$ and then ${\bm y}^k$ in the first and second stages, respectively. Thus, in Problem~\eqref{eq:Kadapt_cstr_min_max_min}, the decision-maker first selects ${\bm x}$, ${\bm w}$, and ${\bm y}^k$, $k\in \sets K$. Subsequently, nature commits to the portion of observed uncertain parameters ${\bm w}\circ {\overline{\bm \xi}}$ and to a choice ${\bm \xi}^k$, $k\in \sets K$, associated with each candidate policy ${\bm y}^k$. Finally, the decision-maker chooses one of the candidate policies.
%

In what follows, we provide insights into the theoretical and computational properties of the $K$-adaptability counterpart to two-stage robust problems with DDID and with binary recourse. 

\begin{remark} We note that the results in Section~\ref{sec:problem_formulation_endo} generalize fully to cases where the objective and constraint functions are continuous (not necessarily linear) in ${\bm x}$, ${\bm y}$, and ${\bm \xi}$. Moreover, 
all of the ideas in our paper generalize to the case where the technology and recourse matrices, ${\bm T}$ and ${\bm W}$, depend on ${\bm \xi}$. We do not discuss these cases in detail so as to minimize notational overhead.
\end{remark}


\subsection{$K$-Adaptability for Problems with Objective Uncertainty}
\label{sec:kadaptability_objective}



In this section, we focus our attention on the case where uncertain parameters only appear in the objective of Problem~\eqref{eq:endo_3} and where the recourse decisions are binary, being expressible as
\begin{equation}\renewcommand{\arraystretch}{1.5}
\tag{$\mathcal{PO}$}
    \begin{array}{cl}
         \minimize & \;\; \displaystyle \max_{\overline{\bm \xi} \in \Xi} \;\;\min_{ {\bm y} \in \sets Y } \;  \left\{ \max_{ {\bm \xi} \in \Xi({\bm w},\overline{\bm \xi}) } \; \; {\bm \xi}^\top {\bm C} \; {\bm x} + {\bm \xi}^\top {\bm D} \; {\bm w} + {\bm \xi}^\top {\bm Q} \; {\bm y} \; : \; {\bm T} {\bm x} + {\bm V} {\bm w} + {\bm W}{\bm y} \leq {\bm h}  \right\}  \\
         \subjectto & \;\; {\bm x} \in \sets X, \; {\bm w} \in \sets W, 
    \end{array}
\label{eq:endo_2_obj}
\end{equation}
where ${\bm h} \in \reals^L$, $\sets Y \subseteq \{0,1\}^{N_y}$. We study the $K$-adaptability counterpart of Problem~\eqref{eq:endo_2_obj} given by
\begin{equation}\renewcommand{\arraystretch}{1.5}
\tag{$\mathcal{PO}_K$}
    \begin{array}{cl}
         \minimize & \;\; \displaystyle \max_{\overline{\bm \xi} \in \Xi} \;\;\min_{ k \in \sets K } \;  \left\{ \max_{ {\bm \xi} \in \Xi({\bm w},\overline{\bm \xi}) } \; \; {\bm \xi}^\top {\bm C} \; {\bm x} + {\bm \xi}^\top {\bm D} \; {\bm w} + {\bm \xi}^\top {\bm Q} \; {\bm y}^k \; : \; {\bm T} {\bm x} + {\bm V} {\bm w} + {\bm W}{\bm y}^k \leq {\bm h} \right\}  \\
         \subjectto & \;\; {\bm x} \in \sets X, \; {\bm w} \in \sets W, \; {\bm y}^k \in \sets Y, \; k \in \sets K.
    \end{array}
\label{eq:Kadapt_obj}
\end{equation}
Applying Lemma~\ref{lem:Kadapt_cstr_min_max_min}, we are able to write Problem~\eqref{eq:Kadapt_obj} equivalently as
\begin{equation}\renewcommand{\arraystretch}{1.5}
    \begin{array}{cl}
         \minimize & \quad \displaystyle \max_{ \{ {\bm \xi}^k \}_{k \in \sets K} \in \Xi^K({\bm w}) } \;\;\min_{ k \in \sets K } \;\;  \left\{  ({\bm \xi}^k)^\top {\bm C} \; {\bm x} + ({\bm \xi}^k)^\top {\bm D} \; {\bm w} + ({\bm \xi}^k)^\top {\bm Q} \; {\bm y}^k \; : \; {\bm T} {\bm x} + {\bm V} {\bm w} + {\bm W}{\bm y}^k \leq {\bm h} \right\}  \\
         \subjectto & \quad {\bm x} \in \sets X, \; {\bm w} \in \sets W, \; {\bm y}^k \in \sets Y, \; k \in \sets K ,
    \end{array}
\label{eq:Kadapt_min_max_min_2}
\end{equation}
where $\Xi^K({\bm w})$ is defined as in Lemma~\ref{lem:Kadapt_cstr_min_max_min}. In the absence of uncertainty in the constraints, the constraints in the $K$-adaptability problem can be moved to the first stage, as summarized by the following observation.

\begin{observation}
The $K$-adaptability counterpart of the two-stage robust optimization problem with decision-dependent information discovery, Problem~\eqref{eq:Kadapt_obj}, is equivalent to
\begin{equation}\renewcommand{\arraystretch}{1.5}
    \begin{array}{cl}
         \minimize & \quad \displaystyle \max_{ \{ {\bm \xi}^k \}_{k \in \sets K} \in \Xi^K({\bm w}) } \;\;\min_{ k \in \sets K } \;\;  \left\{  ({\bm \xi}^k)^\top {\bm C} \; {\bm x} + ({\bm \xi}^k)^\top {\bm D} \; {\bm w} + ({\bm \xi}^k)^\top {\bm Q} \; {\bm y}^k  \right\}  \\
         \subjectto & \quad {\bm x} \in \sets X, \; {\bm w} \in \sets W, \; {\bm y}^k \in \sets Y, \; k \in \sets K \\
         & \quad {\bm T} {\bm x} + {\bm V} {\bm w} + {\bm W}{\bm y}^k \leq {\bm h} \quad \forall k \in \sets K,
    \end{array}
\label{eq:Kadapt_obj_move_constraints_2}
\end{equation}
where $\Xi^K({\bm w})$ is as defined in Equation~\eqref{eq:Xi_K}.
\label{obs:Kadapt_obj_move_constraints_2}
\end{observation}

Note that for all ${\bm w} \in \sets W$, the set $\Xi^K({\bm w})$ is non-empty and bounded. Thus, $({\bm x},{\bm w},\{{\bm y}^k\}_{k \in \sets K}) \in \sets X \times \sets W \times \sets Y^K$ is feasible in Problem~\eqref{eq:Kadapt_obj_move_constraints_2} if ${\bm T} {\bm x} + {\bm V} {\bm w} + {\bm W}{\bm y}^k \leq {\bm h}$ for all $k \in \sets K$, whereas to be feasible in Problem~\eqref{eq:Kadapt_min_max_min_2} (and accordingly in Problem~\eqref{eq:Kadapt_obj}), it need only satisfy ${\bm T} {\bm x} + {\bm V} {\bm w} + {\bm W}{\bm y}^k \leq {\bm h}$ for some $k \in \sets K$. Thus, a triplet $({\bm x},{\bm w},{\bm y}^k)$ feasible in~\eqref{eq:Kadapt_min_max_min_2} (and thus in~\eqref{eq:Kadapt_obj}) need not be feasible in Problem~\eqref{eq:Kadapt_obj_move_constraints_2}. However, the proof of Observation~\ref{obs:Kadapt_obj_move_constraints_2}, provides a way to construct a feasible solution for Problem~\eqref{eq:Kadapt_obj_move_constraints_2} from a feasible solution to Problem~\eqref{eq:Kadapt_min_max_min_2} that achieves the same optimal value. 

Lemma~\ref{lem:Kadapt_cstr_min_max_min} and Observation~\ref{obs:Kadapt_obj_move_constraints_2} are key to reformulating Problem~\eqref{eq:Kadapt_obj} as a finite program. They also enable us to analyze the complexity of evaluating the objective function of the $K$-adaptability problem under a fixed decision. Indeed, from Problem~\eqref{eq:Kadapt_obj_move_constraints_2}, it can be seen that for any fixed choice $({\bm x}, {\bm w}, \{ {\bm y}^k \}_{k \in \sets K})$, the objective value of~\eqref{eq:Kadapt_obj} can be evaluated by solving a linear program (LP) obtained by writing~\eqref{eq:Kadapt_obj_move_constraints_2} in epigraph form. We formalize this result in the following. 
\begin{observation}
For any fixed $K$ and decision $({\bm x}, {\bm w}, \{ {\bm y}^k \}_{k \in \sets K})$, the optimal objective value of the $K$-adaptability problem~\eqref{eq:Kadapt_obj} can be evaluated in polynomial time in the size of the input.
\label{obs:evaluation_Kadapt_obj}
\end{observation}



In Observation~\ref{obs:evaluation_Kadapt_obj}, we showed that for any fixed $K$, ${\bm x}$, ${\bm w}$, and ${\bm y}^k$, the objective function in Problem~\eqref{eq:Kadapt_obj} can be evaluated by means of a polynomially sized LP. By dualizing this LP , we can obtain an equivalent reformulation of Problem~\eqref{eq:Kadapt_obj} in the form of a bilinear problem.
\begin{theorem}
Problem~\eqref{eq:Kadapt_obj} is equivalent to the bilinear problem
\begin{equation}\renewcommand{\arraystretch}{1.5}
    \begin{array}{cl}
         \minimize & \quad   {\bm b}^\top {\bm \beta}  + \sum_{k\in \mathcal K} {\bm b}^\top {\bm \beta}^k\\
         \subjectto & \quad {\bm x} \in \sets X, \; {\bm w} \in \sets W, \; {\bm y}^k \in \sets Y, \; k \in \sets K \\
         & \quad {\bm \alpha} \in \mathbb R^K_+, \; {\bm \beta} \in \mathbb R^R_+, \; {\bm \beta}^k \in \mathbb R^R_+, \; {\bm \gamma}^k \in \mathbb R^{N_\xi}, \; k \in \mathcal K \\
        & \quad \1^\top {\bm \alpha}  = 1  \\
        & \quad     {\bm A}^\top{\bm \beta}^k +  {\bm w} \circ {\bm \gamma}^k = {\bm \alpha}_k \left( {\bm C} {\bm x} + {\bm D} {\bm w} + {\bm Q} {\bm y}^k \right) \quad \forall k \in \mathcal K \\
        & \quad  {\bm A}^\top{\bm \beta}  =  \displaystyle \sum_{k \in \mathcal K} {\bm w} \circ  {\bm \gamma}^k \\
         & \quad {\bm T} {\bm x} + {\bm V} {\bm w} + {\bm W}{\bm y}^k \leq {\bm h} \quad \forall k \in \sets K.
    \end{array}
\label{eq:Kadapt_obj_bilinear}
\end{equation}
\label{thm:Kadapt_obj}
\end{theorem}
Although Problem~\eqref{eq:Kadapt_obj} is generally non-convex (bilinear), there exist several techniques in the literature for solving such problems exactly. In fact, this is an extremely active area of research, see e.g., \cite{Tsoukalas2014} and~\cite{Gupte2017}. Moreover, problems of the form~\eqref{eq:Kadapt_obj} can now be solved with state-of-the-art off-the-shelf solvers like Gurobi. Indeed, Gurobi recently released its 9th version that can tackle non-convex quadratic programs.\footnote{See e.g., \url{https://www.gurobi.com/documentation/9.0/refman/nonconvex.html}} If $\sets X \subseteq \{0,1\}^{N_x}$ and $\sets Y \subseteq \{0,1\}^{N_y}$, the bilinear terms in the formulation above can be linearized using standard techniques and we can obtain an equivalent reformulation of Problem~\eqref{eq:Kadapt_obj} in the form of an MBLP.
\begin{corollary}
Suppose $\sets X \subseteq \{0,1\}^{N_x}$ and $\sets Y \subseteq \{0,1\}^{N_y}$. Then, Problem~\eqref{eq:Kadapt_obj} is equivalent an MBLP involving a suitably chosen ``big-$M$'' constant.
\label{cor:Kadapt_obj_MBLP}
\end{corollary}

We emphasize that the size of the MBLP in Corollary~\ref{cor:Kadapt_obj_MBLP} is polynomial in the size of the input data for the $K$-adaptability problem~\eqref{eq:Kadapt_obj}. Note that, contrary to~\cite{Hanasusanto2015}, we require that $\sets X \subseteq \{0,1\}^{N_x}$. This is to ensure that we are able to linearize the bilinear terms involving the ${\bm x}$ variables that arise from the dualization step. We note that formulation~\eqref{eq:Kadapt_obj_bilinear} and its equivalent MBLP can be augmented with symmetry breaking constraints to speed-up solution, see Section~\ref{sec:symmetry_breaking_general} for details.

\begin{remark}
Most MBLP solvers\footnote{See e.g., \url{https://www.ibm.com/analytics/cplex-optimizer} and \url{https://www.gurobi.com/}.} allow reformulating the bilinear terms without the use of ``big-$M$'' constants, which are known to suffer from numerical instability. These include, for example, so-called \texttt{SOS} or \texttt{IfThen} constraints.
\end{remark}

\begin{observation}
Suppose that we are only in the presence of exogenous uncertainty, i.e., ${\bm w}=\1$, ${\bm D}={\bm 0}$, and ${\bm V}={\bm 0}$. Then, Problem~\eqref{eq:Kadapt_MBLP} reduces to the MBLP formulation of the $K$-adaptability problem with only exogenous uncertainty from~\cite{Hanasusanto2015}. 
\label{obs:grani_equivalence_DVzero_objective}
\end{observation}


\subsection{$K$-Adaptability for Problems with Constraint Uncertainty}
\label{sec:kadaptability_constraint}

The starting point of our analysis is the reformulation of Problem~\eqref{eq:Kadapt} as the min-max-min problem~\eqref{eq:Kadapt_cstr_min_max_min}. Unfortunately, this problem is generally hard as testified by the following theorem.


\begin{theorem}
Evaluating the objective of Problem~\eqref{eq:Kadapt_cstr_min_max_min} if $K$ is not fixed is strongly NP-hard.
\label{thm:hardness_endo_2}
\end{theorem}

We reformulate Problem~\eqref{eq:Kadapt_cstr_min_max_min} equivalently by shifting the second-stage constraints ${\bm T} {\bm x} + {\bm V} {\bm w} + {\bm W}{\bm y}^k \leq {\bm H}{\bm \xi}^k$ from the objective function to the definition of the uncertainty set. We thus replace $\Xi^K({\bm w})$ with a family of uncertainty sets parameterized by a vector ${\bm \ell}$.

\begin{proposition}
The $K$-adaptability problem with decision-dependent information discovery, Problem~\eqref{eq:Kadapt_cstr_min_max_min}, is equivalent to
\begin{equation}\renewcommand{\arraystretch}{1.5}
    \begin{array}{cl}
         \minimize & \;\;  \displaystyle  \max_{{\bm \ell} \in \sets L} \;\; \max_{ \{ {\bm \xi}^k \}_{ k \in \sets K} \in \Xi^K({\bm w},{\bm \ell}) } \;\;\min_{ \begin{smallmatrix} k \in \sets K : \\ {\bm \ell}_k = 0 \end{smallmatrix} } \;  \left\{  ({\bm \xi}^k)^\top {\bm C} \; {\bm x} + ({\bm \xi}^k)^\top {\bm D} \; {\bm w} + ({\bm \xi}^k)^\top {\bm Q} \; {\bm y}^k \right\}  \\
         \subjectto & \;\; {\bm x} \in \sets X, \; {\bm w} \in \sets W, \; {\bm y}^k \in \sets Y, \; k \in \sets K,
    \end{array}
\label{eq:Kadapt_ell}
\end{equation}
where $\sets L := \{0,\ldots,L\}^K$, $L$ is the number of second-stage constraints in Problem~\eqref{eq:endo_3}, and the uncertainty sets $\Xi^K({\bm w},{\bm \ell})$, ${\bm \ell} \in \sets L$, are defined as
$$\renewcommand{\arraystretch}{1.5}
\Xi^K({\bm w},{\bm \ell}) := 
    \left\{ \{ {\bm \xi}^k \}_{k \in \sets K} \in \Xi^K \; : \;
    \begin{array}{ll}
         {\bm w} \circ {\bm \xi}^k = {\bm w} \circ \overline {\bm \xi} & \quad \forall k \in \sets K \text{ for some } \overline {\bm \xi} \in \Xi \\
         {\bm T} {\bm x} + {\bm V} {\bm w} + {\bm W}{\bm y}^k \leq {\bm H}{\bm \xi}^k & \quad \forall k \in \sets K : {\bm \ell}_k = 0 \\
         \left[{\bm T} {\bm x} + {\bm V} {\bm w} + {\bm W}{\bm y}^k \right]_{{\bm \ell}_k} > [{\bm H}{\bm \xi}^k]_{{\bm \ell}_k} & \quad \forall k \in \sets K : {\bm \ell}_k \neq 0
    \end{array}
    \right\},
$$
where, for convenience, we have suppressed the dependence of $\Xi^K({\bm w},{\bm \ell})$ on ${\bm x}$ and ${\bm y}^k$, $k \in \sets K$.
\label{prop:Kadapt_ell}
\end{proposition}

The elements of vector ${\bm \ell} \in \sets L$ in Proposition~\ref{prop:Kadapt_ell} encode which second-stage policies are feasible for the parameter realizations $\{ {\bm \xi}^k \}_{k \in \sets K} \in \Xi^K({\bm w},{\bm \ell})$. Indeed, recall that ${\bm \xi}^k$ can be viewed as the recourse action that nature will take if the decision-maker acts according to ${\bm y}^k$ in response to seeing $\overline{\bm \xi}$. Thus, policy ${\bm y}^k$ is feasible in Problem~\eqref{eq:Kadapt_cstr_min_max_min} (and thus in Problem~\eqref{eq:Kadapt}) if ${\bm \ell}_k = 0$. On the other hand, policy ${\bm y}^k$ violates the ${\bm \ell}_k$-th constraint in Problem~\eqref{eq:Kadapt_cstr_min_max_min} if ${\bm \ell}_k \neq 0$. Thus, if ${\bm \ell}_k \neq 0$, this implies that the ${\bm \ell}_k$-th constraint in~\eqref{eq:Kadapt} is violated for some ${\bm \xi} \in \Xi({\bm w},\overline{\bm \xi})$ and therefore ${\bm y}^k$ is not feasible in~\eqref{eq:Kadapt}. Note that, in contrast to the case with exogenous uncertainty discussed by~\cite{Hanasusanto2016}, ${\bm \ell}_k=0$ if and only if policy ${\bm y}^k$ is \emph{robustly} feasible in~\eqref{eq:Kadapt}.



Having brought Problem~\eqref{eq:Kadapt} to the form~\eqref{eq:Kadapt_ell}, it now presents a similar structure to a problem with objective uncertainty (see Section~\ref{sec:kadaptability_objective}) with the caveats that the problem involves multiple uncertainty sets that are also open. Next, we employ closed inner approximations $\Xi^K_\epsilon({\bm w},{\bm \ell})$ of the sets $\Xi^K({\bm w},{\bm \ell})$ that are parameterized by a scalar $\epsilon > 0$:
\begin{equation}\renewcommand{\arraystretch}{1.5}
\tag{$\text{\ref{eq:Kadapt_ell}}_\epsilon$}
    \begin{array}{cl}
         \minimize & \;\;  \displaystyle  \max_{{\bm \ell} \in \sets L} \;\; \max_{ \{ {\bm \xi}^k \}_{ k \in \sets K} \in \Xi^K_\epsilon({\bm w},{\bm \ell}) } \;\;\min_{ \begin{smallmatrix} k \in \sets K : \\ {\bm \ell}_k = 0 \end{smallmatrix} } \;  \left\{  ({\bm \xi}^k)^\top {\bm C} \; {\bm x} + ({\bm \xi}^k)^\top {\bm D} \; {\bm w} + ({\bm \xi}^k)^\top {\bm Q} \; {\bm y}^k \right\}  \\
         \subjectto & \;\; {\bm x} \in \sets X, \; {\bm w} \in \sets W, \; {\bm y}^k \in \sets Y, \; k \in \sets K,
    \end{array}
\label{eq:Kadapt_ell_eps}
\end{equation}
where the uncertainty sets $\Xi^K_\epsilon({\bm w},{\bm \ell})$ are defined as
$$\renewcommand{\arraystretch}{1.5}
\Xi^K_\epsilon({\bm w},{\bm \ell}) := 
    \left\{ \{ {\bm \xi}^k \}_{k \in \sets K} \in \Xi^K \; : \;
    \begin{array}{ll}
         {\bm w} \circ {\bm \xi}^k = {\bm w} \circ \overline {\bm \xi} & \quad \forall k \in \sets K \text{ for some } \overline {\bm \xi} \in \Xi \\
         {\bm T} {\bm x} + {\bm V} {\bm w} + {\bm W}{\bm y}^k \leq {\bm H}{\bm \xi}^k & \quad \forall k \in \sets K : {\bm \ell}_k = 0 \\
         \left[{\bm T} {\bm x} + {\bm V} {\bm w} + {\bm W}{\bm y}^k \right]_{{\bm \ell}_k} \geq [{\bm H}{\bm \xi}^k]_{{\bm \ell}_k} + \epsilon & \quad \forall k \in \sets K : {\bm \ell}_k \neq 0
    \end{array}
    \right\}.
$$
Using this definition, we next reformulate the approximate Problem~\eqref{eq:Kadapt_ell_eps} equivalently as an MBLP.

\begin{theorem}
The approximate problem~\eqref{eq:Kadapt_ell_eps} is equivalent to the mixed binary bilinear program 
\begin{equation}\renewcommand{\arraystretch}{1.5}
    \begin{array}{cll}
         \min & \;\;  \tau  \\
         \st & \;\; \tau \in \mathbb R, \; {\bm x} \in \sets X, \; {\bm w} \in \sets W, \; {\bm y}^k \in \sets Y, \; k \in \sets K \\
         & \;\;  {\bm \alpha}({\bm \ell}) \in \reals^R_+, \; {\bm \alpha}^k({\bm \ell}) \in \reals^R_+, \; k \in \sets K, \; {\bm \gamma}({\bm \ell}) \in \reals^K_+, \; {\bm \eta}^k({\bm \ell}) \in \reals^{N_\xi}, \; k \in \sets K, \; {\bm \ell}\in \sets L \\
         &   \left. \begin{array}{l}
         \; {\bm \lambda}({\bm \ell}) \in \Lambda_K({\bm \ell}) ,\; {\bm \beta}^k({\bm \ell}) \in \reals^L_+, \; k \in \sets K,  \\
        \displaystyle {\bm A}^\top {\bm \alpha}({\bm \ell}) = \sum_{k\in \sets K} {\bm w} \circ {\bm \eta}^k({\bm \ell})  \\
        {\bm A}^\top {\bm \alpha}^k({\bm \ell}) - {\bm H}^\top {\bm \beta}^k({\bm \ell})  + {\bm w} \circ {\bm \eta}^k({\bm \ell}) = {\bm \lambda}_k({\bm \ell}) \left[ {\bm C} \; {\bm x} + {\bm D} \; {\bm w} + {\bm Q} \; {\bm y}^k  \right]  \quad \forall k \in \sets K : {\bm \ell}_k = 0 \\
        {\bm A}^\top {\bm \alpha}^k({\bm \ell}) + [ {\bm H} ]_{ {\bm \ell}_k } {\bm \gamma}_k({\bm \ell}) + {\bm w} \circ {\bm \eta}^k({\bm \ell})   = {\bm \lambda}_k({\bm \ell}) \left[ {\bm C} \; {\bm x} + {\bm D} \; {\bm w} + {\bm Q} \; {\bm y}^k  \right]   \quad \forall k \in \sets K : {\bm \ell}_k \neq 0 \\
         \tau \geq \displaystyle {\bm b}^\top  \left( {\bm \alpha}({\bm \ell}) + \sum_{k \in \sets K}  {\bm \alpha}^k({\bm \ell}) \right) - \sum_{ \begin{smallmatrix} k \in \sets K  : \\ {\bm \ell}_k = 0 \end{smallmatrix}} ( {\bm T} {\bm x} + {\bm V} {\bm w} + {\bm W}{\bm y}^k )^\top {\bm \beta}^k({\bm \ell}) \\
         \qquad \qquad \qquad \displaystyle + \sum_{ \begin{smallmatrix} k \in \sets K  : \\ {\bm \ell}_k \neq 0 \end{smallmatrix}} \left( \left[{\bm T} {\bm x} + {\bm V} {\bm w} + {\bm W}{\bm y}^k \right]_{{\bm \ell}_k} - \epsilon \right) {\bm \gamma}_k({\bm \ell})
         \end{array} \right\} \quad \forall {\bm \ell} \in \partial \sets L  \\
         %
         %
         &   \left. \begin{array}{l}
          \displaystyle {\bm A}^\top {\bm \alpha}({\bm \ell}) = \sum_{k \in \sets K} {\bm w} \circ {\bm \eta}^k({\bm \ell})   \\
          {\bm A}^\top {\bm \alpha}^k({\bm \ell}) +  [ {\bm H} ]_{ {\bm \ell}_k }  {\bm \gamma}_k({\bm \ell}) + {\bm w} \circ {\bm \eta}^k ({\bm \ell})  = {\bm 0} \quad \forall k \in \sets K \\
          \displaystyle {\bm b}^\top  \left( {\bm \alpha}({\bm \ell}) + \sum_{k \in \sets K} {\bm \alpha}^k({\bm \ell}) \right) + \sum_{ k \in \sets K } \left( \left[{\bm T} {\bm x} + {\bm V} {\bm w} + {\bm W}{\bm y}^k \right]_{{\bm \ell}_k} - \epsilon \right) {\bm \gamma}_k({\bm \ell}) \leq -1 
         \end{array} \right\}
         \quad \forall {\bm \ell}\in \sets L_+,
    \end{array}
\label{eq:Kadapt_MBLP}
\end{equation}
where $\Lambda_K({\bm \ell}) := \{ {\bm \lambda} \in \mathbb R_+^K \; : \;  \1^\top {\bm \lambda} = 1, \; {\bm \lambda}_k = 0  \;\; \forall k \in \sets K : {\bm \ell}_k \neq 0 \}$, $\partial \sets L:= \{ {\bm \ell} \in \sets L : {\bm \ell} \ngtr {\bm 0}  \}$ and $\sets L_+ := \{ {\bm \ell} \in \sets L : {\bm \ell} > {\bm 0} \}$ denote the sets for which the decision $({\bm x},{\bm w},\{{\bm y}_k\}_{k\in \sets K})$ satisfies or violates the second-stage constraints in Problem~\eqref{eq:Kadapt_ell}, respectively.
\label{thm:Kadapt_MBLP}
\end{theorem}

\newpv{As in the case of objective uncertainty, if $\sets X \subseteq \{0,1\}^{N_x}$ and $\sets Y \subseteq \{0,1\}^{N_y}$, then Problem~\eqref{eq:Kadapt_MBLP} is equivalent an MBLP involving a suitably chosen ``big-$M$'' constant.} Similar to the robust counterpart resulting from the decision rule approximation proposed in~\cite{DDI_VKB}, Problem~\eqref{eq:Kadapt_MBLP} presents a number of constraints and decision variables that is \emph{exponential} in the approximation parameter, in this case~$K$. Relative to the prepartitioning approach from~\cite{DDI_VKB}, our method does however present a number of distinct advantages. First, the trade-off between approximation quality and computational tractability is controlled using a \emph{single} design parameter; in contrast, in the prepartitioning approach, the number of design parameters equals the number of observable uncertain parameters. Second, as we increase $K$, the quality of the approximation improves in our case, whereas increasing the number of breakpoints along a given direction does not necessarily yield to improvements in the prepartitioning approach. Finally, to identify breakpoint configurations resulting in low optimality gap, a large number of optimization problems need to be solved. 

\begin{remark}
Theorem~\ref{thm:Kadapt_MBLP} directly generalizes to instances of Problem~\eqref{eq:Kadapt} where the technology and recourse matrices ${\bm T}$, ${\bm V}$, and ${\bm W}$ depend on ${\bm \xi}$. Indeed, it suffices to absorb the coefficients of any uncertain terms in ${\bm T}$, ${\bm V}$, and ${\bm W}$ in the right-hand side matrix ${\bm H}$.
\label{rmk:random_rec}
\end{remark}


\begin{observation}
Suppose that we are only in the presence of exogenous uncertainty, i.e., ${\bm w}=\1$, ${\bm D}={\bm 0}$, and ${\bm V}={\bm 0}$. Then, Problem~\eqref{eq:Kadapt_MBLP} reduces to the MBLP formulation of the K-adaptability problem with constraint uncertainty and with only exogenous uncertain parameters from~\cite{Hanasusanto2015}. In particular, in the case of constraint uncertainty, \cite{Hanasusanto2015} also require that the first stage variables ${\bm x}$ be binary.
\label{obs:grani_equivalence_DVzero}
\end{observation}
%


\section{The Case of Piecewise Linear Convex Objective}
\label{sec:pwl_convex_objective}

In this section, we investigate two-stage robust optimization problems with DDID and objective uncertainty where the objective function is given as the maximum of finitely many linear functions.

\subsection{Problem Formulation}
\label{sec:max_objective_formulation}

A piecewise linear convex objective function can be written compactly as the maximum of finitely many linear functions of ${\bm \xi}$ and $({\bm x},{\bm w},{\bm y})$, being expressible as
\begin{equation}
\displaystyle \max_{i \in \sets I} \;\; {\bm \xi}^\top {\bm C}^i \; {\bm x} + {\bm \xi}^\top {\bm D}^i \; {\bm w} + {\bm \xi}^\top {\bm Q}^i \; {\bm y},    
\label{eq:regret_max_objective}
\end{equation}
where ${\bm C}^i \in \reals^{N_\xi \times N_x}$, ${\bm D}^i \in \reals^{N_\xi \times N_\xi}$, and ${\bm Q}^i \in \reals^{N_\xi \times N_y}$, $i \in \sets I$, $\sets I \subseteq \mathbb N$. A two-stage robust optimization problem with DDID, objective function given by~\eqref{eq:regret_max_objective}, and objective uncertainty is then expressible as
\begin{equation}\renewcommand{\arraystretch}{1.5}
\tag{$\mathcal{PO}^{\rm{PWL}}$}
    \begin{array}{cl}
         \min & \;\; \displaystyle \max_{\overline {\bm \xi} \in \Xi} \;\; \min_{ {\bm y} \in \sets Y } \; 
         \max_{ {\bm \xi} \in \Xi({\bm w},\overline {\bm \xi}) } \; \; 
         \left\{  \max_{i \in \sets I} \;\; {\bm \xi}^\top {\bm C}^i \; {\bm x} + {\bm \xi}^\top {\bm D}^i \; {\bm w} + {\bm \xi}^\top {\bm Q}^i \; {\bm y} 
         \right\}  \\
         \st & \;\; {\bm x} \in \sets X, \; {\bm w} \in \sets W.
    \end{array}
\label{eq:min_max_regret_general_mo}
\end{equation}
Note that, as in Section~\ref{sec:kadaptability_objective}, our framework remains applicable in the presence of joint deterministic constraints on the first and second stage variables. We omit these to minimize notational overhead.


\subsection{$K$-Adaptability Approximation \&  MBLP Reformulation}
\label{sec:WCAR_Kadapt}

The $K$-adaptability counterpart of Problem~\eqref{eq:min_max_regret_general_mo} reads
\begin{equation}\renewcommand{\arraystretch}{1.5}
\tag{$\mathcal{PO}^{\rm{PWL}}_K$}
    \begin{array}{cl}
         \min & \;\; \displaystyle \max_{\overline {\bm \xi} \in \Xi} \;\;\min_{ k \in \sets K } \; 
         \max_{ {\bm \xi} \in \Xi({\bm w},\overline {\bm \xi}) } \; \;  \left\{
         \max_{i \in \sets I} \;\; {\bm \xi}^\top {\bm C}^i \; {\bm x} + {\bm \xi}^\top {\bm D}^i \; {\bm w} + {\bm \xi}^\top {\bm Q}^i \; {\bm y}^k
         \right\}  \\
         \st & \;\; {\bm x} \in \sets X, \; {\bm w} \in \sets W, \; {\bm y}^k \in \sets Y, \; k \in \sets K.
    \end{array}
\label{eq:min_max_regret_general_mo_Kadapt}
\end{equation}
We begin this reformulation by the following lemma, which parallels Lemma~\ref{lem:Kadapt_cstr_min_max_min}, and shows that we can exchange the order of the inner min and max in formulation~\eqref{eq:min_max_regret_general_mo_Kadapt}, by indexing ${\bm \xi}$ by $k$.
\begin{lemma}
The $K$-adaptability counterpart of Problem~\eqref{eq:min_max_regret_general_mo_Kadapt} is equivalent to
\begin{equation}\renewcommand{\arraystretch}{1.5}
    \begin{array}{cl}
         \minimize & \quad \displaystyle \max_{ \{ {\bm \xi}^k \}_{k \in \sets K} \in \Xi^K({\bm w}) } \;\;\min_{ k \in \sets K } \;\;  \left\{  \displaystyle \max_{i \in \sets I} \;\;\; ({\bm \xi}^k)^\top {\bm C}^i \; {\bm x} + ({\bm \xi}^k)^\top {\bm D}^i \; {\bm w} + ({\bm \xi}^k)^\top {\bm Q}^i \; {\bm y}^k  \right\}  \\
         \subjectto & \quad {\bm x} \in \sets X, \; {\bm w} \in \sets W, \; {\bm y}^k \in \sets Y, \; k \in \sets K.
    \end{array}
\label{eq:Kadapt_min_max_regret_mo}
\end{equation}
\label{lem:Kadapt_min_max_regret_mo}
\end{lemma}
Next, by leveraging Lemma~\ref{lem:Kadapt_min_max_regret_mo}, we are able to reformulate Problem~\eqref{eq:Kadapt_min_max_regret_mo} exactly as an MBLP.


\begin{theorem}
Problem~\eqref{eq:min_max_regret_general_mo_Kadapt} is equivalent to the bilinear program
\begin{equation}\renewcommand{\arraystretch}{1.5}
\begin{array}{cl}
     \displaystyle \minimize  & \quad \tau \\
      \subjectto & \quad \tau \in \reals, \;  {\bm x} \in \sets X, \; {\bm w} \in \sets W, \; {\bm y}^k \in \sets Y, \; k \in \sets K \\
      & \quad {\bm \alpha}^{\bm i} \in \mathbb R^K_+, \; {\bm \beta}^{\bm i} \in \mathbb R^R_+, \; {\bm \beta}^{{\bm i},k} \in \mathbb R^R_+, \; {\bm \gamma}^{{\bm i},k} \in \mathbb R^{N_\xi}, \; \forall k \in \mathcal K, \; {\bm i} \in \sets I^K \\
 & \!\! \quad  \left.  \begin{array}{l}
  \displaystyle \tau \; \geq \; \bm b^\top\bm \beta^{\bm i} + \sum_{k\in\sets K}\bm b^\top\bm\beta^{{\bm i},k} \\
  \1^\top {\bm \alpha}^{\bm i}  = 1  \\
 {\bm A}^\top{\bm \beta}^{{\bm i}, k} +  {\bm w} \circ {\bm \gamma}^{{\bm i}, k} = {\bm \alpha}^{\bm i}_k \left( {\bm C^{{\bm i}_k}} {\bm x} + {\bm D^{{\bm i}_k}} {\bm w} + {\bm Q^{{\bm i}_k}} {\bm y}^k \right) \quad \forall k \in \mathcal K \\
 {\bm A}^\top{\bm \beta}^{\bm i}  =  \displaystyle \sum_{k \in \mathcal K}  {\bm w} \circ {\bm \gamma}^{{\bm i}, k}
 \end{array} \quad \right\} \quad  \forall {\bm i} \in \sets I^K,
\end{array}
\label{eq:Kadapt_mo_MBLP}
\end{equation}
which can be written as an MBLP, provided $\sets X\subseteq\{0,1\}^{N_x}$.
\label{thm:Kadapt_mo_MBLP}
\end{theorem}

Albeit Problem~\eqref{eq:Kadapt_mo_MBLP} is an MBLP, it presents an exponential number of decision variables and constraints making it difficult to solve directly using off-the-shelf solvers even when~$K$ is only moderately large ($K \gtrapprox 4$). In the remainder of this section, we exploit the specific structure of Problem~\eqref{eq:min_max_regret_general_mo} to solve its $K$-adaptability counterpart \emph{exactly} by reformulating it as an MBLP that presents an attractive structure amenable to decomposition techniques.


\subsection{``Column-and-Constraint Generation'' Algorithm}
\label{sec:max_objective_ccg}

Column-and-constraint generation techniques are a popular approach for addressing problems that possess an exponential number of decision variables and constraints while presenting a decomposable structure, see e.g., \cite{Fischetti1997,Lobel1998,Carvalho1999,Mamer2000,Feillet2010,SADYKOV2011357,zeng2013solving,Muter2013}, and \cite{Muter2018}. 
We propose a new column-and-constraint generation algorithm to solve the $K$-adaptability counterpart~\eqref{eq:min_max_regret_general_mo_Kadapt} based on its reformulation~\eqref{eq:Kadapt_mo_MBLP}. The key idea is to decompose the problem into a relaxed master problem and a series of subproblems indexed by ${\bm i} \in \sets I^K$. The master problem initially only involves the first stage constraints and a \emph{single auxiliary MBLP} is used to iteratively identify indices~${\bm i} \in \sets I^K$ for which the solution to the relaxed master problem becomes infeasible when plugged into subproblem~${\bm i}$. Constraints associated with infeasible subproblems are added to the master problem and the procedure continues until convergence. We detail this procedure in Electronic Companion~\ref{sec:EC_pwl} where we also show that certain classes of two-stage robust optimization problems that seek to minimize the ``worst-case absolute regret'' criterion can be written in the form~\eqref{eq:min_max_regret_general_mo}. In Section~\ref{sec:preference_elicitation_KAS}, we leverage the column and constraint generation algorithm and this observation to solve an active preference elicitation problem that seeks to recommend kidney allocation policies with least possible worst-case regret.


\section{The Multi-Stage Case with Objective Uncertainty}
\label{sec:multistage}

We now show that many of our results generalize to the multi-stage case. To this end, we propose a novel formulation of multi-stage robust optimization problems with DDID. This formulation will enable us to generalize the $K$-adaptability approximation approach to the multi-stage setting.


\subsection{Multi-Stage Robust Optimization with Exogenous Uncertainty}
\label{sec:multi_exogenous}

In the literature, and similar to the two-stage case, there are (broadly speaking) two formulations of a generic multi-stage robust optimization problem with exogenous uncertainty over the planning horizon $\sets T:=\{1,\ldots,T\}$. These differ in the way in which the ability for the time $t$ decisions to adapt to the history of observed parameter realizations is modeled.

\paragraph{Decision Rule Formulation.}

In the first model, one optimizes today over all recourse actions ${\bm y}^t({\bm \xi}) \in \reals^{N_{y_t}}$ that will be taken in each realization of ${\bm \xi} \in \Xi$. Under this modeling paradigm, a multi-stage robust optimization problem with \emph{exogenous} uncertainty is expressible as
\begin{equation}\renewcommand{\arraystretch}{1.5}
    \begin{array}{cl}
    \text{minimize} &  \quad  \displaystyle \max_{ {\bm \xi} \in \Xi } \;\;\; \sum_{t \in \sets T} {\bm \xi}^\top {\bm Q}^t  \; {\bm y}^t ({\bm \xi}) \\
    \text{subject to} & \quad {\bm y}^t \in \mathcal L_{N_\xi}^{N_{y_t}} \quad \forall t \in \sets T \\
        & \quad \!\!\! \left. 
        \begin{array}{l}
        {\bm y}^t({\bm \xi}) \in \sets Y_t \quad \forall t \in \sets T \\
        \displaystyle \sum_{t \in \sets T} {\bm W}^{t} {\bm y}^t({\bm \xi}) \leq {\bm H}{\bm \xi} \\
        {\bm y}^t ({\bm \xi}) = {\bm y}^t( {\bm w}^{t-1} \circ {\bm \xi} )  \quad \forall t \in \sets T
        \end{array} \quad \right\} \quad \forall {\bm \xi} \in \Xi,
    \end{array}
\label{eq:exo_multistage_1}
\end{equation}
where ${\bm Q}^t \in \reals^{N_\xi \times N_{y_t}}$, ${\bm W}^{t} \in \reals^{L_t \times N_{y_t}}$, and ${\bm H}^t \in \reals^{L_t \times N_\xi}$. The \emph{fixed} binary vector ${\bm w}^t \in \{0,1\}^{N_\xi}$ represents the \emph{information base} at time $t+1$, i.e., it encodes the information revealed up to (and including) time $t$. Thus, ${\bm w}^t_i = 1$ if and only if ${\bm \xi}_i$ has been observed at some time $\tau \in \{0,\ldots,t\}$. As information cannot be forgotten, it holds that ${\bm w}^t \geq {\bm w}^{t-1}$ for all $t\in \sets T$. The last constraint in Problem~\eqref{eq:exo_multistage_1} ensures that the decisions ${\bm y}^t$, $t\in \sets T$, are non-anticipative: it stipulates that ${\bm y}^t$ can only depend on those parameters that have been observed up to and including time $t-1$.

\paragraph{Dynamic Formulation.}

In the second model, the recourse decisions ${\bm y}^t$ are optimized explicitly \emph{after} nature is done making a decision. Under this modeling paradigm, a generic multi-stage robust problem with exogenous uncertainty is expressible as:
\begin{equation}\renewcommand{\arraystretch}{1.5}
    \begin{array}{ccl}
\displaystyle \min_{ {\bm y}^1 \in \sets Y_1 } \;\; \max_{ {\bm \xi}^1 \in \Xi } \;\; \min_{ {\bm y}^2 \in \sets Y_2 }  \;\; \max_{ {\bm \xi}^2 \in \Xi ({\bm w}^1, {\bm \xi}^1) } \;\;  \cdots  \;\;  
& \displaystyle \min_{ {\bm y}^T \in \sets Y_T } & \displaystyle \max_{ {\bm \xi}^T \in \Xi ({\bm w}^{T-1}, {\bm \xi}^{T-1}) } \;\;\; \sum_{t \in \sets T} ({\bm \xi}^T)^\top {\bm Q}^t  \; {\bm y}^t \\
& \text{s.t.} &  \displaystyle \sum_{t\in \sets T} {\bm W}^{t} {\bm y}^t \leq {\bm H}({\bm \xi}^T) \quad \forall {\bm \xi}^T \in \Xi( {\bm w}^{T-1}, {\bm \xi}^{T-1}),
    \end{array}
\label{eq:exo_multistage_2}
\end{equation}
where  
$
\Xi({\bm w}^{t-1},{\bm \xi}^{t-1}) : = \left\{ {\bm \xi}^t \in \Xi \; : \; {\bm w}^{t-1} \circ {\bm \xi}^t = {\bm w}^{t-1} \circ {\bm \xi}^{t-1} \right\} \quad \forall t \in \sets T. 
$
This set stipulates that, given the information base ${\bm w}^{t-1} \in \mathbb R^{n_\xi}$ for time $t$ and the associated uncertainty vector ${\bm \xi}^{t-1} \in \Xi \subseteq \mathbb R^{n_\xi}$, nature can select any vector ${\bm \xi}^t \in \Xi$ at time $t$ whose elements associated with parameters it chose prior to time $t$ do not change (i.e., are equal to the corresponding elements chosen in the past).

We state the following theorem without proof.
\begin{theorem}
Problems~\eqref{eq:exo_multistage_1} and~\eqref{eq:exo_multistage_2} are equivalent.
\label{thm:exo_multistage_equiv}
\end{theorem}

\subsection{Multi-Stage Robust Optimization with DDID}
\label{sec:multi_endogenous}

In this section, we investigate a variant of Problem~\eqref{eq:exo_multistage_1} (and accordingly~\eqref{eq:exo_multistage_2}) that enjoys much greater modeling flexibility since the time of information discovery (i.e., the information base) is kept flexible. Thus, we interpret the information base ${\bm w}^t \in \sets W_t \subseteq \{0,1\}^{N_\xi}$ as a decision variable, which is allowed to depend on~${\bm \xi}$. The set $\sets W_t$ may incorporate constraints stipulating, for example, that a specific uncertain parameter can only be observed after a certain stage or that an uncertain parameter can only be observed if another one has, etc. We assume that a cost is incurred for including uncertain parameters in the information base (equivalently, for observing uncertain parameters) and that the observation decisions ${\bm w}^t$ also impact the constraints through the additional term $\sum_{t\in \sets T} {\bm V}^t {\bm w}^t$, where ${\bm V}^t \in \reals^{L_t \times N_\xi}$. As before, we propose two equivalent models for multi-stage robust problems with DDID which differ in the way the ability for the time $t$ decisions to depend on the history of parameter realizations is modeled.

\paragraph{Decision Rule Formulation.} In the first model, one optimizes today over all recourse actions ${\bm w}^t({\bm \xi}) \in \reals^{N_\xi}$ and ${\bm y}^t({\bm \xi}) \in \reals^{N_{y_t}}$ that will be taken in each realization of ${\bm \xi} \in \Xi$. Under this modeling paradigm, a multi-stage robust optimization problem with decision-dependent information discovery, originally proposed in~\cite{DDI_VKB}, reads
\begin{equation}\renewcommand{\arraystretch}{1.5}
    \begin{array}{cl}
    \text{minimize} &  \quad  \displaystyle \max_{ {\bm \xi} \in \Xi } \;\;\; \sum_{t \in \sets T} {\bm \xi}^\top {\bm D}^t  \; {\bm w}^t ({\bm \xi})  + {\bm \xi}^\top {\bm Q}^t  \; {\bm y}^t ({\bm \xi}) \\
    \text{subject to} & \quad {\bm w}^t \in \mathcal L_{N_\xi}^{N_\xi}, \; {\bm y}^t \in \mathcal L_{N_\xi}^{N_{y_t}} \quad \forall t \in \sets T \\
        & \quad \!\!\! \left. 
        \begin{array}{l}
        {\bm w}^t({\bm \xi}) \in \sets W_t, \; {\bm y}^t({\bm \xi}) \in \sets Y_t \quad \forall t \in \sets T  \\
        \displaystyle  \sum_{t \in \sets T} {\bm V}^{t} {\bm w}^t({\bm \xi}) + {\bm W}^{t} {\bm y}^t({\bm \xi}) \leq {\bm H}{\bm \xi} \\
        {\bm w}^t ({\bm \xi}) \geq {\bm w}^{t-1} ({\bm \xi}) \quad \forall t \in \sets T  \\
        {\bm y}^t ({\bm \xi}) = {\bm y}^t( {\bm w}^{t-1}({\bm \xi}) \circ {\bm \xi} ) \quad \forall t \in \sets T   \\
        {\bm w}^t ({\bm \xi}) = {\bm w}^t( {\bm w}^{t-1}({\bm \xi}) \circ {\bm \xi} ) \quad \forall t \in \sets T  
        \end{array} \quad \right\} \quad \forall {\bm \xi} \in \Xi,
    \end{array}
\label{eq:endo_multistage_1}
\end{equation}
where ${\bm w}^0({\bm \xi})={\bm w}^0$ for all ${\bm \xi} \in \Xi$ and ${\bm w}^0$ is given and encodes the information available at the beginning of the planning horizon. \newpv{The above formulation can be used to model problems involving also some exogenous uncertain parameters by restricting~${\bm w}^t_i({\bm \xi})$ to equal either 1 or 0 for all ${\bm \xi}$ depending on whether or not the \emph{exogenous} uncertain parameter~${\bm \xi}_i$ is observed on or before time~$t$. These restrictions can be conveniently added as constraints to the sets~$\mathcal W_t$.}

\paragraph{Dynamic Formulation.} In the second model, the recourse decisions ${\bm w}^t$ and ${\bm y}^t$ are optimized explicitly \emph{after} nature is done selecting the parameters we have chosen to observe in the past. Under this modeling paradigm, a generic multi-stage robust problem with DDID is expressible as:
\begin{equation}\renewcommand{\arraystretch}{1.5}
\tag{$\mathcal{MP}$}
    \begin{array}{ccl}
\displaystyle \min_{ \begin{smallmatrix} {\bm y}^1 \in \sets Y_1 \\ {\bm w}^1 \in \sets W_1  \end{smallmatrix} } \;\; \max_{ {\bm \xi}^1 \in \Xi } \;\; \min_{ \begin{smallmatrix}  {\bm y}^2 \in \sets Y_2 \\ {\bm w}^2 \in \sets W_2 \\ {\bm w}^2 \geq {\bm w}^1 \end{smallmatrix} }  \;\; \max_{ {\bm \xi}^2 \in \Xi ({\bm w}^1, {\bm \xi}^1) }  \;\;  \cdots  \;\; & \displaystyle \min_{ {\bm y}^T \in \sets Y_T } & \displaystyle \max_{ {\bm \xi}^T \in \Xi ({\bm w}^{T-1}, {\bm \xi}^{T-1}) } \;\;\; \sum_{t \in \sets T}  ({\bm \xi}^T)^\top {\bm D}^t  \; {\bm w}^t + ({\bm \xi}^T)^\top {\bm Q}^t  \; {\bm y}^t \\
& \text{s.t.} &  \displaystyle \sum_{t\in \sets T}{\bm V}^{t} {\bm w}^t + {\bm W}^{t} {\bm y}^t \leq {\bm H}{\bm \xi}^T \quad \forall {\bm \xi}^T \in \Xi( {\bm w}^{T-1}, {\bm \xi}^{T-1}).
    \end{array}
\label{eq:endo_multistage_2}
\end{equation}
%

Similarly to the exogenous case, it can be shown that the two models above are equivalent.
\begin{theorem}
Problems~\eqref{eq:endo_multistage_1} and~\eqref{eq:endo_multistage_2} are equivalent.
\label{thm:endo_multistage_equiv}
\end{theorem}
The proof of Theorem~\ref{thm:endo_multistage_equiv} follows by applying Theorem~\ref{thm:endo_equiv} recursively and we thus omit it.

\subsection{$K$-Adaptability for Multi-Stage Problems with DDID}
\label{sec:Kadapt_multi_endo}

We henceforth propose to approximate Problem~\eqref{eq:endo_multistage_2} with its $K$-adaptability counterpart, whereby $K$ candidate policies are selected here-and-now (for each time period) and the best of these policies is selected, in an adaptive fashion, at each stage. To streamline presentation, we focus on the case where Problem~\eqref{eq:endo_multistage_2} presents only objective uncertainty. Thus, the $K$-adaptability counterpart of the multi-stage robust problem~\eqref{eq:endo_multistage_2} with DDID is expressible as
\begin{equation}\renewcommand{\arraystretch}{1.5}
\tag{$\mathcal{MPO}_K$}
    \begin{array}{cl}
\displaystyle \min_{ k_1 \in \sets K } & \;\; \displaystyle \max_{ {\bm \xi}^1 \in \Xi } \;\; \min_{ k_2 \in \sets K }  \;\; \max_{ {\bm \xi}^2 \in \Xi ({\bm w}^{1,k_1}, {\bm \xi}^1) }  \;\; \cdots \\
& \qquad \quad \;\; \cdots 
\displaystyle \min_{ k_T \in \sets K } \;\; \displaystyle \max_{ {\bm \xi}^T \in \Xi ({\bm w}^{T-1,k_1 \ldots k_{T-1}}, {\bm \xi}^{T-1}) } \;\; \sum_{t \in \sets T}  ({\bm \xi}^T)^\top {\bm D}^t  \; {\bm w}^{t,k_1 \ldots k_t} + ({\bm \xi}^T)^\top {\bm Q}^t  \; {\bm y}^{t,k_1 \ldots k_t} \\
 \text{s.t.} &  {\bm y}^{t,k_1 \ldots k_t} \in \sets Y_t, \; {\bm w}^{t,k_1 \ldots k_t} \in \sets W_t \quad \forall t \in \sets T, \; k_1 , \ldots , \; k_{t} \in \sets K \\
 & {\bm w}^{t,k_1 \ldots k_t} \geq {\bm w}^{t-1,k_1 \ldots k_{t-1}} \quad \forall t \in \sets T, \; k_1 , \ldots , \; k_{t} \in \sets K \\
 & \displaystyle \sum_{t\in \sets T}{\bm V}^{t} {\bm w}^{t,k_1 \ldots k_t} + {\bm W}^{t} {\bm y}^{t,k_1 \ldots k_t} \leq {\bm h} \quad \forall k_1, \ldots, k_T \in \sets K ,
    \end{array}
\label{eq:endo_multistage_2_Kadapt}
\end{equation}
%
where we have defined ${\bm w}^{0,k}={\bm w}^0$ for all $k\in \sets K$ with ${\bm w}^0_i=1$ if and only if ${\bm \xi}_i$ is observed at the beginning of the planning horizon and, as in the two-stage case, we have moved the deterministic constraints to the first stage. We note that using the same value of $K$ to approximate the decisions in all periods is without loss of generality and is used to minimize notational overhead. In our experiments in Section~\ref{sec:computational_studies}, we allow for different choices of $K_t$ for each $t \in \sets T$. 
\begin{observation}
For any fixed $K$ and decision $(\{ {\bm w}^{t,k_t} \}_{t \in \sets T, k_t \in \sets K} , \{ {\bm y}^{t,k_t} \}_{t \in \sets T, k_t \in \sets K})$, the optimal objective value of the $K$-adaptability problem~\eqref{eq:endo_multistage_2_Kadapt} can be evaluated by solving an LP whose size is exponential in the size of the input; and in particular exponential in $T$.
\label{obs:evaluation_Kadapt_obj_multistage}
\end{observation}
The proof of Observation~\ref{obs:evaluation_Kadapt_obj_multistage} exactly parallels that of Observation~\ref{obs:evaluation_Kadapt_obj} and thus we omit it.


\subsection{Reformulation as a Mixed Binary Linear Program}
\label{sec:Kadapt_multi_reformulation}

In Observation~\ref{obs:evaluation_Kadapt_obj_multistage}, we showed that for any fixed $K$, ${\bm x}$, $\{ {\bm w}^{t,k_t} \}_{t \in \sets T, k_t \in \sets K}$, and $\{ {\bm y}^{t,k_t} \}_{t \in \sets T, k_t \in \sets K}$, the objective function in Problem~\eqref{eq:endo_multistage_2_Kadapt} can be evaluated by means of an exponentially sized LP. By dualizing this LP and linearizing the resulting bilinear terms, we can obtain an equivalent reformulation of Problem~\eqref{eq:endo_multistage_2_Kadapt} in the form of a mixed-binary linear program.

\begin{theorem}
Suppose $\sets Y_t \subseteq \{0,1\}^{N_{y_t}}$ for all $t \in \sets T$. Then, Problem~\eqref{eq:endo_multistage_2_Kadapt} is equivalent to a bilinear program 
that can be readily linearized using standard techniques.
\label{thm:Kadapt_obj_MBLP_multistage}
\end{theorem}
While the MBLP reformulation of Problem~\eqref{eq:endo_multistage_2_Kadapt} presents a number of decision variables and constraints that are exponential in $T$, it presents an attractive decomposable structure that can be leveraged to solve the problem using e.g., nested Bender's decomposition.

\section{Computational Studies on Stylized Instances}
\label{sec:computational_studies}


We investigate the performance of our approach on a variety of robust optimization problems with decision-dependent information discovery. We solve these problems with our proposed methods discussed in Sections~\ref{sec:kadaptability_objective}, \ref{sec:kadaptability_constraint}, and \ref{sec:multistage}. To speed-up computation, for the two-stage problems, we employ a conservative greedy heuristic that uses the solution to problems with smaller~$K$ to solve problems with larger~$K$ more efficiently, see Section~\ref{sec:heuristic_numericals}. This strategy enables us to solve many random instances of problems with large approximation parameters~$K$ (up to~$K=10$). In all our experiments, we compare our method to the state-of-the-art prepartitioning approach from~\cite{DDI_VKB} using the \ROCPP{} platform, see \cite{Vayanos2020_ROCPP}. All of our experiments are performed on the High Performance Computing Cluster of our university. Each job is allotted 64GB of RAM, 16 cores, and a 2.6GHz Xeon processor. All optimization problems are solved using Gurobi v9.0.1. In Sections~\ref{sec:2stage_best_box} and~\ref{sec:randd_portfolio_optimization}, a total time limit of 7,200 seconds is allowed to solve each instance cumulatively across all values of~$K$ for the $K$-adaptability problem and across all breakpoint configurations for the prepartitioning approach. In Section~\ref{sec:Tstage_pandora}, a time limit of 7,200 seconds is imposed on each instance solved. In all our experiments, we set~$M=500$ and~$\epsilon = 0.001$.



\subsection{Two-stage Robust Best Box Selection (Objective Uncertainty)}
\label{sec:2stage_best_box}

The first problem we study is a robust variant of the best box selection problem, see e.g.,~\cite{Gupta2016_SP,Gupta2017_SP} for results on the stochastic version. In this problem, an agent must select one out of~$N$ boxes, indexed in the set $\mathcal{N} := \{1,...,N\}$, each of which contains a prize. The value $\bm\xi_{i}\in \mathbb{R}$ of the prize in each box $i\in\mathcal{N}$ is unknown and will only be revealed if the box is opened. Opening box~$i \in \sets N$ incurs a cost~${\bm c}_i$. In the first stage, the agent can decide whether to open each box $i \in \sets N$ which we indicate with the decision variables $\bm w_i\in\{0,1\}$. Thus, $\bm w_i = 1$ if and only if $\bm \xi_i$ is observed between the first and second decision stages. The total budget available to open boxes is $B$. In the second stage, the agent can choose one of the opened boxes to keep, which we indicate with the decision variable $\bm y_i\in\{0,1\}$, $i \in \sets N$, earning its prize. We assume that the value of box $i\in\mathcal{N}$ is expressible as $\bm\xi_i = (1+\bm\Phi^\top_i\bm\zeta/2)\bm\xi_i^0$, where $\bm\xi_i^0$ corresponds to the nominal value of the prize of box~$i$, $\bm\zeta\in [-1,1]^L$ are~$L$ risk factors, and $\bm\Phi_i\in\mathbb{R}^L$ collects the factor loadings associated with the value of box~$i$. The goal of the agent is to select the boxes to open (first stage decisions) and the box to keep (second stage decision) to maximize the worst-case value of the box kept. The Best Box Selection problem has numerous applications, for example in house purchasing or in candidate interviewing, see e.g., \cite{Singla_PhDThesis}. With the notation above, the problem can be expressed as a two-stage robust optimization problem with decision-dependent information discovery of the form~\eqref{eq:endo_2_obj} as
%
\begin{equation}
\begin{array}{cccl}
     \displaystyle \max_{ {\bm w} \in \{0,1\}^N } & \quad \displaystyle \min_{\overline{\bm \xi} \in\Xi}  & \quad \displaystyle \max_{ {\bm y} \in \{0,1\}^N } & \;\; \left\{ \;\; \displaystyle \min_{ {\bm \xi} \in \Xi( {\bm w} , \overline {\bm \xi})} \;\;  \bm \xi^{\top}\bm y \; : \;  \textbf{e}^\top\bm y = 1, \;  \bm c^\top\bm w\leq B,\; \bm y\leq \bm w  \;\;\right\},
\end{array}
\label{eq:twostage_TheBestBox}
\end{equation}
where
$\Xi \; := \; \left\{{\bm \xi} \in \mathbb R^{N} \; : \; \exists {\bm \zeta} \in [-1,1]^L  \text{ : } {\bm \xi} \; = \; (1 + {\bm \Phi}^\top {\bm \zeta} / 2 ) {\bm \xi}^{0} \right\}$.

We evaluate the performance of our approach on~100 randomly generated instances of Problem~\eqref{eq:twostage_TheBestBox} with $L=4$ risk factors: 20 instances for each $N \in \{10,20,30,40,50\}$. In these instances, $\bm c$ is drawn uniformly at random from the box $[0,10]^N$, we let~$\bm \xi_i^0 = \bm c/5$, and~$B = \textbf{e}^\top \bm c/2$. The matrix~$\bm\Phi$ is sampled uniformly at random from the box $[-1,1]^{N \times L}$. Our computational results across those instances are summarized in Table~\ref{tab:BestBox_table}. From the table, we observe that with the proposed $K$-adaptability approach, all instances (even those involving $N=50$ boxes) solved to optimality with an average solver time no greater than 15.6 seconds across all problem sizes. In contrast, the average solver time of the prepartitioning approach exceeded 475 seconds for $N=20$ boxes and equaled 7028 seconds for $N=50$ boxes, with only 70.1\% of the problems associated with all breakpoint configurations solving within the allotted time on average. In addition, the quality of the best solution identified by the proposed $K$-adaptability solution consistently outperformed that of the best prepartitioning solution. For example, an average improvement of over 148\% over the static solution was exhibited for the $K$-adaptability method for $N=50$, while the prepartitioning solution only resulted in a 128\% improvement. Finally, we note that the smallest value of $K$ needed to achieve saturation in the optimal value of the problem was consistently smaller that the number of subsets needed to obtain the best possible solution in the prepartitioning method, resulting in more interpratable solutions for our proposed approach. For example, for $N=10$ boxes, a value of $K=3$ is sufficient to yield a 177.4\% improvement in optimal value while an average of 8.2 subsets are needed in the prepartitioning approach to achieve a 164.3\% improvement.

\begin{table}[t!] \renewcommand{\arraystretch}{1.1}
\centering
\scriptsize{
\begin{tabular}{@{}|c|c|c|c|c|c|c|c|c|c|c|@{}}
  \toprule 
 & Adapt. & $N=10$, $L=4$ & $N=20$, $L=4$ & $N=30$, $L=4$ & $N=40$, $L=4$ & $N=50$, $L=4$\\ 
  \midrule
\multirow{10}{*}{ \rotatebox[origin=c]{90}{ $K$-adaptability } }  & $K=1$ & 100\%/0.0\%/0s & 100\%/0.0\%/0s & 100\%/0.0\%/0s & 100\%/0.0\%/0.1s & 100\%/0.0\%/0s \\ 
& $K=2$ & 100\%/153.9\%/0s & 100\%/115.6\%/0s & 100\%/116.9\%/1s & 100\%/110.5\%/1s & 100\%/87.3\%/1s \\ 
& $K=3$ & 100\%/177.4\%/0s & 100\%/147.2\%/0s & 100\%/150.1\%/1s & 100\%/135.5\%/3s & 100\%/113.9\%/3s \\ 
& $K=4$ & 100\%/177.4\%/0s & 100\%/156.0\%/1s & 100\%/159.4\%/2s & 100\%/145.8\%/8s & 100\%/122.6\%/7s \\ 
& $K=5$ & 100\%/177.4\%/1s & 100\%/159.0\%/1s & 100\%/164.0\%/3s & 100\%/148.3\%/9s & 100\%/126.2\%/8s \\
& $K=6$ & 100\%/177.4\%/1s & 100\%/159.0\%/2s & 100\%/164.0\%/3s & 100\%/148.3\%/9s & 100\%/126.2\%/9s \\
& $K=7$ & 100\%/177.4\%/1s & 100\%/159.0\%/2s & 100\%/164.0\%/3s & 100\%/148.3\%/10s & 100\%/126.2\%/10s \\
& $K=8$ & 100\%/177.4\%/1s & 100\%/159.0\%/3s & 100\%/164.0\%/4s & 100\%/148.3\%/12s & 100\%/126.2\%/12s \\
& $K=9$ & 100\%/177.4\%/1s & 100\%/159.0\%/3s & 100\%/164.0\%/5s & 100\%/148.3\%/13s & 100\%/126.2\%/14s \\
& $K=10$ & 100\%/177.4\%/2s & 100\%/159.0\%/3s & 100\%/164.0\%/5s & 100\%/148.3\%/15s & 100\%/126.2\%/16s \\ \midrule
    \rotatebox[origin=c]{90}{ \parbox[c]{1.1cm}{\centering  Preparti-tioning } }  & \parbox[c]{1.1cm}{\centering $\leq$ 10 subsets} & \begin{tabular}[c]{@{}c@{}}100\%/164.3\%\\ /24s/8.2\end{tabular} & \begin{tabular}[c]{@{}c@{}}100\%/137.5\%\\/476s/9 \end{tabular}& \begin{tabular}[c]{@{}c@{}}100\%/143.3\%\\/1173s/8.3\end{tabular} & \begin{tabular}[c]{@{}c@{}}99\%/128.1\%\\/4203s/8.6\end{tabular} & \begin{tabular}[c]{@{}c@{}}70\%/109.3\%\\/7028s/8.8\end{tabular} \\
  \bottomrule
\end{tabular}}
\caption{Summary of computational results on the best box selection problem for various choices of~$N$ over 20 randomly generated instances of each size. In the $K$-adaptability part of the table, each entry corresponds to: percentage of instances solved within the time limit / average improvement in the objective value of the $K$-adaptable solution over the static solution / average solution time across all instances. In the \newpv{prepartitioning} part of the table, each entry corresponds to: average percentage of breakpoint configurations that solved within the time limit out of all configurations with cardinality at most 10 / average improvement in the objective value of the best prepartitioning solution found within the time limit relative to that of the static solution / average cumulative solver time / average cardinality of the best solution found within the time limit.} 
\label{tab:BestBox_table}
\end{table}


\newpage
\subsection{\newpv{Preference Elicitation with Real-Valued Recommendations (Real Decisions)}}
\label{sec:pe_rv}

\newpv{
The second problem we consider is a robust active preference elicitation problem where user preferences can be elicited by asking them ``how much'' they like any particular item and where real-valued quantities of multiple items can be recommended after preferences are elicited, see~\cite{Vayanos_ActivePreferences} for a variant where pairwise comparison queries are used instead.} 

\newpv{The building blocks of our framework are candidate items which we index in the set $\sets I:=\{1,\ldots,I\}$. We let ${\bm \phi}^i \in \reals^J$ be the feature vector of item~$i \in \sets I$. We assume that user preferences are cardinal and model them by means of a linear utility function. Specifically, we assume that the utility of item~$i$ is given by $u({\bm \phi}^i)= {\bm u}^\top {\bm \phi}^i + \tilde {\bm \epsilon}_i$, where $\{ \tilde {\bm \epsilon}_i \}_{i\in \sets I}$ are independent identically distributed and ${\bm u}$ is a vector of (unknown) utility function coefficients supported in the set $\mathcal U \subseteq [-1,1]^J$. These assumptions are standard in the literature, see e.g.,~\cite{OHair_LearningPreferences} and~\cite{Boutilier04:Eliciting}. Before making recommendations, the system has the opportunity to make $Q$ queries to the user. Each query is based on one of the candidate items: if query $i \in \sets I$ is chosen, the user is asked ``On a scale from 0 to 1, where 1 is the most anyone could like an item and 0 is the least anyone could like an item, how much do you like policy $i$?'' We denote by ${\bm \xi}_i \in [0,1]$ the answer to query $i$. After the answers to these queries are observed, the system can select $N$ out of the $I$ items to recommend and the quantity ${\bm y}_i \in [0,1]$, $i\in \sets I$, of those items to recommend. The goal of the recommender system is to select $Q$ queries the answers to which will enable the system to recommend a set of items in quantities resulting in greatest possible worst-case utility.


To formulate the preference elicitation problem mathematically we let ${\bm w}_i$, $i \in \sets I$, denote the decision to pose query~$i$, i.e., to observe ${\bm \xi}_i$ before making a recommendation. Thus, $\mathcal W := \left\{ {\bm w} \in \{0,1\}^I \; : \; \textbf{e}^\top {\bm w} = Q  \right\}.$ The set of possible realizations of ${\bm \xi}$ is given by 
$$
\Xi := \left\{ {\bm \xi} \in [0,1]^I \; : \; 
\exists {\bm u} \in [-1,1]^J, \; {\bm \epsilon} \in \mathcal E \text{ such that } {\bm \xi}_i = 
\frac{ \displaystyle {\bm u}^\top {\bm \phi}^i + \max_{j \in \sets I} \| {\bm \phi}^j \|_1 }{  2  \displaystyle \max_{j \in \sets I} \| {\bm \phi}^j \|_1 }
+ {\bm \epsilon}_i \;\; \forall i \in \sets I
\right\},
$$
where the normalization of ${\bm u}^\top {\bm \phi}^i$ ensures that ${\bm \xi}_i$ has the correct interpretation and, in the spirit of modern robust optimization, see e.g.,~\cite{Lorca2016MultistageRU}, we assume that~${\bm \epsilon}$ is valued in the set  
$
\mathcal E := \left\{ {\bm \epsilon} \in \mathbb R^I \; : \;  \sum_{i=1}^I \left|  {\bm \epsilon}_i \right| \; \leq  \; \Gamma \right\},
$
where~$\Gamma$ is a user-specified \emph{budget of uncertainty} parameter. Once the answers to the queries are observed, the recommender system may select the quantity of each item $i \in \sets I$ to recommend which we encode with decisions ${\bm y}_i \in [0,1]$. We let ${\bm z}_i \in \{0,1\}$ indicate if item $i$ is recommended and require that the quantity of items recommended equals 1. Thus, 
$$\mathcal Y := \left\{ {\bm y} \in [0,1]^I \; : \exists{\bm z} \in \{0,1\}^I \text{ such that }  \; \textbf{e}^\top {\bm z} = N, \; {\bm y} \leq {\bm z}, \; \textbf{e}^\top {\bm y} = 1 \right\}.$$ 
With this notation, the preference elicitation problem is expressible as
\begin{equation}
\tag{$\mathcal {WCU}^{\rm{PE}}$}
         \maximize_{{\bm w} \in \sets W} \;\; \;\; \displaystyle \min_{\overline{\bm \xi} \in \Xi}  \; \;\max_{ {\bm y} \in \sets Y } \;  \min_{ {\bm \xi} \in \Xi({\bm w},\overline{\bm \xi}) } \; {\bm \xi}^\top {\bm y}.
    \label{eq:max_min_utility}
\end{equation}
A conservative solution to Problem~\eqref{eq:max_min_utility} can be obtained using the $K$-adaptability approximation scheme discussed in Section~\ref{sec:kadaptability_objective}, by solving the bilinear reformulation~\eqref{eq:Kadapt_obj_bilinear}.}

\begin{table}[t!] \renewcommand{\arraystretch}{1.1}
\centering
\scriptsize{\newpv{
\begin{tabular}{@{}|c|c|c|c|c|c|c|c|c|c|@{}}
  \toprule 
 & Adapt. & $Q=1$, $N=2$, $\Gamma=0.1$ & $Q=3$, $N=3$, $\Gamma=0.3$  & $Q=6$, $N=4$, $\Gamma=0.6$  & $Q=9$, $N=5$, $\Gamma=0.9$ \\ 
  \midrule
\multirow{10}{*}{ \rotatebox[origin=c]{90}{ $K$-adaptability } }  & $K=1$ & 100\%/0.0\%/2s & 100\%/0.0\%/4s & 
100\%/0.0\%/7s & 
100\%/0.0\%/25s  \\ 
 & $K=2$ & 100\%/14.0\%/8s &
 100\%/28.2\%/15s &
 100\%/35.7\%/22.7s & 
 100\%/41.7\%/48s \\
 & $K=3$ & 100\%/15.2\%/22s & 
 100\%/29.0\%/92s & 
 100\%/42.6\%/406.8s &
 100\%/52.2\%/232s  \\ 
 & $K=4$ & 100\%/15.5\%/44s & 
 100\%/29.4\%/221s & 
 100\%/43.5\%/1396s & 
 90\%/53.1\%/2425s   \\ 
 & $K=5$ & 100\%/15.9\%/75s &
 100\%/29.4\%/434s & 
 90\%/43.9\%/3352s & 
 70\%/53.3\%/4839s   \\ 
 & $K=6$ & 100\%/16.1\%/117s & 
 100\%/29.6\%/779s & 
 65\%/44.2\%/5050s & 
 25\%/53.3\%/6295s   \\ 
 & $K=7$ & 100\%/16.6\%/172s & 
 100\%/29.8\%/1279s & 
 35\%/44.3\%/6221s & 
 15\%/53.3\%/6761s  \\ 
 & $K=8$ & 100\%/16.7\%/226s & 
 100\%/30.0\%/1973s & 
 20\%/44.4\%/6594s & 
 10\%/53.3\%/7035s  \\ 
 & $K=9$ & 100\%/17.3\%/294s &
 95\%/30.0\%/2970s & 
 5\%/44.4\%/6941s & 
 5\%/53.3\%/7167s   \\ 
 & $K=10$ & 100\%/17.3\%/375s & 
 85\%/30.0\%/4031s & 
 5\%/44.4\%/7056s & 
 0\%/53.3\%/7200s   \\ \midrule
    \rotatebox[origin=c]{90}{ \parbox[c]{1.1cm}{\centering  Preparti-tioning } } & \parbox[c]{1.1cm}{\centering $\leq$ 10 subsets} & \begin{tabular}[c]{@{}c@{}}3.3\%/15.9\%\\ /7200s/8.1\end{tabular} & \begin{tabular}[c]{@{}c@{}}2.2\%/14.1\%\\/7200s/9.7 \end{tabular}& \begin{tabular}[c]{@{}c@{}}0.3\%/16.0\%\\/7200s/8.9\end{tabular} & \begin{tabular}[c]{@{}c@{}}0.1\%/16.8\%\\/7200s/9.2\end{tabular} \\ 
  \bottomrule
\end{tabular}}}
\caption{\newpv{Summary of computational results on the preference elicitation problem with real-valued recommendations for various choices of~$Q$, $N$, and~$\Gamma$ over 20 randomly generated instances for each setting. The row names and table entries have the same interpretation as in Table~\ref{tab:BestBox_table}.}}
\label{tab:PE_RV_table}
\end{table}

\newpv{
We evaluate the performance of our approach on~80 randomly generated instances of Problem~\eqref{eq:max_min_utility}: 20 instances for each $(Q,N,\Gamma) \in \{(1,2,0.1),(3,3,0.3),(6,4,0.6),(9,5,0.9)\}$. In these instances, $I=30$, $J=15$, and ${\bm \phi}^i$, $i\in \sets I$, are drawn uniformly at random from the box $[-1,1]^{J}$. Our computational results across these instances are summarized in Table~\ref{tab:PE_RV_table}. From the table, we observe that on average the optimal value of our proposed $K$-adaptability method (across all $K$) is greater than that of the best optimal value of the prepartitioning method (across all breakpoint configurations). For example, for the $(Q,N,\Gamma)=(9,5,0.9)$ setting, $K$-adaptability yields an average improvement in optimal value of $53.3\%$ relative to the static solution, whereas prepartitioning only results in an average improvement of $16.8\%$ in the best case. In addition, the solutions obtained by the $K$-adaptability approach in the same time needed to solve for all breakpoint configurations (or to reach the time limit) in the prepartitioning approach are of far better quality. For example, the prepartitioning approach always reached the 7200 seconds time limit for instances of size $(Q,N,\Gamma)=(3,3,0.3)$ with an associated average improvement in optimal value of $14.1\%$. In contrast, within just 15 seconds on average, the $K$-adaptability approach results in an improvement of $28.2\%$ in optimal value on average over the same instances. Finally, and similar to our results on the best box problem in Section~\ref{sec:2stage_best_box}, the average value of~$K$ needed to achieve a solution of quality comparable to that of the best prepartitioning approach is a lot smaller than the number of subsets needed in prepartitioning, implying that $K$-adaptability has more attractive interpretability properties. For example, for $(Q,N,\Gamma)=(6,4,0.6)$, 8.9 subsets are needed by prepartitioning to yield a 16.0\% improvement in optimal value whereas $K=2$ is sufficient for our method to yield an improvement of~35.7\%.}

\subsection{Robust R\&D Project Portfolio Optimization (Constraint Uncertainty)}
\label{sec:randd_portfolio_optimization}

The \newpv{third} problem we investigate is a robust variant of the R\&D project portfolio optimization problem, see e.g., \cite{Solak_RandDportfolios} for a solution approach on the stochastic version. In this problem, an R\&D firm has a pipeline of~$N$ candidate projects indexed in the set $\sets N := \{1,\ldots,N\}$ that it can invest in. The return ${\bm \xi}_i^{\rm r}$ of each project $i \in \sets N$ is uncertain and will only be revealed if the firm chooses to undertake the project. The firm can decide to undertake each project $i\in \sets N$ in year one, indicated by decision ${\bm w}_i^{\rm r} \in \{0,1\}$, in the following year, indicated by decision ${\bm y}_i \in \{0,1\}$, or not at all. Thus, ${\bm w}_i^{\rm r} = 1$ if and only if ${\bm \xi}_i^{\rm r}$ is observed between the first and second years. If the firm chooses to undertake the investment in the second year, it will only realize a known fraction $\theta \in (0,1]$ of the return. Undertaking project~$i$ incurs an unknown cost ${\bm \xi}_i^{\rm c}$ that will only be revealed if the firm chooses to undertake the project. 
The total budget available to invest in projects across the two years is~$B$. We assume that the return and cost of project $i \in \sets N$ are expressible as
$
{\bm \xi}_i^{\rm r} \; = \; (1 + {\bm \Phi}_i^\top {\bm \zeta} / 2 ) {\bm \xi}_i^{\rm r,0}$ and ${\bm \xi}_i^{\rm c} \; = \; (1 + {\bm \Psi}_i^\top {\bm \zeta} / 2 ) {\bm \xi}_i^{\rm c,0},
$
where ${\bm \xi}_i^{\rm r,0}$ and ${\bm \xi}_i^{\rm c,0}$ corresponds to the nominal return and cost for project $i$, respectively, ${\bm \zeta} \in [-1,1]^L$ are~$L$ risk factors, and the vectors ${\bm \Phi}_i \in \reals^L$ and ${\bm \Psi}_i \in \reals^L$ collect the factor loadings for the return and cost of project~$i$, respectively. With this notation, the R\&D project portfolio optimization problem is expressible as a two-stage robust optimization problem with decision-dependent information discovery of the form~\eqref{eq:endo_3} as
\begin{equation}
\begin{array}{cl}
     \maximize & \quad \displaystyle \min_{\overline{\bm \xi} \in\Xi} \; \;  \max_{ {\bm y} \in \{0,1\}^N } \left\{ \min_{ {\bm \xi} \in \Xi( {\bm w} , \overline {\bm \xi})}   ( {\bm w}^{\rm r}  + \theta {\bm y} )^\top {\bm \xi}^{\rm r} \; : \; ( {\bm w}^{\rm r} + {\bm y} )^\top {\bm \xi}^{\rm c} \; \leq \; B,\;\; {\bm w}^{\rm r} + {\bm y} \leq \1 \right\}\\
    \subjectto &\quad {\bm w} = ({\bm w}^{\rm r} , {\bm w}^{\rm r}), \; {\bm w}^{\rm r} \in \{0,1\}^N,
\end{array}
\label{eq:twostage_RandD}
\end{equation}
where
$$
\Xi \; := \; \left\{ ({\bm \xi}^{\rm r}, {\bm \xi}^{\rm c}) \in \mathbb R^{2N} : \exists {\bm \zeta} \in [-1,1]^\newpv{L}  \text{ : } {\bm \xi}_i^{\rm r} \; = \; (1 + {\bm \Phi}_i^\top {\bm \zeta} / 2 ) {\bm \xi}_i^{\rm r,0},\;\; {\bm \xi}_i^{\rm c} \; = \; (1 + {\bm \Psi}_i^\top {\bm \zeta} / 2 ) {\bm \xi}_i^{\rm c,0} , \;\; i =1,\ldots, N \right\}.
$$

We evaluate the performance of our approach on~100 randomly generated instances of Problem~\eqref{eq:twostage_RandD}: 20 instances for each $(N,L) \in \{(5,3),(10,5),(15,8),(20,10),(25,13)\}$. In these instances, $\theta = 0.8$, $\bm \xi^{\rm c,0}$ is drawn uniformly at random from the box $[0,10]^N$, and we let $\bm \xi^{\rm r,0} = \bm \xi^{\rm c,0}/5$ and $B = \textbf{e}^\top \bm c^0/2$. The elements of $\bm\Phi$ and $\bm\Psi$ are uniformly distributed in the interval $[-1,1]$. Our computational results across these instances are summarized in Table~\ref{tab:RDPPO_table}. From the table, we observe that on average the optimal value of our proposed $K$-adaptability method (across all $K$) is greater than that of the best optimal value of the prepartitioning method (across all breakpoint configurations). For example, for the $(N,L)=(20,10)$ setting, $K$-adaptability yields an average improvement in optimal value of $56.6\%$ relative to the static solution, whereas prepartitioning only results in an average improvement of $18.7\%$ in the best case. In addition, the solutions obtained by the $K$-adaptability approach in the same time needed to solve for all breakpoint configurations (or to reach the time limit) in the prepartitioning approach are of far better quality. For example, the prepartitioning approach needed 5771 seconds on average to solve instances of size $(N,L)=(15,8)$ with an associated average improvement in optimal value of $23.8\%$. In contrast, within just 1814 seconds on average, the $K$-adaptability approach results in an improvement of $65.7\%$ in optimal value on average over the same instances. Finally, and similar to our results on the best box problem in Section~\ref{sec:2stage_best_box}, the average value of~$K$ needed to achieve a solution of quality comparable to that of the best prepartitioning approach is a lot smaller than the number of subsets needed in prepartitioning, implying that $K$-adaptability has more attractive interpretability properties. For example, for $(N,L)=(10,5)$, 8.2 subsets are needed by prepartitioning to yield a 35.4\% improvement in optimal value whereas $K=3$ is sufficient for our method to yield an improvement of~46.2\%. 

\begin{table}[t!] \renewcommand{\arraystretch}{1.1}
\centering
\scriptsize{
\begin{tabular}{@{}|c|c|c|c|c|c|c|c|c|c|c|@{}}
  \toprule 
 & Adapt. & $N=5$, $L=3$ & $N=10$, $L=5$ & $N=15$, $L=8$ & $N=20$, $L=10$ & $N=25$, $L=13$ \\ 
  \midrule
\multirow{10}{*}{ \rotatebox[origin=c]{90}{ $K$-adaptability } }  & $K=1$ & 100\%/0.0\%/0s & 100\%/0.0\%/0s & 100\%/0.0\%/1s & 100\%/0.0\%/15s & 100\%/0.0\%/39s \\ 
 & $K=2$ & 100\%/40.8\%/1s & 100\%/32.7\%/2s & 100\%/32.2\%/12s & 100\%/31.0\%/264s & 100\%/31.9\%/1915s \\
 & $K=3$ & 100\%/81.1\%/1s & 100\%/46.2\%/6s & 100\%/50.6\%/55s & 90\%/43.9\%/1543s & 50\%/42.8\%/5633s \\ 
 & $K=4$ & 100\%/102.6\%/3s & 100\%/52.2\%/22s & 100\%/59.4\%/323s & 70\%/51.9\%/3959s & 0\%/44.6\%/7200s \\ 
 & $K=5$ & 100\%/109.2\%/5s & 100\%/57.1\%/103s & 100\%/65.7\%/1814s & 20\%/56.6\%/6661s & 0\%/44.6\%/7200s \\ 
 & $K=6$ & 100\%/110.5\%/10s & 100\%/59.4\%/435s & 40\%/67.6\%/6321s & 0\%/56.6\%/7200s & 0\%/44.6\%/7200s \\ 
 & $K=7$ & 100\%/110.7\%/43s & 100\%/62.5\%/1835s & 5\%/67.8\%/7122s & 0\%/56.6\%/7200s & 0\%/44.6\%/7200s \\ 
 & $K=8$ & 100\%/114.5\%/32s & 60\%/66.6\%/5319s & 0\%/67.8\%/7200s & 0\%/56.6\%/7200s & 0\%/44.6\%/7200s \\ 
 & $K=9$ & 100\%/117.5\%/100s & 10\%/66.9\%/6955s & 0\%/67.8\%/7200s & 0\%/56.6\%/7200s & 0\%/44.6\%/7200s \\ 
 & $K=10$ & 100\%/118.6\%/301s & 0\%/66.9\%/7200s & 0\%/67.8\%/7200s & 0\%/56.6\%/7200s & 0\%/44.6\%/7200s \\ \midrule
    \rotatebox[origin=c]{90}{ \parbox[c]{1.1cm}{\centering  Preparti-tioning } } & \parbox[c]{1.1cm}{\centering $\leq$ 10 subsets} & \begin{tabular}[c]{@{}c@{}}100\%/56.3\%\\ /35s/8.2\end{tabular} & \begin{tabular}[c]{@{}c@{}}100\%/35.4\%\\/810s/8.2 \end{tabular}& \begin{tabular}[c]{@{}c@{}}99\%/23.8\%\\/5771s/9.2\end{tabular} & \begin{tabular}[c]{@{}c@{}}28\%/18.7\%\\/7200s/8.4\end{tabular} & \begin{tabular}[c]{@{}c@{}}5\%/14.3\%\\/7200s/8.2\end{tabular} \\ 
  \bottomrule
\end{tabular}}
\caption{Summary of computational results on the R\&D project portfolio optimization problem for various choices of~$N$ and~$L$ over 20 randomly generated instances of each size. The row names and table entries have the same interpretation as in Table~\ref{tab:BestBox_table}.}
\label{tab:RDPPO_table}
\end{table}




\subsection{Multi-Stage Robust Pandora Box (Multi-Stage Objective Uncertainty)}
\label{sec:Tstage_pandora}

The \newpv{fourth} problem we investigate is a robust variant of the multi-stage Pandora Box problem, see e.g., \cite{Doval2018} and~\cite{Singla_PhDThesis} for the stochastic setting. This problem is similar to best box selection, see Section~\ref{sec:2stage_best_box}, however opening boxes incurs a cost in the objective and no budget constraint is imposed. The planning horizon consists of $T$ periods indexed in the set $\sets T:=\{1,...,T\}$. At the beginning of each period $t \in \sets T$, the agent can select one box to open out of~$N$ available boxes indexed in the set $\sets N:=\{1,\ldots, N\}$, each of which contains a prize. We let $\bm w_i^t\in\{0,1\}$ indicate if box~$i$ has been opened on or before time $t$. The value~${\bm \xi}_i$ of the prize in each box $i$ is unknown at the beginning of the planning horizon and is only observable at time $t\in \sets T$ if the box has been opened before time $t$. Thus, the decisions to open the boxes control the time of information discovery in this problem. Opening box $i \in \sets N$ incurs a cost ${\bm c}_i$. If the agent chooses to not open a box at some time $t\in \sets T$, they may alternatively choose to keep one of the previously opened boxes, earning its prize. We let $\bm y_i^t\in\{0,1\}$ indicate the choice to keep box $i$ at time $t$. Throughout the planning horizon, the agent may keep at most one box. The goal of the agent is to select whether and when to open each box to maximize the worst-case value of the box they choose to keep less the total cost of opening boxes. We assume that $\bm\xi_i = (1+\bm\Phi^\top_i\bm\zeta/2)\bm\xi_i^0$, where $\bm \xi_i^0$ denotes the nominal value of box~$i$, $\bm\zeta\in\mathbb{R}^L$ are $L$ risk factors, and the vectors $\bm\Phi_i\in\mathbb{R}^L$, $i\in\mathcal{N}$, collect the factor loadings for the value of box~$i$. With this notation, Pandora's Box problem is expressible as a multi-stage robust optimization problem with decision-dependent information discovery of the form~\eqref{eq:endo_multistage_2} as
\begin{equation}
\begin{array}{ccl}
\displaystyle \max_{ \begin{smallmatrix} {\bm y}^1 \in \sets Y_1 \\ {\bm w}^1 \in \sets W_1  \end{smallmatrix} } \;\; \min_{ {\bm \xi}^1 \in \Xi } \;\; \max_{ \begin{smallmatrix}  {\bm y}^2 \in \sets Y_2 \\ {\bm w}^2 \in \sets W_2 \\ {\bm w}^2 \geq {\bm w}^1 \end{smallmatrix} }  \;\; \min_{ {\bm \xi}^2 \in \Xi ({\bm w}^1, {\bm \xi}^1) } \;\;  \cdots  \;\; & \displaystyle \max_{ {\bm y}^T \in \sets Y_T } & \;\;\; \displaystyle \min_{ {\bm \xi}^T \in \Xi ({\bm w}^{T-1}, {\bm \xi}^{T-1}) } \;\;\; \sum_{t \in \sets T} \bm\xi^\top \bm y^t-\bm c^\top\bm w^{T-1}\\ 
& \text{s.t.} & \;\;\;  \textbf{e}^\top(\bm w^t-\bm w^{t-1}) \; \leq \; 1-\displaystyle\sum_{\tau\leq t} \textbf{e}^\top\bm y^\tau \quad \forall t \in \sets T\backslash\{T\}  \\
&& \;\;\; \displaystyle\sum_{t\in \sets T}\bm e^\top\bm y^t = 1 ,\quad \bm y^t\leq \bm w^{t-1}\quad \forall t \in \sets T, 
\end{array}
\label{eq:multistage_PandoraBox}
\end{equation}
where $\Xi := \left\{ {\bm \xi} \in \reals^N \; : \; \exists {\bm \zeta} \in [-1,1]^L \; : \; {\bm \xi}_i = (1 + \bm\Phi_i^\top {\bm \zeta}/2) {\bm \xi}_i^0 \;\;\; \forall i \in \sets N \right\}$, $\mathcal Y_t := \{0,1\}^N$, and $\mathcal W_t := \{0,1\}^N$. The first set of constraints ensures that if a box has been kept in the past, no box can be opened. The second set of constraints guarantees that exactly one of the open boxes is kept.

\begin{figure}
    \centering
    \includegraphics[width=1\textwidth]{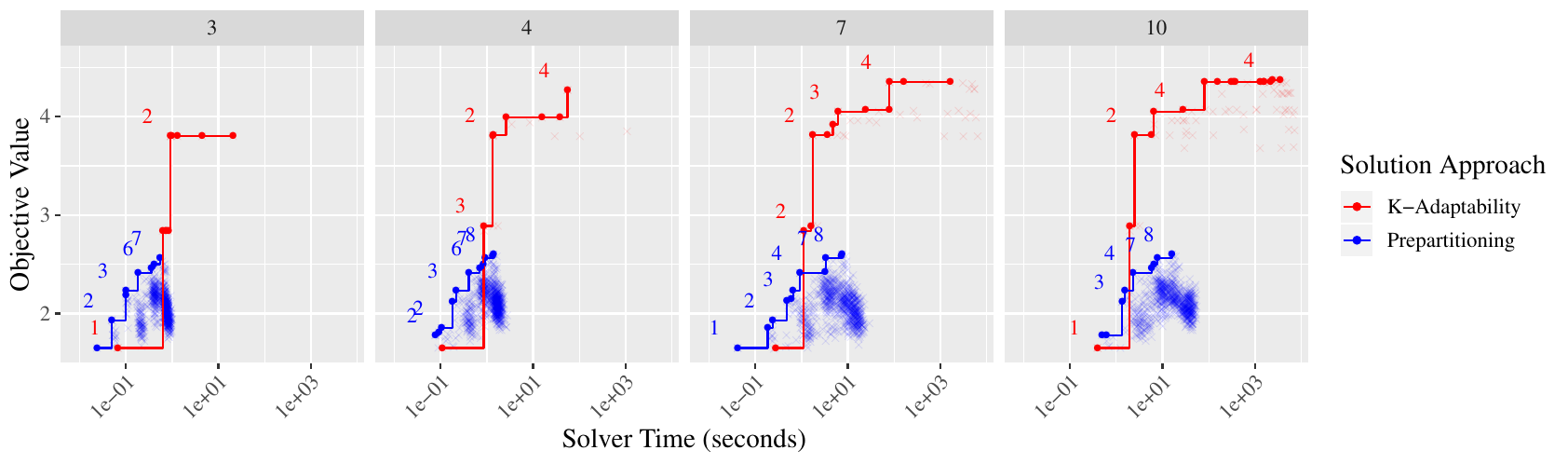}
    \caption{Optimality-scalability results for the multi-stage Pandora box problem~\eqref{eq:multistage_PandoraBox}. Each number on the top of each facet corresponds to the number of time-periods~$T$. Each red dot (and cross) corresponds to a different choice of $(K_1,\ldots,K_T) \in \{1,\ldots,5\}^T$ for the $K$-adaptability problem such that the total number of candidate policies $K_1\times K_2 \times \ldots \times K_T \leq 5$. Each blue dot (and cross) corresponds to a different breakpoint configuration for the prepartitioning approach (we consider 100 different breakpoint configurations drawn randomly from the set of all configurations with cardinality less than 5). The position of each dot is obtained by averaging the solver time and optimal value across the 20 random instances. Whether a point is indicated with a dot or a cross depends on whether it is on the efficient frontier of the configurations that resulted in the highest average optimal value for the given time budget. The numbers next to the efficient points indicate the degree of adaptability in the corresponding solution, i.e., the number of subsets in the prepartitioning approach and the total number of candidate policies in the $K$-adaptability approach.}
    \label{fig:multistage_PandoraBox}
\end{figure}

We evaluate the performance of our approach on 80 randomly generated instances of Problem~\eqref{eq:multistage_PandoraBox} with $N=15$ boxes and $L = 4$ risk factors: 20 instances for each choice of $T\in \{3,4,7,10\}$. In these instances, $\bm \xi^{0}$ is drawn uniformly at random from $[0,10]^N$ and we let $\bm c = \bm \xi^{0}/5$. The elements of $\bm\Phi$ are uniformly distributed in the interval $[-1,1]$. Our computational results across these instances are summarized in Figure~\ref{fig:multistage_PandoraBox}. From the figure, it can be seen that the $K$-adaptability approach results in solutions of far better quality than the prepartitioning approach. For example, the best average optimal value of the prepartitioning approach for the $T=4$ case is 2.60 whereas it is 4.27 for the $K$-adaptability method. Moreover, solutions outperforming the best prepartitioning solution are found on average faster with the $K$-adaptability approach. For example, in the setting $T=10$, the breakpoint configuration that resulted in the best average optimal value had an average solver time of 15.9 seconds, whereas 1.9 seconds on average were needed to solve the $K$-adaptability problem associated with the adaptability configuration that resulted in a higher optimal value on average. Finally, the solutions obtained by the $K$-adaptability approach were consistently more interpretable than those of the prepartitioning method. In particular, across all our experiments, the cardinality of the breakpoint configuration that resulted in the best performance of the prepartitioning method was between 7 and 8 and that solution was consistently outperformed by $K$-adaptable solutions with under~3 candidate policies.

\section{Preference Elicitation to Improve the US Kidney Allocation System}
\label{sec:preference_elicitation_KAS}

In this section, we evaluate our approach on a preference elicitation and recommendation problem that explicitly captures the endogenous nature of the elicitation process.

\subsection{Motivation \& Problem Formulation (Piecewise Linear Convex Objective)}

The motivation for our study is one of the central problems faced by policymakers at the OPTN/UNOS who must periodically make changes to the policy for prioritizing patients on the kidney transplant waiting list for scarce deceased donor kidneys. To tackle this problem, a Kidney Transplantation Committee (KTC) is appointed at the OPTN that examines the outcomes of numerous candidate policies simulated using the Kidney-Pancreas Simulated Allocation Model (KPSAM), a simulator developed by the Scientific Registry of Transplant Recipients (SRTR), see~\cite{KPSAM2015}. The KTC examines the outcomes of the allocation policy alternatives along several dimensions (measures) of fairness and efficiency (e.g., number of recipients by age group, number of deaths by gender) before ultimately committing to one of the alternatives. This process was for example followed in the latest big policy change, see e.g.,~\cite{FinalAnalysis_SRTR2009}. Since selecting one alternative (policy) over many others is a challenging task, in particular when the dimension of each alternative is large, see e.g., \cite{Toubia_2003,Toubia_2004,Toubia_2007} and~\cite{Boutilier04:Eliciting}, we propose a preference elicitation and recommendation framework for identifying a preferred policy using a moderate number of strategically chosen queries.


\newpv{We formulate this problem as a variant of the active preference elicitation problem from Section~\ref{sec:pe_rv} where a single item can be recommended and where we select queries that minimize worst-case regret of the recommendation. Items indexed in the set $\sets I$ correspond to policies where the feature vector ${\bm \phi}^i \in \reals^J$ of policy~$i \in \sets I$ collects various measures of fairness and efficiency of the policy. The problem is expressible mathematically as}
\begin{equation}
\tag{$\mathcal {WCR}^{\rm{PE}}$}
         \minimize_{{\bm w} \in \sets W} \;\; \;\; \displaystyle \max_{\overline{\bm \xi} \in \Xi}  \; \;\min_{ {\bm y} \in \sets Y } \;  \max_{ {\bm \xi} \in \Xi({\bm w},\overline{\bm \xi}) } \left\{  \max_{i \in \mathcal I} \;\; {\bm \xi}_i -   \; {\bm \xi}^\top {\bm y} \;  \right\},
    \label{eq:min_max_regret}
\end{equation}
\newpv{where $\sets W$ and $\Xi$ are as in Section~\ref{sec:pe_rv} and where $\mathcal Y := \left\{ {\bm y} \in \{0,1\}^I \; : \; \textbf{e}^\top {\bm y} = 1 \right\}$.} In this problem, the first part of the objective computes the utility of the best item to offer in hindsight, after the utilities ${\bm \xi}$ have been observed. The second part of the objective corresponds to the worst-case utility of the item recommended when only a portion of the uncertain parameters are observed, as dictated by the vector ${\bm w}$. Problem~\eqref{eq:min_max_regret} can be solved approximately using the $K$-adaptability approximation scheme discussed in Section~\ref{sec:pwl_convex_objective}. Indeed, the regret in Problem~\eqref{eq:min_max_regret} is given as the maximum of finitely many linear functions and Theorem~\ref{thm:Kadapt_mo_MBLP} applies. We note that in this case $| \sets Y | = I$. Thus, solving the $K$-adaptability counterpart of~\eqref{eq:min_max_regret} with $K=I$ recovers an optimal solution to the corresponding original problem.



\subsection{Generating KAS Candidate Policies}
\label{sec:nr_candidate_policies}

We generate the outcomes ${\bm \phi}^i$, $i \in \sets I$, of $I=20$ candidate policies using the KPSAM simulator which we obtained from the SRTR using a modeling window from 01/01/2010 to 12/31/2010. The candidate policies we consider are linear scoring rules that use the patient dialysis time, the life years from transplant score, the Calculated Panel Reactive Antibodies and the age of the patient. For each policy, we record $J=22$ outcomes, including the number of transplants overall, by age, by blood type, by race, and by gender, and the number of deaths by race and by gender. For details on the construction of the policies and for a list of outcomes, see Electronic Companion~\ref{sec:EC_preference_elicitation_KAS}.

\subsection{Numerical Results on KAS Candidate Policies}
\label{sec:nr_real}

We evaluate the performance of our approach on the KAS policies dataset from Section~\ref{sec:nr_candidate_policies}. Throughout our experiments, the $K$-adaptability counterpart of Problem~\eqref{eq:min_max_regret} is solved using the techniques described in Section~\ref{sec:pwl_convex_objective}. To speed-up computation, we also use \newpv{a heuristic adapted from~\cite{subramanyam2017k} and} detailed in Section~\ref{sec:heuristic_numericals}. The tolerance $\delta$ used in the column-and-constraint generation algorithm (see Section~\ref{sec:max_objective_ccg}) is $10^{-5}$. We evaluate the \emph{true} worst-case regret of any given solution ${\bm w}^\star$, which we denote by $r_{\rm{wc}}({\bm w}^\star)$, as follows: we fix ${\bm w}={\bm w}^\star$ in Problem~\eqref{eq:ccg_feas}, where we set $K = I$ and employ all $I$ candidate policies~$\{ {\bm y}^k \}_{k \in \sets K}$ in the set~$\mathcal Y$. As before, we use the \ROCPP{} platform to solve the prepartitioning problem, see \cite{Vayanos2020_ROCPP}. All of our experiments are performed using the same computing resources as in Section~\ref{sec:computational_studies}.

\begin{figure}[t!]
    \begin{center}
    \includegraphics[width=0.7\textwidth]{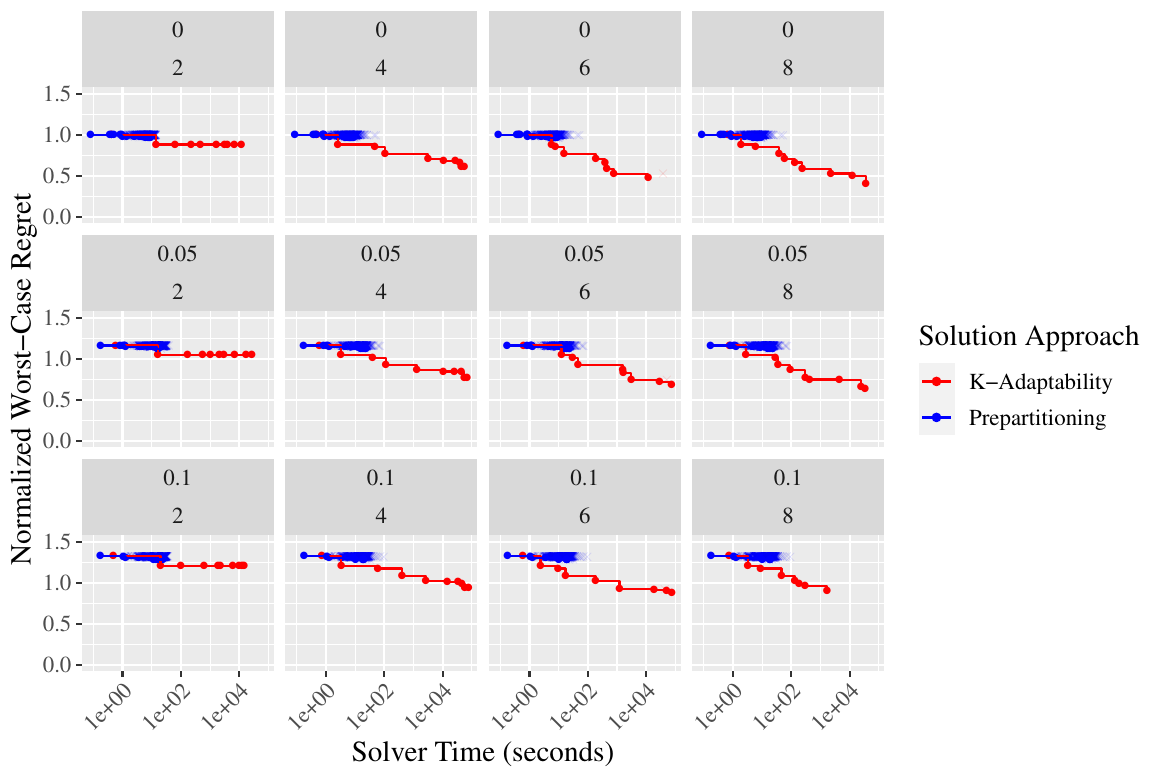}
    \end{center}
    \caption{Optimality-scalability results for the min-max regret preference elicitation problem~\eqref{eq:min_max_regret} on the KAS data. The numbers on each facet correspond to the values of $\Gamma$ (top number) and $Q$ (bottom number). The shapes, lines, and colors have the same interpretation as in Figure~\ref{fig:multistage_PandoraBox}.}
    \label{fig:MinMaxRegret_EF_Real}
\end{figure}

\paragraph{Optimality-Scalability Trade-Off.} We evaluate the trade-off between computational complexity and scalability of our approach. We solve the min-max regret problems as $Q$ and $\Gamma$ are varied in the sets $\{2,4,6,8\}$ and $\{0,0.05,0.1\}$, respectively. The results are summarized in Figure~\ref{fig:MinMaxRegret_EF_Real}. From the figure it can be seen that the $K$-adaptability approach significantly outperforms the prepartitioning approach and static policies are very sub-optimal. In fact, the prepartitioning approach performs comparably to static policies across all settings. On the other hand, with the $K$-adaptability approach, the normalized\footnote{To aid with interpretability, we normalize regret such that the worst-case regret when no question is asked is 1 and the worst-case regret when all questions are asked is 0.} worst-case regret drops to 0.40, 0.68, and 0.9 from 1, 1.16, and 1.32, for $\Gamma=0,0.05,$ and $0.1$, respectively (for $Q=8$). This experiment shows the strength of the $K$-adaptability approach compared to the state of the art.

\begin{figure}[t!]
    \begin{center}
    \includegraphics[width=0.75\textwidth]{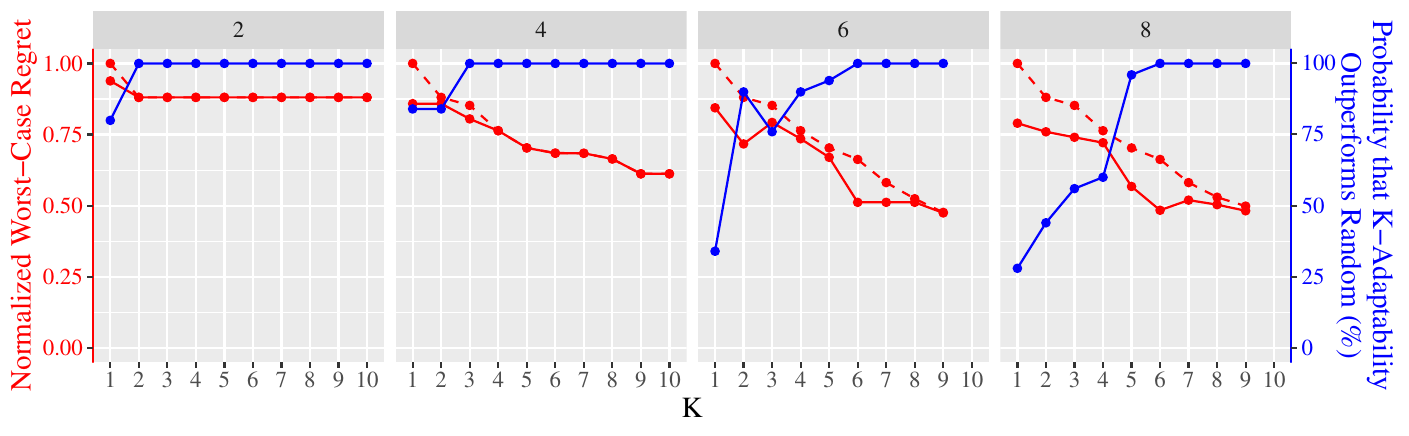}
    \end{center}
    \caption{Results on the performance of the $K$-adaptability approach relative to random elicitation for the min-max regret preference elicitation problem~\eqref{eq:min_max_regret} on the LAHSA dataset. The dashed red line corresponds to the objective value of the $K$-adaptability problem. The red line corresponds to $r_{\rm{wc}}({\bm w}_K^\star)$ where ${\bm w}_K^\star$ is the optimal $K$-adaptable solution. The blue line represents the percentage of time that $r_{wc}({\bm w}_{\rm r})$ was lower than $r_{\rm{wc}}({\bm w}_K^\star)$, where ${\bm w}_{\rm r}$ is a randomly drawn feasible solution.}
    \label{fig:Simulation_MinMaxRegret_Real}
\end{figure}

\paragraph{Performance Relative to Random Elicitation.} We evaluate the benefits of computing near-optimal queries using the $K$-adaptability approximation approach relative to asking questions at random. We compare the \emph{true} performance of a solution to the $K$-adaptability problem, $r_{\rm{wc}}({\bm w}_K^\star)$, to that of~50 questions drawn uniformly at random from the set $\mathcal W$, $r_{wc}({\bm w}_{\rm r})$. The results are summarized on Figure~\ref{fig:Simulation_MinMaxRegret_Real}. From the figure, we see that the probability that the $K$-adaptability solution outperforms random elicitation converges to 1 as~$K$ grows. We observe that, for values of $K$ greater than~5, the $K$-adaptability solution outperforms random elicitation in over 90\% of the cases.

\ACKNOWLEDGMENT{This work was supported primarily by the Operations Engineering Program of the National Science Foundation under NSF Award No.\ 1763108. The authors are grateful to Miss.\ Qing Jin for valuable discussions on implementation issues.} 



\newpage

\typeout{}
\bibliographystyle{informs2014} 
\bibliography{Mendeley.bib}


\ECSwitch


\ECHead{E-Companion}


\section{Companion to Section~\ref{sec:pwl_convex_objective}}
\label{sec:EC_pwl}


\subsection{Column and Constraint Generation Algorithm}
\label{sec:EC_ccg}

We define the following relaxed master problem parameterized by the index set $\widetilde{\sets I} \subseteq \sets I^K$
\begin{equation}\renewcommand{\arraystretch}{1.5}
\tag{$\mathcal{CCG}_{\rm{mstr}}(\widetilde{\sets I})$}
\begin{array}{cl}
     \displaystyle \min  & \quad \tau \\
      \st & \quad \tau \in \reals, \;  {\bm x} \in \sets X, \; {\bm w} \in \sets W, \; {\bm y}^k \in \sets Y, \; k \in \sets K \\
      & \quad {\bm \alpha}^{\bm i} \in \mathbb R^K_+, \; {\bm \beta}^{\bm i} \in \mathbb R^R_+, \; {\bm \beta}^{{\bm i},k} \in \mathbb R^R_+, \; {\bm \gamma}^{{\bm i},k} \in \mathbb R^{N_\xi}, \; \forall k \in \mathcal K, \; {\bm i} \in \widetilde{\sets I} \\
 & \!\! \quad  \left.  \begin{array}{l}
  \displaystyle \tau \; \geq \; \bm b^\top\bm \beta^{\bm i} + \sum_{k\in\sets K}\bm b^\top\bm\beta^{{\bm i},k} \\
  \1^\top {\bm \alpha}^{\bm i}  = 1  \\
 {\bm A}^\top{\bm \beta}^{{\bm i}, k} +  {\bm w} \circ {\bm \gamma}^{{\bm i}, k} = {\bm \alpha}^{\bm i}_k \left( {\bm C^{{\bm i}_k}} {\bm x} + {\bm D^{{\bm i}_k}} {\bm w} + {\bm Q^{{\bm i}_k}} {\bm y}^k \right) \quad \forall k \in \mathcal K \\
 {\bm A}^\top{\bm \beta}^{\bm i}  =  \displaystyle \sum_{k \in \mathcal K}  {\bm w} \circ {\bm \gamma}^{{\bm i}, k}
 \end{array} \quad \right\} \quad  \forall {\bm i} \in \widetilde{\sets I}.
\end{array}
\label{eq:ccg_master}
\end{equation}
Given variables $(\tau,{\bm x},{\bm w}, \{{\bm y}^k\}_{k\in \sets K})$ feasible in the master problem, we define the ${\bm i}$th subproblem, ${\bm i} \in \sets I$, through
\begin{equation}\renewcommand{\arraystretch}{1.5}
\tag{$\mathcal{CCG}^{\bm i}_{\rm{sub}}(\tau,{\bm x},{\bm w}, \{{\bm y}^k\}_{k\in \sets K})$}
\begin{array}{cl}
     \displaystyle \min  & \quad 0 \\
      \st & \quad {\bm \alpha}^{\bm i} \in \mathbb R^K_+, \; {\bm \beta}^{\bm i} \in \mathbb R^R_+, \; {\bm \beta}^{{\bm i},k} \in \mathbb R^R_+, \; {\bm \gamma}^{{\bm i},k} \in \mathbb R^{N_\xi}, \; \forall k \in \mathcal K \\
 & \!\! \quad   \begin{array}{l}
  \displaystyle \tau \; \geq \; \bm b^\top\bm \beta^{\bm i} + \sum_{k\in\sets K}\bm b^\top\bm\beta^{{\bm i},k} \\
  \1^\top {\bm \alpha}^{\bm i}  = 1  \\
 {\bm A}^\top{\bm \beta}^{{\bm i}, k} +  {\bm w} \circ {\bm \gamma}^{{\bm i}, k} = {\bm \alpha}^{\bm i}_k \left( {\bm C^{{\bm i}_k}} {\bm x} + {\bm D^{{\bm i}_k}} {\bm w} + {\bm Q^{{\bm i}_k}} {\bm y}^k \right) \quad \forall k \in \mathcal K \\
 {\bm A}^\top{\bm \beta}^{\bm i}  =  \displaystyle \sum_{k \in \mathcal K}  {\bm w} \circ {\bm \gamma}^{{\bm i}, k}.
 \end{array}
\end{array}
\label{eq:ccg_subproblem}
\end{equation}
An inspection of the Proof of Theorem~\ref{thm:Kadapt_mo_MBLP} reveals that the last three constraints in Problem~\eqref{eq:ccg_subproblem} define the feasible set of the dual of a linear program that is feasible and bounded. Thus, for $\tau$ sufficiently large, Problem~\eqref{eq:ccg_subproblem} will be feasible. 

To identify indices of subproblems~\eqref{eq:ccg_subproblem} that, given a solution $(\tau,{\bm x},{\bm w}, \{{\bm y}^k\}_{k\in \sets K})$ to the relaxed master problem, are infeasible, we solve a \emph{single} feasibility MBLP defined through
\begin{equation}\renewcommand{\arraystretch}{1.5}
\tag{$\mathcal{CCG}_{\rm{feas}}({\bm x},{\bm w}, \{{\bm y}^k\}_{k\in \sets K})$}
\begin{array}{cl}
    \max & \quad \theta \\
    \st & \quad \theta\in \reals,\; \overline{\bm \xi}\in\Xi,\; {\bm \xi^k} \in \Xi({\bm w}, \overline {\bm \xi}), \; k \in \sets K \\
    & \quad {\bm \eta} \in \reals^K, \; {\bm \zeta}^k \in \{0,1\}^I, \; k \in \sets K  \\
    & \quad \theta \; \leq \; {\bm \eta}_k   \quad\forall k \in \sets K \\
    & \!\! \left. \begin{array}{l}
    \quad {\bm \eta}_k \; \geq \; (\bm\xi^k)^\top \bm C^{i}\;{\bm x} + {(\bm \xi^k)}^\top {\bm D^i} \; {\bm w} + {(\bm \xi^k)}^\top {\bm Q^i} \; {\bm y^k} \\
    \quad {\bm \eta}_k \; \leq \;  (\bm\xi^k)^\top \bm C^{i}\;{\bm x} + {(\bm \xi^k)}^\top {\bm D^i} \; {\bm w} + {(\bm \xi^k)}^\top {\bm Q^i} \; {\bm y^k}  + M (1-{\bm \zeta}^k_i ) 
    \end{array} \right\} \; \begin{matrix} \forall i \in \sets I, \\ k \in \sets K \end{matrix} \\
    & \quad \1^\top {\bm \zeta}^k  = 1 \quad \forall k \in \sets K.
\end{array}
\label{eq:ccg_feas}
\end{equation}
The following proposition enables us to bound the optimality gap associated with a given feasible solution to the relaxed master problem.
\begin{proposition}
Let $({\bm x},{\bm w},\{{\bm y}^k\}_{k\in \sets K})$ be feasible in the relaxed master problem~\eqref{eq:ccg_master}. Then, $({\bm x},{\bm w},\{{\bm y}^k\}_{k\in \sets K})$ is feasible in Problem~\eqref{eq:min_max_regret_general_mo_Kadapt} and the objective value of $({\bm x},{\bm w},\{{\bm y}^k\}_{k\in \sets K})$ in Problem~\eqref{eq:min_max_regret_general_mo_Kadapt} is given by the optimal objective value of Problem~\eqref{eq:ccg_feas}.
\label{prop:algo_correct}
\end{proposition}
Proposition~\ref{prop:algo_correct} implies that, for any $({\bm x},{\bm w},\{{\bm y}^k\}_{k\in \sets K})$ feasible in the relaxed master problem~\eqref{eq:ccg_master}, the optimal value of~\eqref{eq:ccg_feas} yields an upper bound to the optimal value of the $K$-adaptability problem~\eqref{eq:min_max_regret_general_mo_Kadapt}. At the same time, it is evident that for any index set $\widetilde{\sets I} \subseteq \sets I^K$, the optimal value of Problem~\eqref{eq:ccg_master} yields a lower bound to the optimal objective value of Problem~\eqref{eq:min_max_regret_general_mo_Kadapt}. The lemma below is key to identify indices of subproblems ${\bm i} \in \sets I^K$ that are infeasible.
\begin{lemma}
Let $(\tau,{\bm x},{\bm w},\{{\bm y}^k\}_{k\in \sets K}, \{ {\bm \alpha}^{\bm i}, {\bm \beta}^{\bm i}\}_{{\bm i}\in \widetilde{\sets I}} , \{ {\bm \beta}^{{\bm i},k} ,  {\bm \gamma}^{{\bm i},k}\}_{{\bm i}\in \widetilde{\sets I}, k\in \sets K} )$ be optimal in the relaxed master problem~\eqref{eq:ccg_master}. Let $(\theta,\overline{\bm \xi}, \{{\bm \xi}^{k}\}_{k \in \sets K}, {\bm \eta}, \{{\bm \zeta}^{k}\}_{k \in \sets K})$ be optimal in Problem~\eqref{eq:ccg_feas}. Then, the following hold:
\begin{enumerate}[label=(\roman*)]
    \item $\theta \geq \tau$;
    \item If $\theta = \tau$, then Problem~\eqref{eq:ccg_subproblem} is feasible for all ${\bm i} \in \sets I^K$;
    \item If $\theta > \tau$, then the index~${\bm i}$ defined through
$$
{\bm i}_k := \sum_{i \in \sets I} i \cdot \mathbb I(  {\bm \zeta}_i^k = 1 ) \quad \forall  k \in \sets K
$$
corresponds to an infeasible subproblem, i.e., Problem~\eqref{eq:ccg_subproblem} is infeasible. 
\end{enumerate}
\label{lem:algo_correct2}
\end{lemma}

Propositions~\ref{prop:algo_correct} and Lemma~\ref{lem:algo_correct2} culminate in Algorithm~\ref{alg:ccg} whose convergence is guaranteed by the following theorem.
\begin{theorem}
Algorithm~\ref{alg:ccg} terminates in a final number of steps with a feasible solution to Problem~\eqref{eq:min_max_regret_general_mo_Kadapt}. The objective value~$\theta$ attained by this solution is within $\delta$ of the optimal objective value of the problem.
\label{thm:algo_converges}
\end{theorem}


\begin{algorithm}[t!]
\SetAlgoLined
\textbf{Inputs:} Optimality tolerance $\delta$; $K$-adaptability parameter $K$\; 
\textbf{Output:} Near optimal solution $({\bm x},{\bm w}, \{{\bm y}^{k}\}_{k\in \sets K})$ to Problem~\eqref{eq:min_max_regret_general_mo_Kadapt} with associated objective $\theta$\;
\textbf{Initialization:} 

 Initialize upper and lower bounds: ${\rm{LB}} \leftarrow -\infty$ and ${\rm{UB}} \leftarrow +\infty$\;
 Initialize index set: $\widetilde{\sets I} \leftarrow \{ \1 \}$\;
 
 \While{${\rm{UB}}-{\rm{LB}} > \delta$}{
  Solve the master problem~\eqref{eq:ccg_master}, let $(\tau,{\bm x},{\bm w},\{{\bm y}^k\}_{k\in \sets K}, \{ {\bm \alpha}^{\bm i}, {\bm \beta}^{\bm i}\}_{{\bm i}\in \widetilde{\sets I}} , \{ {\bm \beta}^{{\bm i},k} ,  {\bm \gamma}^{{\bm i},k}\}_{{\bm i}\in \widetilde{\sets I}, k\in \sets K} )$ be an optimal solution\;
  Let ${\rm{LB}} \leftarrow \tau$\;
  Solve the feasibility subproblem~\eqref{eq:ccg_feas}, let $(\theta,\overline{\bm \xi}, \{{\bm \xi}^{k}\}_{k \in \sets K}, {\bm \eta}, \{{\bm \zeta}^{k}\}_{k \in \sets K})$ denote an optimal solution\;
  Let ${\rm{UB}} \leftarrow \theta$\;
    \If{$\theta > \tau$}{
        ${\bm i}_k \leftarrow  \sum_{i \in \sets I} i \cdot \mathbb I(  {\bm \zeta}_i^k = 1 )$ for all $k\in \sets K$\;
        $\widetilde{\sets I} \leftarrow \widetilde{\sets I}  \cup \{ {\bm i}\}$\;
    }
   
 }
 \textbf{Result:} $({\bm x},{\bm w},\{{\bm y}^k\}_{k\in \sets K})$ is near-optimal in~\eqref{eq:min_max_regret_general_mo_Kadapt} with objective value $\theta$.
 \caption{``Column-and-Constraint'' Generation Procedure.}
 \label{alg:ccg}
\end{algorithm}


\subsection{Worst-Case Absolute Regret Minimization}
\label{sec:EC_min_max_regret}

In this section, we show that certain classes of two-stage robust optimization problems that seek to minimize the ``worst-case absolute regret'' criterion can be written in the form~\eqref{eq:min_max_regret_general_mo}. According to the ``worst-case absolute regret'' criterion, the performance of a decision is evaluated with respect to the worst-case regret that is experienced, when comparing the performance of the decision taken relative to the performance of the best decision that should have been taken \emph{in hindsight}, after all uncertain parameters are revealed, see e.g., \cite{Savage51:Theory}. The minimization of worst-case absolute regret is often believed to mitigate the conservatism of classical robust optimization and is thus attractive in practical applications, see also Section~\ref{sec:computational_studies} for corroborating evidence.

Mathematically, we are given a utility function
\begin{equation}
u({\bm x},{\bm w},{\bm y},{\bm \xi}) \; := \; {\bm \xi}^\top {\bm C} \; {\bm x} + {\bm \xi}^\top {\bm D} \; {\bm w} + {\bm \xi}^\top {\bm Q} \; {\bm y}
\label{eq:utility_function}
\end{equation}
for which high values are preferred. This function depends on both the decisions ${\bm x}$, ${\bm w}$, and ${\bm y}$, and on the uncertain parameters ${\bm \xi}$. Given a realization of ${\bm \xi}$, we can measure the \emph{absolute regret} of a decision $({\bm x},{\bm w},{\bm y})$ as the difference between the utility of the best decision in hindsight (i.e., after ${\bm \xi}$ becomes known) and the utility of the decision actually taken, i.e.,
\begin{equation*}
\left\{
         \max_{ {\bm x}', {\bm w}', {\bm y}' } \;\;\; u({\bm x}',{\bm w}',{\bm y}',{\bm \xi}) - 
         u({\bm x},{\bm w},{\bm y},{\bm \xi}) \; : \;
         {\bm x}' \in \sets X, \; {\bm w}' \in \sets W, \; {\bm y}' \in \sets Y
\right\}.
\end{equation*}
Regret averse decision-makers seek to minimize the worst-case (maximum) absolute regret
\begin{equation}
\max_{ {\bm \xi} \in \Xi({\bm w},\overline {\bm \xi}) } \; \;  \left\{
         \max_{ {\bm x}', {\bm w}', {\bm y}' } \;\;\; u({\bm x}',{\bm w}',{\bm y}',{\bm \xi}) - 
         u({\bm x},{\bm w},{\bm y},{\bm \xi}) \; : \;
         {\bm x}' \in \sets X, \; {\bm w}' \in \sets W, \; {\bm y}' \in \sets Y \right\}.
\end{equation}

A two-stage robust optimization problem with DDID in which the decision-maker seeks to minimize his worst-case absolute regret is then expressible as
\begin{equation}
\tag{$\mathcal{WCAR}$}
    \begin{array}{cl}
         \min & \;\; \displaystyle \max_{\overline {\bm \xi} \in \Xi} \;\;  \min_{ {\bm y} \in \sets Y } \; 
         \max_{ {\bm \xi} \in \Xi({\bm w},\overline {\bm \xi}) } \; \;  \left\{
         \max_{ {\bm x}', {\bm w}', {\bm y}' } \;\;\; u({\bm x}',{\bm w}',{\bm y}',{\bm \xi}) - 
         u({\bm x},{\bm w},{\bm y},{\bm \xi}) \; : \;
         {\bm x}' \in \sets X, \; {\bm w}' \in \sets W, \; {\bm y}' \in \sets Y \right\}   \\
         \st & \;\; {\bm x} \in \sets X, \; {\bm w} \in \sets W.
    \end{array}
\label{eq:min_max_regret_general}
\end{equation}
The following observation shows that under certain assumptions, Problem~\eqref{eq:min_max_regret_general} can be written in the form~\eqref{eq:min_max_regret_general_mo}.

\begin{observation}
Suppose that the utility function~$u$ in~\eqref{eq:utility_function} and the feasible sets $\mathcal X$, $\mathcal W$, and $\mathcal Y$ in Problem~\eqref{eq:min_max_regret_general} are such that:
\begin{enumerate}[label=(\roman*)]
    \item Either ${\bm C}={\bm 0}$ or $\sets X := \{ {\bm x}: \1^\top {\bm x}=1\}$, and
    \item Either ${\bm D}={\bm 0}$ or $\sets W := \{ {\bm w}: \1^\top {\bm w}=1\}$, and
    \item Either ${\bm Q}={\bm 0}$ or $\sets Y := \{ {\bm y} : \1^\top {\bm y}=1\}$.
\end{enumerate}
Then, Problem~\eqref{eq:min_max_regret_general} can be written in the form~\eqref{eq:min_max_regret_general_mo}.
\label{obs:min_max_regret_as_mo}
\end{observation}



In Section~\ref{sec:preference_elicitation_KAS}, we leverage Observation~\ref{obs:min_max_regret_as_mo} and use Theorems~\ref{thm:Kadapt_mo_MBLP} and~\ref{thm:algo_converges}, and Algorithm~\ref{alg:ccg} to solve an active preference learning problem that seeks to recommend kidney allocation policies with least possible worst-case regret.


\section{Companion to Section~\ref{sec:preference_elicitation_KAS}: Generating Candidate Policies}
\label{sec:EC_preference_elicitation_KAS}

Based on the analysis of~\cite{bertsimas2013fairness}, we considered scoring policies where the score obtained by patient~$p$ for organ~$o$ is given by
$$
\begin{array}{l}
\text{Score}(p,o) \; = \; \alpha \text{LYFT}(p,o) + g(\text{DT}(p)) + \gamma \text{CPRA}(p) + \delta_1 \mathbb{I} ( \text{Age}(p) \leq 2 ) +  \cdots \\
\qquad \qquad \qquad \qquad \qquad \cdots  + \delta_2 \mathbb{I} ( 2 < \text{Age}(p) \leq 10 )  + \delta_3 \mathbb{I} ( 10 < \text{Age}(p) \leq 18 )  +  \cdots \\
\qquad \qquad \qquad \qquad \qquad \cdots   + \delta_4 \mathbb I ( 18 < \text{Age}(p) \leq 35 ) +\delta_5 \mathbb I ( 35 < \text{Age}(p) \leq 50 ) +  \cdots \\
\qquad \qquad \qquad \qquad \qquad \cdots   + \delta_6 \mathbb I ( 50 < \text{Age}(p) \leq 65 )  + \delta_7 \mathbb{I} ( 65 < \text{Age}(p) ),
\end{array}
$$
where LYFT represents life years from transplant, DT corresponds to the patient dialysis time, CPRA denotes Calculated Panel Reactive Antibodies, i.e., the percentage of donors with whom a particular recipient would be incompatible, and $g$ is a piecewise linear function with breakpoints at 5 and 10 years given by
$$
g(\text{DT}) := 
\begin{cases}
g_1 \text{DT} & \text{if } \text{DT} \leq 5 \\
5 g_1 + g_2(  \text{DT}-5 )   & \text{if } 5< \text{DT} \leq 10 \\
5 g_1 + 5 g_2 + g_3(  \text{DT}-10 )   & \text{if } 10< \text{DT}.
\end{cases}
$$
We generate $I=20$ policies of the form above. The first policy considered prioritizes based on LYFT only, the second and third policies are the KAS and pre-KAS policies of OPTN, respectively, the fourth policy prioritizes based on dialysis time only, and the remaining policies are generated randomly. For the random policies, they each use LYFT, DT, CPRA, and Age with probability 0.75. If LYFT is used, $\alpha=1$. If CPRA is used, $\gamma$ is sampled uniformly from the set $\{0,0.1,0.2,\ldots,5\}$, and the same is true for $\delta_i$, $i=1,\ldots,7$ if Age is used. If DT is used, the parameters of the function $g$ are sampled uniformly from $\{0,0.25,0.5,\ldots,5\}$.

For each of these policies, we record $J=22$ outcomes from the KPSAM simulator: the number of kidney transplants, the number of discarded kidneys, the number of pediatric kidney transplants, the number of ABO (i.e., blood type) identical kidney transplants, the number of kidney waitlist deaths, the number of zero BDr (antigen) mismatch isolated kidney transplants, the number of one BDr mismatch isolated kidney transplants, the number of zero Dr mismatch isolated kidney transplants, the number of one Dr (antigen) mismatch isolated kidney transplants, the number of kidney transplants to white people, black people, and hispanics, the number of transplants to female candidates, the number of transplants to candidates in 4 age groups, the number of deaths of white, black, and hispanic people, the number of deaths of female candidates, and the average wait time.


\section{Speed-Up Strategies}
\label{sec:extensions}

This section proposes several strategies for speeding-up the solution of the $K$-adaptability counterpart of problems with exogenous and/or endogenous uncertainty.


\subsection{Symmetry Breaking Constraints}
\label{sec:symmetry_breaking_general}

The $K$-adaptability problem~\eqref{eq:Kadapt} presents a large amount of symmetry since indices of the candidate policies can be permuted to yield another, distinct, feasible solution with identical cost. This symmetry yields to significant slow down of the branch-and-bound procedure, see e.g.,~\cite{BertsimasBook_OptIntegers}, in particular as $K$ grows. Thus, we propose to eliminate the symmetry in the problem by introducing symmetry breaking constraints. Specifically, we constrain the candidate policies $\{ {\bm y}^k \}_{k \in \sets K}$ to be lexicographically decreasing. For this purpose, we introduce auxiliary binary variables ${\bm z}^{k,k+1} \in \{0,1\}^{N_y}$ for all $k \in \sets K \backslash \{ K \}$ such that ${\bm z}^{k,k+1}_i=1$ if and only if policies ${\bm y}^k$ and ${\bm y}^{k+1}$ differ in their $i$th component. These variables can be defined by means of a moderate number of linear inequality constraints, as follows
\begin{equation}
\left.
    \begin{array}{l}
    {\bm z}^{k,k+1}_i \; \leq \; {\bm y}_i^k + {\bm y}_i^{k+1} \\
    {\bm z}^{k,k+1}_i \; \leq \; 2- {\bm y}_i^k - {\bm y}_i^{k+1} \\
    {\bm z}^{k,k+1}_i \; \geq \; {\bm y}_i^{k} - {\bm y}_i^{k+1} \\
    {\bm z}^{k,k+1}_i \; \geq \;  {\bm y}_i^{k+1} - {\bm y}_i^{k}
    \end{array} \right\} \quad \forall i \in \sets I, \; k \in \sets K \backslash \{K\}.
\label{eq:symmetry_breaking_zvars_definition}
\end{equation}
%
%
The first set of constraints above ensures that if ${\bm y}_i^k = {\bm y}_i^{k+1}$, then ${\bm z}^{k,k+1}_i = 0$. Conversely, the second set of constraints guarantees that ${\bm z}^{k,k+1}_i = 1$ whenever ${\bm y}_i^k \neq {\bm y}_i^{k+1}$. Using the variables ${\bm z}^{k,k+1}_i$, the lexicographic ordering constraints can be written as
\begin{equation}
    {\bm y}_i^k \; \geq \; {\bm y}_i^{k+1} - \sum_{i' < i} {\bm z}_{i'}^{k,k+1} \quad \forall i \in \sets I, \; k \in \sets K \backslash \{K\}.
\label{eq:lexicographic}
\end{equation}
These stipulate that if ${\bm y}_{i'}^k = {\bm y}_{i'}^{k+1}$ for all $i' < i$, then ${\bm y}_i^k \geq {\bm y}_i^{k+1}$. Since the symmetry breaking constraints in~\eqref{eq:symmetry_breaking_zvars_definition} and~\eqref{eq:lexicographic} are deterministic, they can be added to the $K$-adaptability problem without affecting the solution procedure.

\subsection{Heuristic $K$-Adaptability Solution Approach}
\label{sec:heuristic_numericals}

To speed-up computation in our numerical experiments to be able to investigate the performance of our approach for high values of $K$ (up to $\sim 10$), we employ a conservative solution approach, as detailed in Algorithm~\ref{alg:heuristic}. A variant of this approach has been previously used by~\cite{subramanyam2017k}. This algorithm returns a feasible but potentially suboptimal solution to the $K$-adaptability counterpart of the problem to be solved.

\begin{algorithm}[t!]
\SetAlgoLined
\textbf{Inputs:} 
    Instance of Problem~\eqref{eq:endo_2_obj},~\eqref{eq:endo_3}, or~\eqref{eq:min_max_regret_general_mo}; $K$-adaptability parameter $K$\; 
\textbf{Output:} Conservative solution $({\bm x},{\bm w}, \{{\bm y}^{k}\}_{k\in \sets K})$ to the $K$-adaptability counterpart of the input instance (\eqref{eq:Kadapt_obj},~\eqref{eq:Kadapt}, or~\eqref{eq:min_max_regret_general_mo_Kadapt}, respectively)\;

 \For{$k \in \{1,\ldots,K\}$}{
 
    \eIf{$k=1$}
    {
        Solve the the $k$-adaptability counterpart of the input instance (using its MBLP reformulation)\;
        Let $({\bm x}^\star,{\bm w}^\star,{\bm y}^{\star,1})$ denote an optimal solution\;
    }
    {
        Solve the $k$-adaptability counterpart of the input instance (using its MBLP reformulation) with the added constraints that ${\bm y}^\kappa={\bm y}^{\star,\kappa}$ for all $\kappa \in \{1,\ldots,k-1\}$\;
        Let $({\bm x}^\star,{\bm w}^\star,\{{\bm y}^{\star,\kappa}\}_{\kappa=1}^{k})$ denote an optimal solution\;
    }
 }
 \textbf{Result:} Return $({\bm x}^\star,{\bm w}^\star,\{{\bm y}^{\star,\kappa}\}_{\kappa\in \sets K})$.
 \caption{Heuristic algorithm for solving the $K$-adaptability counterpart of a problem; adapted from~\cite{subramanyam2017k}.}
 \label{alg:heuristic}
\end{algorithm}


\section{Companion to Section~\ref{sec:problem_formulation_endo}: Unattained Optimal Value}
\label{sec:EC_optval_notattained}


\begin{figure}
    \centering
    \includegraphics[height=0.205\textheight]{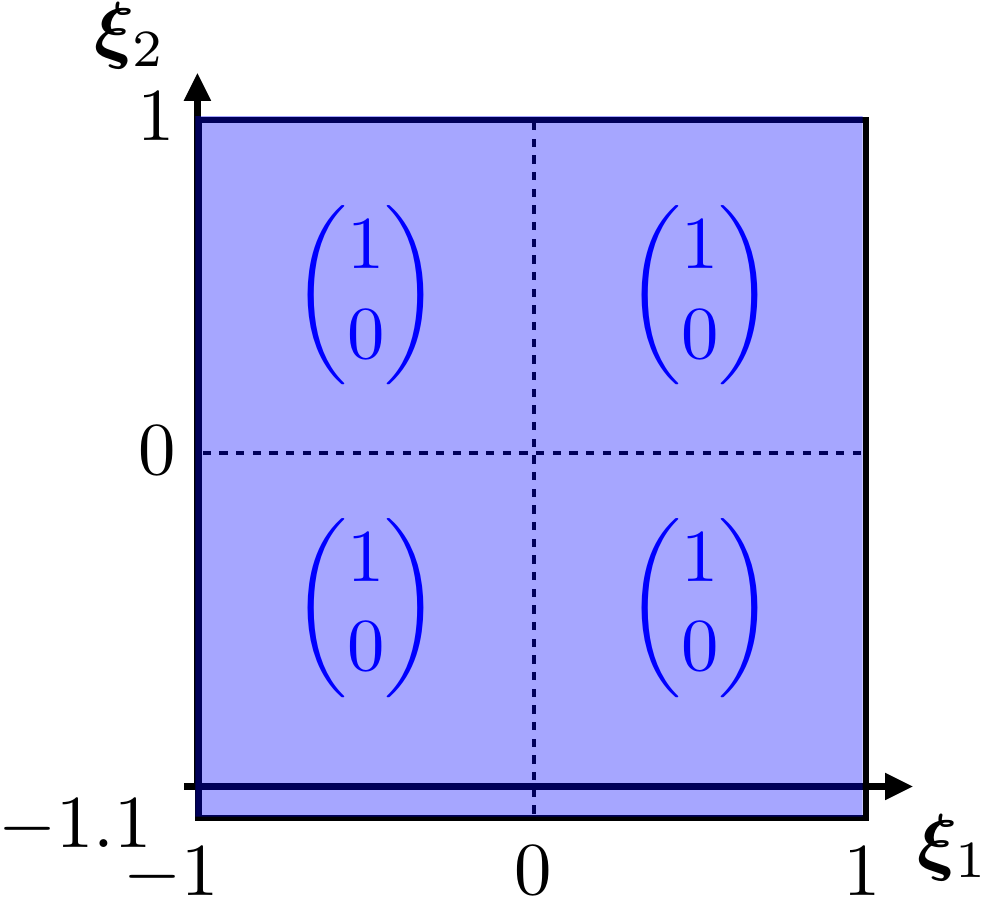} \hspace{1cm}
    \includegraphics[height=0.205\textheight]{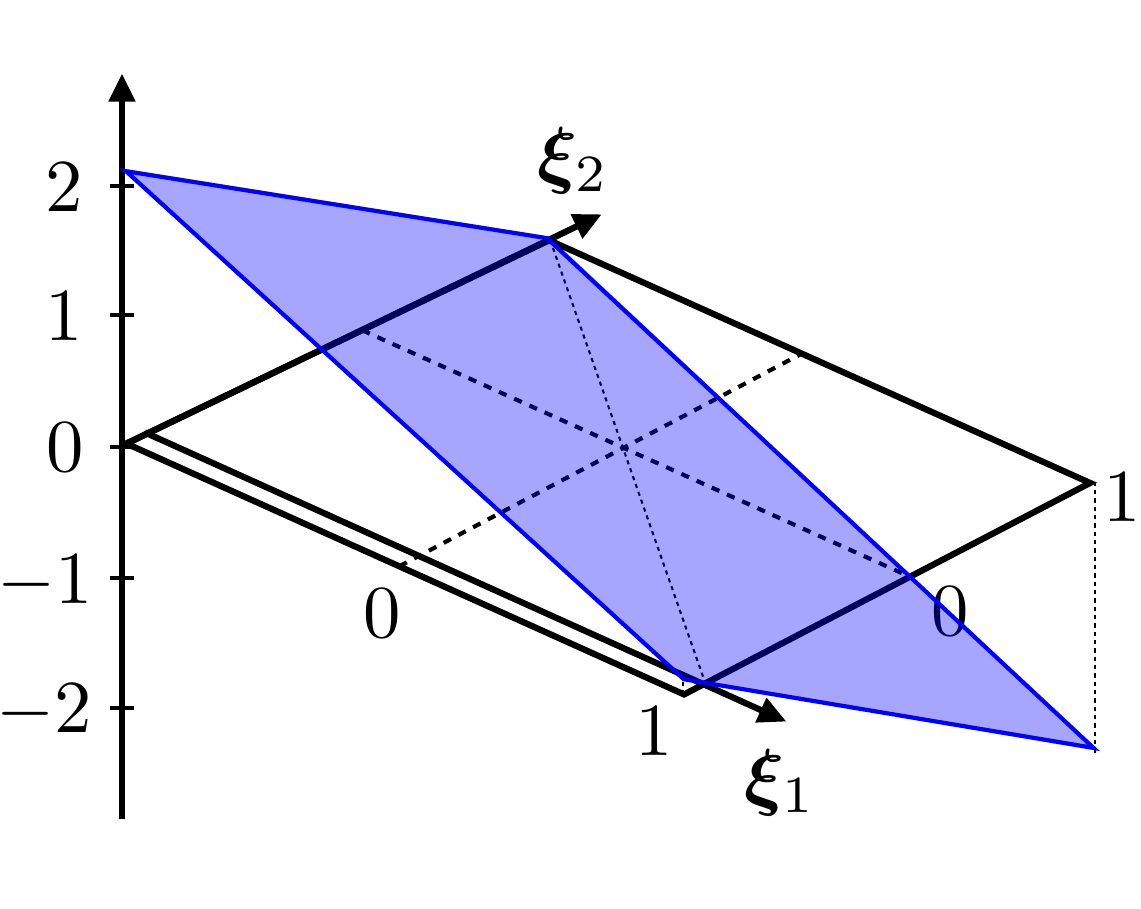} \\
    \includegraphics[height=0.205\textheight]{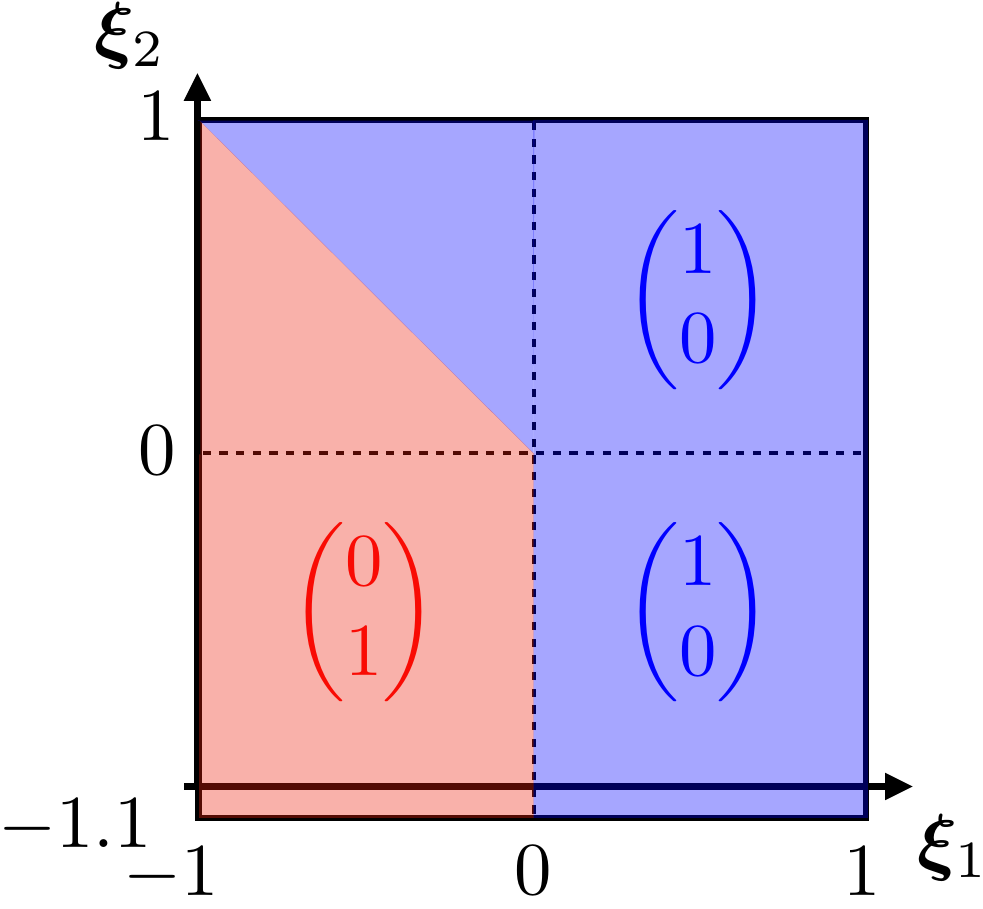} \hspace{1cm}
    \includegraphics[height=0.205\textheight]{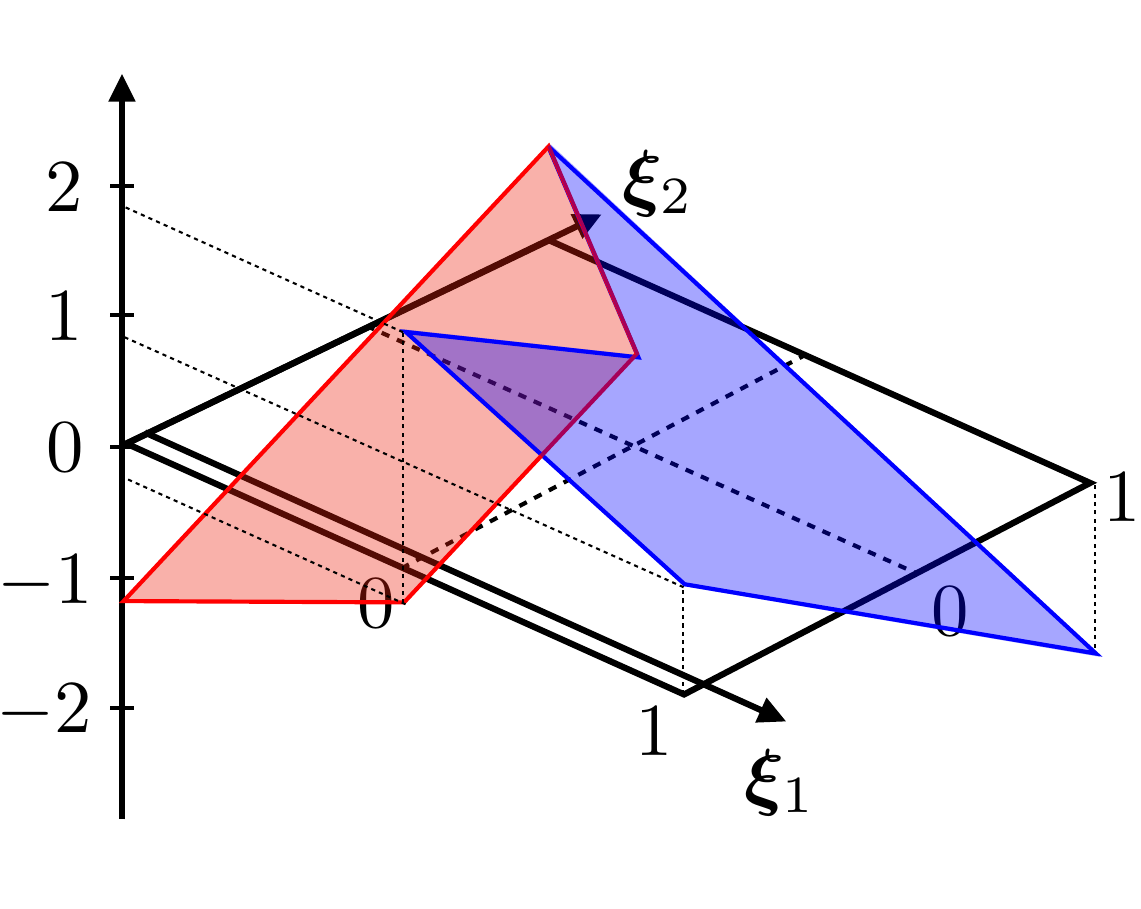} \\
    \includegraphics[height=0.205\textheight]{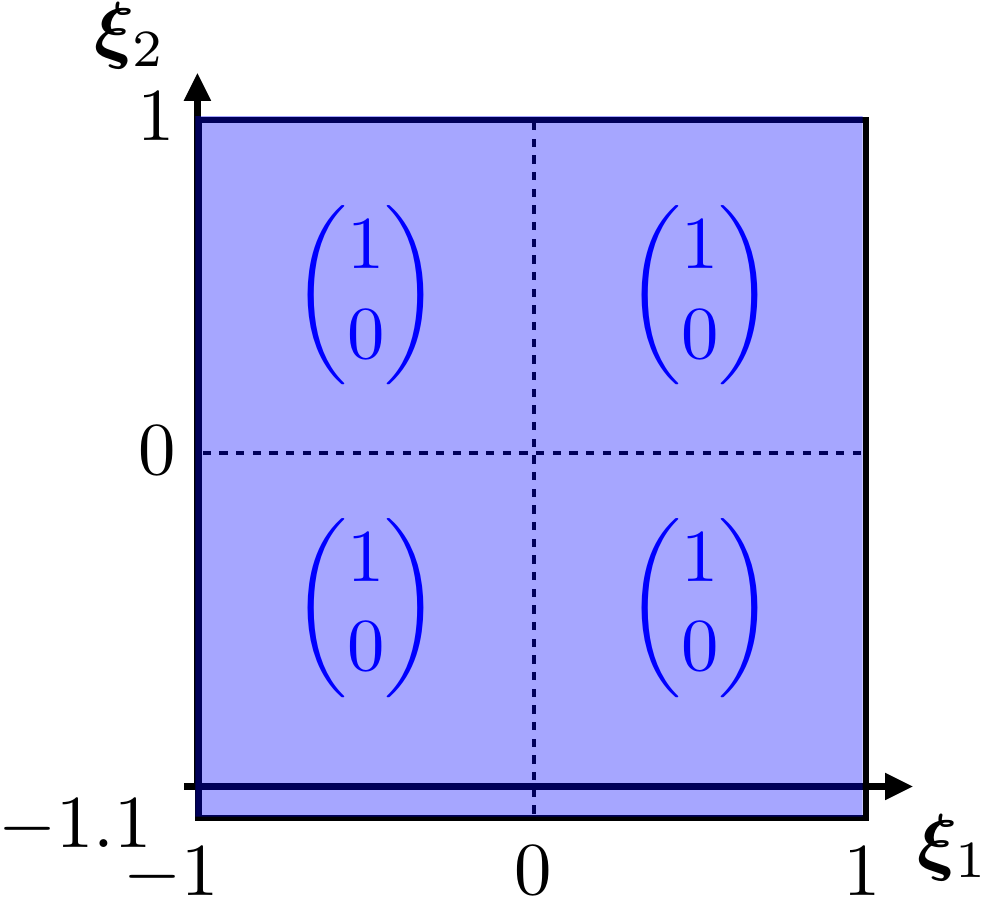} \hspace{1cm}
    \includegraphics[height=0.205\textheight]{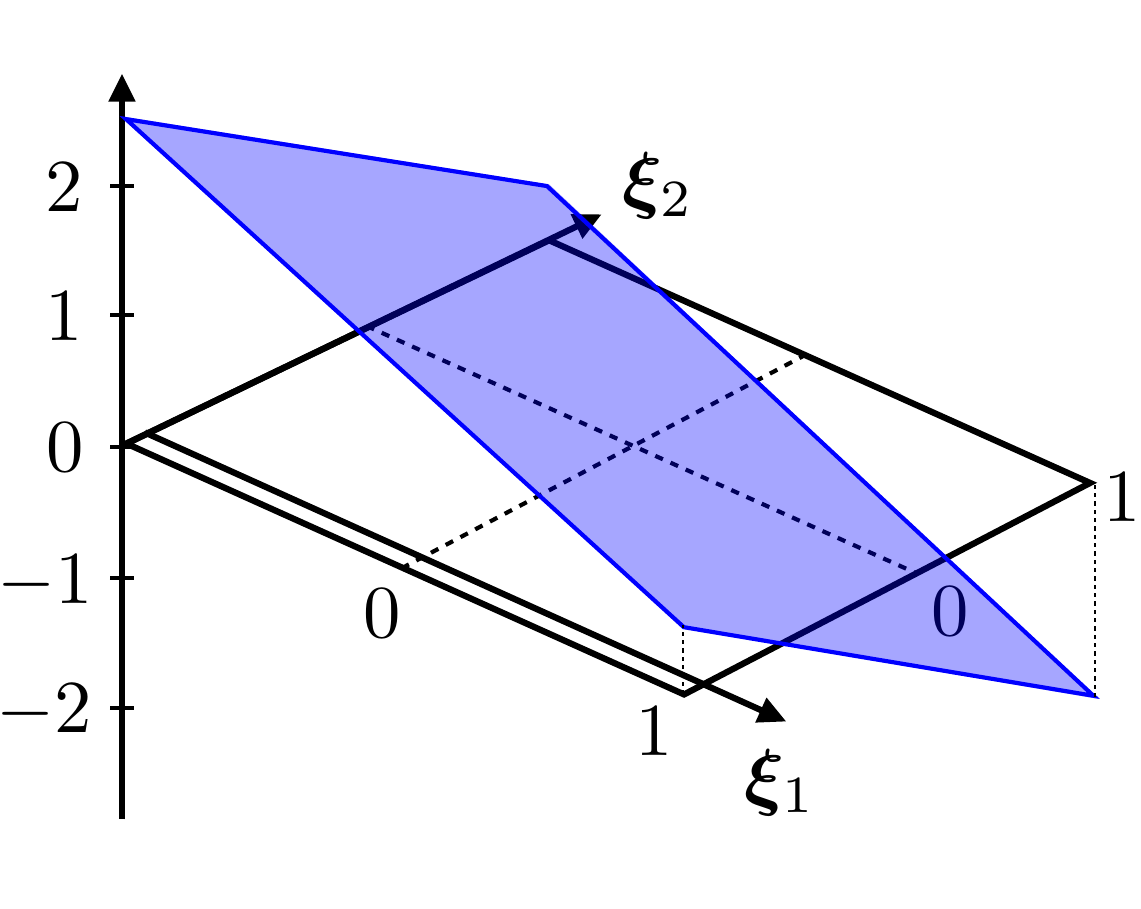} \\
    \includegraphics[height=0.205\textheight]{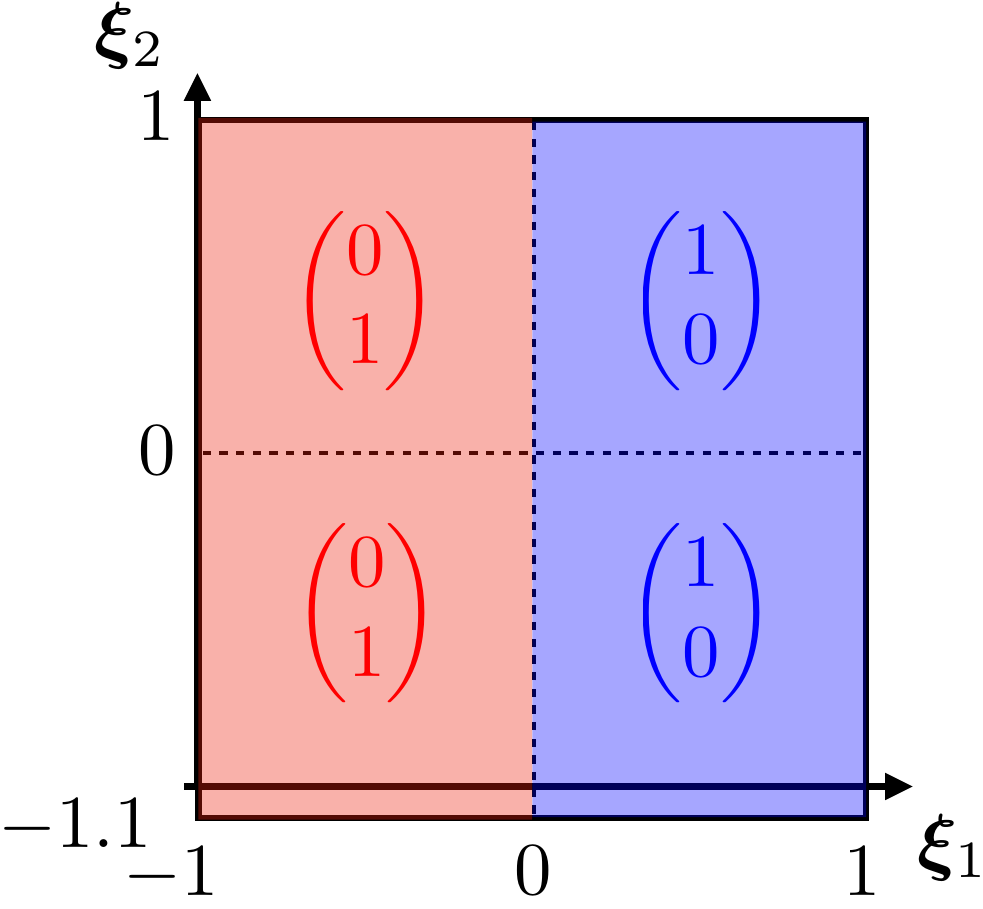} \hspace{1cm}
    \includegraphics[height=0.205\textheight]{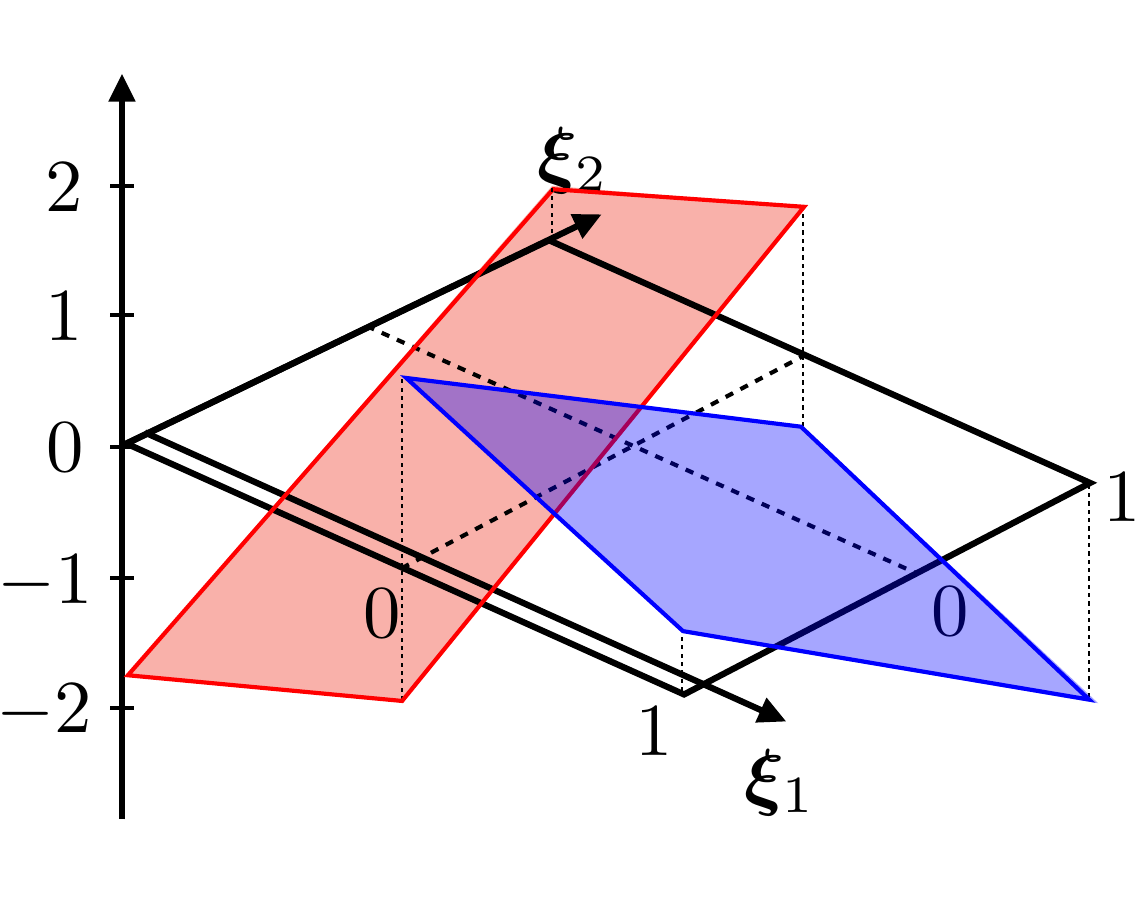}
    \caption{Companion figure for the example in Electronic Companion~\ref{sec:EC_optval_notattained}, assuming ${\bm d}_1={\bm d}_2=0.4$. Optimal wait-and-see decision ${\bm y}$ (left) and associated objective function (right) for the cases when ${\bm w}={\bm 0}$ (first row), ${\bm w}=\1$ (second row), ${\bm w}=(0,1)$ (third row), and ${\bm w}=(1,0)$ (last row) in Problem~\eqref{eq:example_general}. The optimal solution is given by ${\bm w}^\star = (1,0)$. For the optimal solution ${\bm w}^\star$, the objective function is discontinuous on the set $\{ {\bm \xi} \in \Xi \; : \; {\bm \xi}_1=0$\} and in particular the optimal objective value is not attained.}
    \label{fig:example_general}
\end{figure}


In this section, we show that the optimal value of Problem~\eqref{eq:endo_3} is in general not attained. Consider the following instance of Problem~\eqref{eq:endo_3}, adapted from~\cite{Hanasusanto2015} to incorporate decision-dependent information discovery.
\begin{equation}
\begin{array}{cccl}
\minimize_{{\bm w} \in \{0,1\}^ 2} & \;\; \displaystyle \max_{ \overline{\bm \xi} \in \Xi } & \;\; \displaystyle \min_{ {\bm y} \in \{0,1\}^2 } \;\; & \;\; \displaystyle \max_{ {\bm \xi} \in \Xi({\bm w},\overline{\bm \xi}) } \;\; ( {\bm \xi}_1 + {\bm \xi}_2) ({\bm y}_2 - {\bm y}_1 ) + {\bm d}_1 {\bm w}_1 + {\bm d}_2 {\bm w}_2 \\
&&\;\;\; \st \;\; & \;\; {\bm y}_1 \geq {\bm \xi}_1 \quad \forall {\bm \xi} \in \Xi({\bm w},\overline{\bm \xi}) \\
&&& \;\; {\bm y}_1 + {\bm y}_2 = 1,
\end{array}
\label{eq:example_general}
\end{equation}
where ${\bm d}_1$, ${\bm d}_2 \in (0,1)$ are given scalars representing the observation costs associated with ${\bm w}_1$ and ${\bm w}_2$, respectively, and $\Xi:= \{ {\bm \xi} \in \reals^2 \; : \; -1 \leq {\bm \xi}_1 \leq 1 ,\; -1.1 \leq {\bm \xi}_2 \leq 1 \}$. For each feasible choice of ${\bm w}$, we investigate the associated optimal wait-and-see decision, as well as the corresponding objective function value, see Figure~\ref{fig:example_general}.

Consider the choice ${\bm w}={\bm 0}$, whereby no uncertain parameter is observed between the first and second decision stages. Then, $\Xi( {\bm w},\overline{\bm \xi}) = \Xi$. Under this here-and-now decision, Problem~\eqref{eq:example_general} is expressible as a single-stage robust problem as follows
$$
\displaystyle \left\{ \min_{ {\bm y} \in \{0,1\}^2 } \;\; \left(\displaystyle \max_{ {\bm \xi} \in \Xi } \;\; ( {\bm \xi}_1 + {\bm \xi}_2) ({\bm y}_2 - {\bm y}_1 ) \right) \;\; : \;\;
 {\bm y}_1 \geq {\bm \xi}_1 \;\; \forall {\bm \xi} \in \Xi, \;\; {\bm y}_1 + {\bm y}_2 = 1 \right\}.
$$
It can can be readily verified that the only feasible (and therefore optimal) wait-and-see action in this case is ${\bm y}=( 1 , 0 )$, a {static} decision. The associated objective function and corresponding value is
$$
\max_{ {\bm \xi} \in \Xi } \;\; -( {\bm \xi}_1 + {\bm \xi}_2) \; = \; 2.1.
$$

Consider the choice ${\bm w} = \1$, whereby both uncertain parameters are observed between the first and second decision stages. Then, $\Xi( {\bm w},\overline{\bm \xi} ) = \{ \overline{\bm \xi} \}$. Under this here-and-now decision, Problem~\eqref{eq:example_general} reduces to
$$
\displaystyle \max_{ {\bm \xi} \in \Xi } \;\; \left\{ \;\; \displaystyle \min_{ {\bm y} \in \{0,1\}^2 } \;\;  \;\; ( {\bm \xi}_1 + {\bm \xi}_2) ({\bm y}_2 - {\bm y}_1 ) + {\bm d}_1 + {\bm d}_2 
\;\;\; : \;\;  \;\; {\bm y}_1 \geq {\bm \xi}_1 , \;\; {\bm y}_1 + {\bm y}_2 = 1 \right\}.
$$
The constraints in the problem imply that ${\bm y}=(1,0)$ is the only feasible (and therefore optimal) solution whenever ${\bm \xi}_1 > 0$. For ${\bm \xi}_1 \leq 0$, the optimal choices are ${\bm y}=(1,0)$ if ${\bm \xi}_1 + {\bm \xi}_2 \geq 0$, and ${\bm y} = (0,1)$, else. The associated objective function is
$$
 {\bm d}_1 + {\bm d}_2 + \max_{{\bm \xi} \in \Xi} \begin{cases}
- ({\bm \xi}_1 + {\bm \xi}_2) & \text{if } ( {\bm \xi}_1 > 0 ) \text{ or } ({\bm \xi}_1 \leq 0 \text{ and } {\bm \xi}_1 + {\bm \xi}_2 \geq 0 ) \\
\;\;\;\; {\bm \xi}_1 + {\bm \xi}_2 & \text{else,} 
\end{cases}
$$
yielding an objective value of $(1.1+{\bm d}_1 + {\bm d}_2)$ that is not attained.

Consider the choice ${\bm w}=(0,1)$, whereby only ${\bm \xi}_2$ is observed between the first and second decision stages. Then, Problem~\eqref{eq:example_general} reduces to
$$
\displaystyle \max_{ {\bm \xi}_2 \in [-1.1,1] } \;\; \left\{ \displaystyle \min_{ {\bm y} \in \{0,1\}^2 } \;\;  \;\; \displaystyle \max_{ {\bm \xi}_1 \in [-1,1] } \;\; ( {\bm \xi}_1 + {\bm \xi}_2) ({\bm y}_2 - {\bm y}_1 ) + {\bm d}_2  
\;\;\; : \;\;  {\bm y}_1 \geq {\bm \xi}_1 \quad \forall {\bm \xi}_1 \in [ -1 , 1] , 
 \;\; {\bm y}_1 + {\bm y}_2 = 1 \right\}.
$$
For any choice of ${\bm \xi}_2$, the only option for the wait-and-see decision is ${\bm y}=(1,0)$ (since ${\bm \xi}_1$ remains uncertain). The associated objective function and corresponding objective value is
$$
{\bm d}_2 + \max \;\;\left\{  -( {\bm \xi}_1 + {\bm \xi}_2) \; : \; {\bm \xi}_1 \in [-1,1], \; {\bm \xi}_2  \in [-1.1,1] \right\}\; = \; {\bm d}_2 + 2.1. 
$$

Lastly, consider the choice ${\bm w}=(1,0)$, whereby only ${\bm \xi}_1$ is observed between the first and second decision stages. Then, Problem~\eqref{eq:example_general} reduces to
$$
\displaystyle \max_{ {\bm \xi}_1 \in [-1,1] } \;\; \displaystyle \left\{ \min_{ {\bm y} \in \{0,1\}^2 } \;\; \;\; \displaystyle \max_{ {\bm \xi}_2 \in [-1.1,1] } \;\; ( {\bm \xi}_1 + {\bm \xi}_2) ({\bm y}_2 - {\bm y}_1 ) + {\bm d}_1  
\;\;\; : \;\;  \;\; {\bm y}_1 \geq {\bm \xi}_1,
 \;\; {\bm y}_1 + {\bm y}_2 = 1 \right\}.
$$
For ${\bm \xi}_1 >0$, the only feasible (and therefore optimal) choice is ${\bm y}=(1,0)$. If ${\bm \xi}_1 \leq 0$, then the optimal wait-and-see decision is ${\bm y}=(0,1)$. The associated objective function is
$$
 {\bm d}_1 + \max_{{\bm \xi} \in \Xi} \begin{cases}
- ({\bm \xi}_1 + {\bm \xi}_2) & \text{if } {\bm \xi}_1 > 0  \\
\;\;\;\; {\bm \xi}_1 + {\bm \xi}_2 & \text{else,} 
\end{cases}
$$
yielding an objective value of $(1.1+{\bm d}_1)$ that is not attained by any feasible solution.

We conclude that, since ${\bm d}_1$, ${\bm d}_2 \in (0,1)$, the optimal solution to Problem~\eqref{eq:example_general} is ${\bm w}^\star=(1,0)$ with associated optimal objective value $(1.1+{\bm d}_1)$ which is never attained.

\section{Proofs of Statements in Sections~\ref{sec:problem_formulation_exo} and~\ref{sec:problem_formulation_endo}}
\label{sec:EC_problem_formulation}

\proof{Proof of Theorem~\ref{thm:endo_equiv}}
Let $({\bm x},{\bm w})$, ${\bm y}'(\cdot)$, and ${\bm y}(\cdot)$ be defined as in the premise of claim \textit{(i)}. 
Then, ${\bm x} \in \sets X$, ${\bm w} \in \sets W$, and for each ${\bm \delta}$ such that ${\bm \delta} = {\bm w} \circ \overline{\bm \xi}$ for some $\overline{\bm \xi} \in \Xi$, we have that ${\bm y}'({\bm \delta}) \in \sets Y$ and ${\bm T} {\bm x} + {\bm V} {\bm w} + {\bm W} {\bm y}'({\bm \delta}) \leq {\bm H}{\bm \xi}$ for all ${\bm \xi} \in \Xi({\bm w},{\bm \delta})$. 
We show that $({\bm x},{\bm w},{\bm y}(\cdot))$ is feasible in Problem~\eqref{eq:endo_1}. Fix any ${\bm \xi} \in \Xi$. First, ${\bm y}({\bm \xi}) \in \sets Y$. Second, we have
$$
{\bm T} {\bm x} + {\bm V} {\bm w} + {\bm W} {\bm y}({\bm \xi}) 
\; = \; {\bm T} {\bm x} + {\bm V} {\bm w} + {\bm W} {\bm y}'({\bm w} \circ {\bm \xi}) 
\; \leq \; {\bm H}{\bm \xi},
$$
where the equality follows by definition of ${\bm y}(\cdot)$ and the inequality follows from the fact that ${\bm \xi} \in \Xi({\bm w},{\bm w} \circ {\bm \xi})$ and from the definition of ${\bm y}'(\cdot)$. Fix ${\bm \xi}' \in \Xi : {\bm w} \circ {\bm \xi} = {\bm w} \circ {\bm \xi}'$. Then, ${\bm y}({\bm \xi}) = {\bm y}({\bm \xi}')$, so that the decision-dependent non-anticipativity constraints are also satisfied. Since the choice of ${\bm \xi} \in \Xi$ was arbitrary, $({\bm x},{\bm w},{\bm y}(\cdot))$ is feasible in Problem~\eqref{eq:endo_1}. The objective value attained by $({\bm x},{\bm w})$ in Problem~\eqref{eq:endo_3} is given by
\begin{equation}
\displaystyle \max_{\begin{smallmatrix} \overline{\bm \xi} \in \Xi,\\ {\bm \xi} \in \Xi({\bm w},\overline{\bm \xi}) \end{smallmatrix}} \;\; \;\; {\bm \xi}^\top {\bm C} \; {\bm x} +  {\bm \xi}^\top {\bm D} \; {\bm w} + {\bm \xi}^\top {\bm Q} \; {\bm y}'({\bm w} \circ \overline {\bm \xi})
\; = \; 
\displaystyle \max_{\begin{smallmatrix} \overline{\bm \xi} \in \Xi,\\ {\bm \xi} \in \Xi({\bm w},\overline{\bm \xi}) \end{smallmatrix}} \;\; \;\; {\bm \xi}^\top {\bm C} \; {\bm x} +  {\bm \xi}^\top {\bm D} \; {\bm w} + {\bm \xi}^\top {\bm Q} \; {\bm y}({\bm \xi}),
\label{eq:thm2_obj1}
\end{equation}
where we have grouped the two maximization problems in a single one and where the equality follows from the definition of ${\bm y}(\cdot)$. The value attained by $({\bm x},{\bm w},{\bm y}(\cdot))$ in Problem~\eqref{eq:endo_1} is
\begin{equation}
\displaystyle \max_{{\bm \xi} \in \Xi } \;\; \;\; {\bm \xi}^\top {\bm C} \; {\bm x} +  {\bm \xi}^\top {\bm D} \; {\bm w} + {\bm \xi}^\top {\bm Q} \; {\bm y}({\bm \xi}).
\label{eq:thm2_obj2}
\end{equation}
Since $
 \left\{ {\bm \xi} \in \Xi({\bm w},\overline {\bm \xi}) \; : \; \overline {\bm \xi} \in \Xi \right\} \; = \;  \Xi
 $, it follows that the optimal objective values of the Problems~\eqref{eq:thm2_obj1} and~\eqref{eq:thm2_obj2} are equal. We have thus shown that Problem~\eqref{eq:endo_1} lower bounds Problem~\eqref{eq:endo_3} and that if $({\bm x},{\bm w})$ is optimal in Problem~\eqref{eq:endo_3}, then the triple $({\bm x},{\bm w},{\bm y}(\cdot))$ is feasible in Problem~\eqref{eq:endo_1} with the two solutions attaining the same cost in their respective problems.

Next, let $({\bm x},{\bm w},{\bm y}(\cdot))$ be defined as in the premise of claim \textit{(ii)}, i.e., let it be optimal in Problem~\eqref{eq:endo_1}. The here-and-now decision $({\bm x},{\bm w})$ is feasible in Problem~\eqref{eq:endo_3} and, for each $\overline{\bm \xi} \in \Xi$, we can define
$$
 {\bm y}'( \overline{\bm \xi} ) \; \in \; \argmin_{ {\bm y} \in \sets Y } \; \; \left\{ \max_{ {\bm \xi} \in \Xi({\bm w},\overline{\bm \xi}) } \; \; {\bm \xi}^\top {\bm C} \; {\bm x} + {\bm \xi}^\top {\bm D} \; {\bm w} + {\bm \xi}^\top {\bm Q} \; {\bm y} \; : \; {\bm T} {\bm x}  + {\bm V}{\bm w} + {\bm W}{\bm y} \leq {\bm H}{\bm \xi}  \; \; \; \; \forall {\bm \xi} \in \Xi({\bm w},\overline{\bm \xi}) \right\}.
$$
By construction, $({\bm x},{\bm w}) \in \sets X \times \sets W$. Moreover, it holds that 
$$
\begin{array}{cl}
& \quad \displaystyle \max_{\overline{\bm \xi} \in \Xi} \;\; \displaystyle \max_{{\bm \xi} \in \Xi({\bm w},\overline{\bm \xi})} \;\; \;\; {\bm \xi}^\top {\bm C} \; {\bm x} +  {\bm \xi}^\top {\bm D} \; {\bm w} + {\bm \xi}^\top {\bm Q} \; {\bm y}'(\overline{\bm \xi}) \\
 =  & \quad 
\displaystyle \max_{{\bm \xi} \in \Xi} \;\; \;\; {\bm \xi}^\top {\bm C} \; {\bm x} +  {\bm \xi}^\top {\bm D} \; {\bm w} + {\bm \xi}^\top {\bm Q} \; {\bm y}'({\bm \xi}) \\
 \leq & \quad
\displaystyle \max_{{\bm \xi} \in \Xi}  \;\; \;\; {\bm \xi}^\top {\bm C} \; {\bm x} +  {\bm \xi}^\top {\bm D} \; {\bm w} + {\bm \xi}^\top {\bm Q} \; {\bm y}({\bm \xi}).
\end{array}
$$
Thus, $({\bm x},{\bm w})$ is feasible in Problem~\eqref{eq:endo_3} with a cost no greater than that of $({\bm x},{\bm w},{\bm y}(\cdot))$ in Problem~\eqref{eq:endo_1}. We have thus shown that Problem~\eqref{eq:endo_3} lower bounds Problem~\eqref{eq:endo_1} and that if $({\bm x},{\bm w},{\bm y}(\cdot))$ is optimal in Problem~\eqref{eq:endo_1}, then $({\bm x},{\bm w})$ is feasible in Problem~\eqref{eq:endo_3} with the cost attained by ${\bm x}$ in Problem~\eqref{eq:endo_3} being no greater than the cost of $({\bm x},{\bm y}(\cdot))$ in Problem~\eqref{eq:endo_1}.

We conclude that the optimal costs of Problems~\eqref{eq:endo_1} and~\eqref{eq:endo_3} are equal, and that claims \textit{(i)} and \textit{(ii)} hold.~\Halmos 
\endproof

\section{Proofs of Statements in Section~\ref{sec:kadaptability_overview}}
\label{sec:EC_kadaptability_overview}

\proof{Proof of Lemma~\ref{lem:Kadapt_cstr_min_max_min}} Fix ${\bm x} \in \sets X$, ${\bm w} \in \sets W$, and ${\bm y}^k$, $k \in \sets K$, and $\overline {\bm \xi} \in \Xi$. It suffices to show that the problems
\begin{equation}
\min_{ k \in \sets K } \;  \left\{ \max_{ {\bm \xi} \in \Xi({\bm w},\overline{\bm \xi}) } \; \; {\bm \xi}^\top {\bm C} \; {\bm x} + {\bm \xi}^\top {\bm D} \; {\bm w} + {\bm \xi}^\top {\bm Q} \; {\bm y}^k \; : \; {\bm T} {\bm x} + {\bm V} {\bm w} + {\bm W}{\bm y}^k \leq {\bm H}{\bm \xi}  \; \; \; \; \forall {\bm \xi} \in \Xi({\bm w},\overline{\bm \xi}) \right\}
\label{eq:Kadapt_min_max_min_sp1}
\end{equation}
and
\begin{equation}
\max_{   {\bm \xi}^k \in \Xi({\bm w},\overline {\bm \xi}), \; k \in \sets K} \;\;\min_{ k \in \sets K } \;  \left\{  ({\bm \xi}^k)^\top {\bm C} \; {\bm x} + ({\bm \xi}^k)^\top {\bm D} \; {\bm w} + ({\bm \xi}^k)^\top {\bm Q} \; {\bm y}^k \; : \; {\bm T} {\bm x} + {\bm V} {\bm w} + {\bm W}{\bm y}^k \leq {\bm H}{\bm \xi}^k \right\}
\label{eq:Kadapt_min_max_min_sp2}
\end{equation}
have the same optimal objective.

Problem~\eqref{eq:Kadapt_min_max_min_sp1} is either infeasible or has a finite objective value. Indeed, it cannot be unbounded below since, if it is feasible, its objective value is given as the minimum of finitely many terms each of which is bounded, by virtue of the compactness of the non-empty set $\Xi({\bm w},\overline{\bm \xi})$. Similarly, Problem~\eqref{eq:Kadapt_min_max_min_sp2} is either unbounded above or has a finite objective value. It cannot be infeasible since $\Xi({\bm w},\overline{\bm \xi})$ is non-empty.

We proceed in two steps. First, we show that Problem~\eqref{eq:Kadapt_min_max_min_sp1} is infeasible if and only if Problem~\eqref{eq:Kadapt_min_max_min_sp2} is unbounded above, in which case both problems have an optimal objective value of $+\infty$. Second, we show that if the problems have a finite optimal objective value, then their optimal values are equal.

For the first claim, we have
$$
\begin{array}{cl}
& \text{Problem~\eqref{eq:Kadapt_min_max_min_sp1} is infeasible} \\
\Leftrightarrow & \nexists k \in \sets K \; : \; {\bm T} {\bm x} + {\bm V} {\bm w} + {\bm W}{\bm y}^k \leq {\bm H}{\bm \xi}  \; \; \; \; \forall {\bm \xi} \in \Xi({\bm w},\overline{\bm \xi})\\
\Leftrightarrow & \forall k \in \sets K, \; \exists \tilde{\bm \xi}^k \in \Xi({\bm w},\overline{\bm \xi}) \; : \; {\bm T} {\bm x} + {\bm V} {\bm w} + {\bm W}{\bm y}^k \nleq {\bm H}\tilde{\bm \xi}^k\\
\Leftrightarrow & \text{Problem~\eqref{eq:Kadapt_min_max_min_sp2} is unbounded}.
\end{array}
$$

For the second claim, we proceed in two steps. First, we show that the optimal objective value of Problem~\eqref{eq:Kadapt_min_max_min_sp2} can be no greater than the optimal objective value of Problem~\eqref{eq:Kadapt_min_max_min_sp1}. Then, we show that the converse is also true.

For the first part, let $\tilde k$ be feasible in Problem~\eqref{eq:Kadapt_min_max_min_sp1} and $\{ \tilde{\bm \xi}^k \}_{k \in \sets K}$ be feasible in Problem~\eqref{eq:Kadapt_min_max_min_sp2}. The objective value attained by $\tilde k$ in Problem~\eqref{eq:Kadapt_min_max_min_sp1} is given by
$$
\max_{ {\bm \xi} \in \Xi({\bm w},\overline{\bm \xi}) } \; \; {\bm \xi}^\top {\bm C} \; {\bm x} + {\bm \xi}^\top {\bm D} \; {\bm w} + {\bm \xi}^\top {\bm Q} \; {\bm y}^{\tilde k}.
$$
Accordingly, the objective value attained by $\{ \tilde{\bm \xi}^k \}_{k \in \sets K}$ in Problem~\eqref{eq:Kadapt_min_max_min_sp2} is given by
$$
\min_{ k \in \sets K } \;  \left\{  (\tilde{\bm \xi}^k)^\top {\bm C} \; {\bm x} + (\tilde{\bm \xi}^k)^\top {\bm D} \; {\bm w} + (\tilde{\bm \xi}^k)^\top {\bm Q} \; {\bm y}^k \; : \; {\bm T} {\bm x} + {\bm V} {\bm w} + {\bm W}{\bm y}^k \leq {\bm H}\tilde{\bm \xi}^k \right\}.
$$
Next, note that
$$
\begin{array}{cl}
& \displaystyle \min_{ k \in \sets K } \;  \left\{  (\tilde{\bm \xi}^k)^\top {\bm C} \; {\bm x} + (\tilde{\bm \xi}^k)^\top {\bm D} \; {\bm w} + (\tilde{\bm \xi}^k)^\top {\bm Q} \; {\bm y}^k \; : \; {\bm T} {\bm x} + {\bm V} {\bm w} + {\bm W}{\bm y}^k \leq {\bm H}\tilde{\bm \xi}^k \right\} \\
\leq &  (\tilde{\bm \xi}^{\tilde k})^\top {\bm C} \; {\bm x} + (\tilde{\bm \xi}^{\tilde k})^\top {\bm D} \; {\bm w} + (\tilde{\bm \xi}^{\tilde k})^\top {\bm Q} \; {\bm y}^{\tilde k} \\
\leq & \displaystyle \max_{ {\bm \xi} \in \Xi({\bm w},\overline{\bm \xi}) } \; \; {\bm \xi}^\top {\bm C} \; {\bm x} + {\bm \xi}^\top {\bm D} \; {\bm w} + {\bm \xi}^\top {\bm Q} \; {\bm y}^{\tilde k},
\end{array}
$$
where the first inequality follows by feasibility of $\tilde k$ in Problem~\eqref{eq:Kadapt_min_max_min_sp1} since $\tilde{\bm \xi}^{\tilde k} \in \Xi({\bm w},\overline{\bm \xi})$ and the second inequality follows by feasibility of $\tilde{\bm \xi}^{\tilde k}$ in the maximization problem. Since the choices of $\tilde k$ and $\{ \tilde{\bm \xi}^k \}_{k \in \sets K}$ were arbitrary, it follows that the optimal objective of Problem~\eqref{eq:Kadapt_min_max_min_sp1} upper bounds the optimal objective of Problem~\eqref{eq:Kadapt_min_max_min_sp2}. 

For the second part, we show that the converse also holds. For each $k \in \sets K$, let
$$
{\bm \xi}^{k,\star} \; \in \; \argmax_{ {\bm \xi} \in \Xi({\bm w},\overline{\bm \xi}) } \;\; {\bm \xi}^\top {\bm C} \; {\bm x} + {\bm \xi}^\top {\bm D} \; {\bm w} + {\bm \xi}^\top {\bm Q} \; {\bm y}^k.
$$
Then, the optimal objective value of Problem~\eqref{eq:Kadapt_min_max_min_sp1} is expressible as
\begin{equation}
\min_{ k \in \sets K } \;  \left\{ ( {\bm \xi}^{k,\star} )^\top {\bm C} \; {\bm x} + ( {\bm \xi}^{k,\star} )^\top {\bm D} \; {\bm w} + ( {\bm \xi}^{k,\star} )^\top {\bm Q} \; {\bm y}^k \; : \; {\bm T} {\bm x} + {\bm V} {\bm w} + {\bm W}{\bm y}^k \leq {\bm H}{\bm \xi}  \; \; \; \; \forall {\bm \xi} \in \Xi({\bm w},\overline{\bm \xi}) \right\}.
\label{eq:Kadapt_obj_all_xi}
\end{equation}
Since ${\bm \xi}^{k,\star} \in \Xi({\bm w},\overline{\bm \xi})$, the solution $\{ {\bm \xi}^{k,\star} \}_{k \in \sets K}$ is feasible in Problem~\eqref{eq:Kadapt_min_max_min_sp2} with objective 
\begin{equation}
\min_{ k \in \sets K } \;  \left\{  ({\bm \xi}^{k,\star} )^\top {\bm C} \; {\bm x} + ({\bm \xi}^{k,\star})^\top {\bm D} \; {\bm w} + ({\bm \xi}^{k,\star})^\top {\bm Q} \; {\bm y}^k \; : \; {\bm T} {\bm x} + {\bm V} {\bm w} + {\bm W}{\bm y}^k \leq {\bm H}{\bm \xi}^{k,\star} \right\}.
\label{eq:Kadapt_obj_one_xi}
\end{equation}
If the optimal objective values of Problems~\eqref{eq:Kadapt_obj_all_xi} and~\eqref{eq:Kadapt_obj_one_xi} are equal, then we can directly conclude that the optimal objective value of Problem~\eqref{eq:Kadapt_min_max_min_sp2} exceeds that of Problem~\eqref{eq:Kadapt_min_max_min_sp1}. Suppose to the contrary that the optimal objective value of Problems~\eqref{eq:Kadapt_obj_one_xi} is strictly lower than that of Problem~\eqref{eq:Kadapt_obj_all_xi}. Then, there exists (at least one) $k \in \sets K$ that is feasible in~\eqref{eq:Kadapt_obj_one_xi} but infeasible in~\eqref{eq:Kadapt_obj_all_xi} and for each such $k$, there exists ${\bm \xi}^{k,'} \in \Xi({\bm w},\overline{\bm \xi})$ such that ${\bm T} {\bm x} + {\bm V} {\bm w} + {\bm W}{\bm y}^{\tilde k} \nleq {\bm H}{\bm \xi}^{k,'}$. We can construct a feasible solution $\{ \tilde {\bm \xi}^{k,\star} \}_{k \in \sets K}$ to Problem~\eqref{eq:Kadapt_min_max_min_sp2} with the same objective as Problem~\eqref{eq:Kadapt_obj_all_xi} as follows:
$$
\tilde{\bm \xi}^{k,\star} := 
\begin{cases}
{\bm \xi}^{k,\star} & \text{if } k \; : \; {\bm T} {\bm x} + {\bm V} {\bm w} + {\bm W}{\bm y}^k \leq {\bm H}{\bm \xi}  \; \; \; \; \forall {\bm \xi} \in \Xi({\bm w},\overline{\bm \xi}), \\
{\bm \xi}^{k,'} & \text{else}.
\end{cases}
$$
Indeed, the objective value attained by $\{\tilde{\bm \xi}^{k,\star}\}_{k \in \sets K}$ in Problem~\eqref{eq:Kadapt_min_max_min_sp2} is
$$
\begin{array}{cl}
& \displaystyle \min_{ k \in \sets K } \;  \left\{  (\tilde{\bm \xi}^{k,\star} )^\top {\bm C} \; {\bm x} + (\tilde{\bm \xi}^{k,\star})^\top {\bm D} \; {\bm w} + (\tilde{\bm \xi}^{k,\star})^\top {\bm Q} \; {\bm y}^k \; : \; {\bm T} {\bm x} + {\bm V} {\bm w} + {\bm W}{\bm y}^k \leq {\bm H}\tilde{\bm \xi}^{k,\star} \right\} \\
= & \displaystyle \min_{ k \in \sets K } \;  \left\{  (\tilde{\bm \xi}^{k,\star} )^\top {\bm C} \; {\bm x} + (\tilde{\bm \xi}^{k,\star})^\top {\bm D} \; {\bm w} + (\tilde{\bm \xi}^{k,\star})^\top {\bm Q} \; {\bm y}^k \; : \; {\bm T} {\bm x} + {\bm V} {\bm w} + {\bm W}{\bm y}^k \leq {\bm H}{\bm \xi}  \; \; \; \; \forall {\bm \xi} \in \Xi({\bm w},\overline{\bm \xi}) \right\}  \\
= &  \displaystyle \min_{ k \in \sets K } \;  \left\{  ({\bm \xi}^{k,\star} )^\top {\bm C} \; {\bm x} + ({\bm \xi}^{k,\star})^\top {\bm D} \; {\bm w} + ({\bm \xi}^{k,\star})^\top {\bm Q} \; {\bm y}^k \; : \; {\bm T} {\bm x} + {\bm V} {\bm w} + {\bm W}{\bm y}^k \leq {\bm H}{\bm \xi}  \; \; \; \; \forall {\bm \xi} \in \Xi({\bm w},\overline{\bm \xi}) \right\}, 
\end{array}
$$
where the first equality follows by construction since 
$$
\{ k \in \sets K \; : \; {\bm T} {\bm x} + {\bm V} {\bm w} + {\bm W}{\bm y}^k \leq {\bm H}{\bm \xi}  \; \; \; \; \forall {\bm \xi} \in \Xi({\bm w},\overline{\bm \xi})  \} 
\; = \;
\{ k \in \sets K \; : \; {\bm T} {\bm x} + {\bm V} {\bm w} + {\bm W}{\bm y}^k \leq {\bm H}\tilde{\bm \xi}^{k,\star}  \}
$$
and the second equality follows since 
$$
\tilde{\bm \xi}^{k,\star} = {\bm \xi}^{k,\star} \quad \forall k \in \sets K \; : \; {\bm T} {\bm x} + {\bm V} {\bm w} + {\bm W}{\bm y}^k \leq {\bm H}{\bm \xi}  \; \; \; \; \forall {\bm \xi} \in \Xi({\bm w},\overline{\bm \xi}).
$$
We have thus shown that the optimal objective value of Problem~\eqref{eq:Kadapt_min_max_min_sp2} is at least as large as that of Problem~\eqref{eq:Kadapt_min_max_min_sp1}.

Combining the first and second parts of the proof, we conclude that Problems~\eqref{eq:Kadapt_min_max_min_sp1} and~\eqref{eq:Kadapt_min_max_min_sp2} have the same optimal objective values, which concludes the proof. \Halmos
\endproof


\subsection{Proofs of Statements in Section~\ref{sec:kadaptability_objective}}
\label{sec:EC_kadaptability_objective}

\proof{Proof of Obervation~\ref{obs:Kadapt_obj_move_constraints_2}} Since Problem~\eqref{eq:Kadapt_obj} is equivalent to Problem~\eqref{eq:Kadapt_min_max_min_2} (by Lemma~\ref{lem:Kadapt_cstr_min_max_min}), it suffices to show that Problems~\eqref{eq:Kadapt_min_max_min_2} and~\eqref{eq:Kadapt_obj_move_constraints_2} are equivalent. 

First, note that for any choice of ${\bm w} \in \sets W$, the set $\Xi^K({\bm w})$ is non-empty. If there is no ${\bm x} \in \sets X$, ${\bm w} \in \sets W$, and ${\bm y} \in \sets Y$ such that ${\bm T} {\bm x} + {\bm V} {\bm w} + {\bm W}{\bm y} \leq {\bm h}$, then Problem~\eqref{eq:Kadapt_obj_move_constraints_2} is infeasible and has an optimal objective value of $+\infty$. Accordingly, Problem~\eqref{eq:Kadapt_min_max_min_2} also has an objective value of $+\infty$ since either its outer or inner minimization problems are infeasible. 

Suppose now that there exists ${\bm x} \in \sets X$, ${\bm w} \in \sets W$, and ${\bm y} \in \sets Y$ such that ${\bm T} {\bm x} + {\bm V} {\bm w} + {\bm W}{\bm y} \leq {\bm h}$. Then, Problems~\eqref{eq:Kadapt_min_max_min_2} and \eqref{eq:Kadapt_obj_move_constraints_2} are both feasible. Let $({\bm x},{\bm w},\{{\bm y}\}_{k\in\sets K})$ be a feasible solution for~\eqref{eq:Kadapt_obj_move_constraints_2}. Then, it is feasible in~\eqref{eq:Kadapt_min_max_min_2} and attains the same objective value in both problems since all second stage policies ${\bm y}^k$, $k \in \sets K$, satisfy the second-stage constraints in Problem~\eqref{eq:Kadapt_min_max_min_2}. Conversely, let $({\bm x},{\bm w},\{{\bm y}\}_{k\in\sets K})$ be feasible in Problem~\eqref{eq:Kadapt_min_max_min_2}. Since $\Xi^K({\bm w})$ is non-empty, there must exist $k^\star \in \sets K$ such that ${\bm T} {\bm x} + {\bm V} {\bm w} + {\bm W}{\bm y}^{k^\star} \leq {\bm h}$ (else the problem would have an optimal objective value of $+\infty$ and thus be infeasible, a contradiction). If ${\bm T} {\bm x} + {\bm V} {\bm w} + {\bm W}{\bm y}^{k} \leq {\bm h}$ for all $k \in \sets K$, then $({\bm x},{\bm w},\{{\bm y}\}_{k\in\sets K})$ is feasible in~\eqref{eq:Kadapt_obj_move_constraints_2} and attains the same objective value in both problems. On the other hand, if ${\bm T} {\bm x} + {\bm V} {\bm w} + {\bm W}{\bm y}^{k} > {\bm h}$ for some $k \in \sets K$, define
$$
\overline{\bm y}^k =
\begin{cases}
{\bm y}^k & \text{if } {\bm T} {\bm x} + {\bm V} {\bm w} + {\bm W}{\bm y}^{k} \leq {\bm h} \\
{\bm y}^{k^\star} & \text{else.}
\end{cases}
$$
Then, $({\bm x},{\bm w},\{\overline{\bm y}\}_{k\in\sets K})$ is feasible in~\eqref{eq:Kadapt_obj_move_constraints_2} and attains the same objective value in both problems. \Halmos 
\endproof

\proof{Proof of Observation~\ref{obs:evaluation_Kadapt_obj}} Fix $K \in \mathbb N$ and $({\bm x}, {\bm w}, \{ {\bm y}^k \}_{k \in \sets K})$ such that ${\bm x} \in \sets X$, ${\bm w} \in \sets W$, ${\bm y}^k \in \sets Y$. Assume, w.l.o.g.\ (see the Proof of Observation~\ref{obs:Kadapt_obj_move_constraints_2}) that ${\bm T}{\bm x} + {\bm V}{\bm w} + {\bm W}{\bm y}^k \leq {\bm h}$ for all $k\in \sets K$. From Observation~\ref{obs:Kadapt_obj_move_constraints_2}, the objective value of~\eqref{eq:Kadapt_obj} under this decision is equal to 
$$
    \begin{array}{cl}
         \maximize & \quad \displaystyle \min_{ k \in \sets K } \;\;  \left\{  ({\bm \xi}^k)^\top {\bm C} \; {\bm x} + ({\bm \xi}^k)^\top {\bm D} \; {\bm w} + ({\bm \xi}^k)^\top {\bm Q} \; {\bm y}^k  \right\}  \\
         \subjectto & \quad \overline{\bm \xi} \in \Xi, \; {\bm \xi}^k \in \Xi, \; k \in \sets K \\
         & \quad  {\bm w} \circ {\bm \xi}^k = {\bm w} \circ \overline{\bm \xi}  \quad \forall k \in \sets K.
    \end{array}
$$
We can write the problem above in epigraph form as an LP:
$$
    \begin{array}{cl}
         \maximize & \quad \tau   \\
         \subjectto & \quad \tau \in \reals, \; \overline{\bm \xi} \in \Xi, \; {\bm \xi}^k \in \Xi, k \in \sets K \\
         & \quad  \tau \leq   ({\bm \xi}^k)^\top {\bm C} \; {\bm x} + ({\bm \xi}^k)^\top {\bm D} \; {\bm w} + ({\bm \xi}^k)^\top {\bm Q} \; {\bm y}^k \quad \forall k \in \sets K \\
         & \quad  {\bm w} \circ {\bm \xi}^k = {\bm w} \circ \overline{\bm \xi}  \quad \forall k \in \sets K.
    \end{array}
$$
For any fixed $K$, the size of this LP is polynomial in the size of the input.\Halmos
\endproof

\proof{Proof of Theorem~\ref{thm:Kadapt_obj}} For any fixed $({\bm x}, {\bm w}, \{ {\bm y}^k \}_{k \in \sets K})$, we can express the inner maximization problem in~\eqref{eq:Kadapt_obj_move_constraints_2} in epigraph form as
$$
    \begin{array}{cl}
         \maximize & \quad \tau   \\
         \subjectto & \quad \tau \in \reals, \; \overline{\bm \xi} \in \reals^{N_\xi}, \; {\bm \xi}^k \in \reals^{N_\xi}, k \in \sets K \\
         & \quad  \tau \; \leq \;  ( {\bm C} \; {\bm x} +  {\bm D} \; {\bm w} + {\bm Q} \; {\bm y}^k )^\top {\bm \xi}^k \quad \forall k \in \sets K \\
         & \quad {\bm A}\overline{\bm \xi} \leq {\bm b} \\
         & \quad {\bm A}{\bm \xi}^k \leq {\bm b}  \quad \forall k \in \sets K \\
         & \quad  {\bm w} \circ {\bm \xi}^k = {\bm w} \circ \overline{\bm \xi}  \quad \forall k \in \sets K.
    \end{array}
$$
Strong LP duality (which applies since the feasible set is non-empty and since the problem is bounded by virtue of the boundedness of $\Xi$) implies that the optimal objective value of this problem coincides with the optimal objective value of its dual
$$
\begin{array}{cl}
\displaystyle \minimize & \displaystyle \quad {\bm b}^\top {\bm \beta}  + \sum_{k\in \mathcal K} {\bm b}^\top {\bm \beta}^k  \\
\text{subject to} & \quad {\bm \alpha} \in \mathbb R^K_+, \; {\bm \beta} \in \mathbb R^R_+, \; {\bm \beta}^k \in \mathbb R^R_+, \; {\bm \gamma}^k \in \mathbb R^{N_\xi}, \; k \in \mathcal K \\
& \quad \1^\top {\bm \alpha}  = 1  \\
& \quad     {\bm A}^\top{\bm \beta}^k +  {\bm w} \circ {\bm \gamma}^k = {\bm \alpha}_k \left( {\bm C} {\bm x} + {\bm D} {\bm w} + {\bm Q} {\bm y}^k \right) \quad \forall k \in \mathcal K \\
& \quad  {\bm A}^\top{\bm \beta}  =  \displaystyle \sum_{k \in \mathcal K}  {\bm w} \circ {\bm \gamma}^k.
\end{array}
$$
We can now group the outer minimization with the minimization above to obtain
$$
    \begin{array}{cl}
         \minimize & \quad   {\bm b}^\top {\bm \beta}  + \sum_{k\in \mathcal K} {\bm b}^\top {\bm \beta}^k\\
         \subjectto & \quad {\bm x} \in \sets X, \; {\bm w} \in \sets W, \; {\bm y}^k \in \sets Y, \; k \in \sets K \\
         & \quad {\bm \alpha} \in \mathbb R^K_+, \; {\bm \beta} \in \mathbb R^R_+, \; {\bm \beta}^k \in \mathbb R^R_+, \; {\bm \gamma}^k \in \mathbb R^{N_\xi}, \; k \in \mathcal K \\
        & \quad \1^\top {\bm \alpha}  = 1  \\
        & \quad     {\bm A}^\top{\bm \beta}^k +  {\bm w} \circ {\bm \gamma}^k = {\bm \alpha}_k \left( {\bm C} {\bm x} + {\bm D} {\bm w} + {\bm Q} {\bm y}^k \right) \quad \forall k \in \mathcal K \\
        & \quad  {\bm A}^\top{\bm \beta}  =  \displaystyle \sum_{k \in \mathcal K} {\bm w} \circ  {\bm \gamma}^k \\
         & \quad {\bm T} {\bm x} + {\bm V} {\bm w} + {\bm W}{\bm y}^k \leq {\bm h} \quad \forall k \in \sets K.
    \end{array}
$$
This concludes the proof. \Halmos

\proof{Proof of Corollary~\ref{cor:Kadapt_obj_MBLP}} 
The result follows directly from Theorem~\ref{thm:Kadapt_obj} by replacing the bilinear terms ${\bm w} \circ {\bm \gamma}^k$, ${\bm \alpha}_k {\bm x}$, ${\bm \alpha}_k {\bm w}$, and ${\bm \alpha}_k {\bm y}^k$ with auxiliary variables $\overline{\bm \gamma}^k \in \reals^{N_\xi}$, $\overline{\bm x}^k \in \reals^{N_x}_+$, $\overline {\bm w}^k \in \reals^{N_\xi}_+$, and $\overline {\bm y}^k  \in \reals^{N_y}_+$ such that
$$
\begin{array}{ccl}
\overline{\bm \gamma}^k = {\bm w} \circ {\bm \gamma}^k 
\quad & \Leftrightarrow & \quad
\overline{\bm \gamma}^k  \leq {\bm \gamma}^k + M (\1 - {\bm w}), \; \overline{\bm \gamma}^k  \leq M {\bm w}, \; \overline{\bm \gamma}^k \geq -M{\bm w}, \;  \overline{\bm \gamma}^k \geq {\bm \gamma}^k - M (\1 - {\bm w}), \\
\overline{\bm x}^k = {\bm \alpha}_k {\bm x} 
\quad & \Leftrightarrow & \quad 
\overline{\bm x}^k \leq {\bm x}, \; \overline{\bm x}^k \leq {\bm \alpha}_k \1 ,\; \overline{\bm x}^k \geq ({\bm \alpha}_k - 1) \1 + {\bm x}, \\
\overline{\bm w}^k = {\bm \alpha}_k {\bm w} 
\quad & \Leftrightarrow & \quad 
\overline{\bm w}^k \leq {\bm w}, \; \overline{\bm w}^k \leq {\bm \alpha}_k \1 ,\; \overline{\bm w}^k \geq ({\bm \alpha}_k -1) \1 + {\bm w},\\
\overline{\bm y}^k = {\bm \alpha}_k {\bm y}^k 
\quad & \Leftrightarrow & \quad 
\overline{\bm y}^k \leq {\bm y}^k, \; \overline{\bm y}^k \leq {\bm \alpha}_k \1 ,\; \overline{\bm y}^k \geq ({\bm \alpha}_k -1) \1 + {\bm y}^k,
\end{array}
$$
where in the last three cases we have exploited the fact that ${\bm x}$, ${\bm w}$, and ${\bm y}^k$ are binary and that ${\bm \alpha}^k \in [ {\bm 0} , \1 ]$.~\Halmos
\endproof

\proof{Proof of Observation~\ref{obs:grani_equivalence_DVzero_objective}}
It follows from the Proof of Corollary~\ref{cor:Kadapt_obj_MBLP} that Problem~\eqref{eq:Kadapt_obj} is equivalent to the following MBLP.
\begin{equation}
    \begin{array}{cl}
         \minimize & \quad \displaystyle  {\bm b}^\top {\bm \beta}  + \sum_{k\in \mathcal K} {\bm b}^\top {\bm \beta}^k\\
         \subjectto & \quad {\bm x} \in \sets X, \; {\bm w} \in \sets W, \; {\bm y}^k \in \sets Y, \; k \in \sets K \\
         & \quad {\bm \alpha} \in \mathbb R^K_+, \; {\bm \beta} \in \mathbb R^R_+, \; {\bm \beta}^k \in \mathbb R^R_+, \; {\bm \gamma}^k \in \mathbb R^{N_\xi}, \; k \in \mathcal K \\
        & \quad \overline{\bm \gamma}^k \in \reals^{N_\xi}, \; \overline{\bm x}^k \in \reals^{N_x}_+, \; \overline {\bm w}^k \in \reals^{N_\xi}_+, \; \overline {\bm y}^k  \in \reals^{N_y}_+, \; k \in \mathcal K  \\
        & \quad \1^\top {\bm \alpha}  = 1, \; {\bm A}^\top{\bm \beta}  =  \displaystyle \sum_{k \in \mathcal K} \overline{\bm \gamma}^k  \\
        & \quad     {\bm A}^\top{\bm \beta}^k +  \overline{\bm \gamma}^k =  {\bm C} \overline{\bm x}^k + {\bm D} \overline{\bm w}^k + {\bm Q} \overline{\bm y}^k \quad \forall k \in \mathcal K \\
        %
         & \quad {\bm T} {\bm x} + {\bm V} {\bm w} + {\bm W}{\bm y}^k \leq {\bm h} \quad \forall k \in \sets K \\
        & \left. \!\! \quad  \begin{array}{l}
         \overline{\bm x}^k \leq {\bm x}, \; \overline{\bm x}^k \leq {\bm \alpha}_k \1 ,\; \overline{\bm x}^k \geq ({\bm \alpha}_k - 1) \1 + {\bm x} \\
          \overline{\bm w}^k \leq {\bm w}, \; \overline{\bm w}^k \leq {\bm \alpha}_k \1 ,\; \overline{\bm w}^k \geq ({\bm \alpha}_k - 1) \1 + {\bm w}\\
          \overline{\bm y}^k \leq {\bm y}^k, \; \overline{\bm y}^k \leq {\bm \alpha}_k \1 ,\; \overline{\bm y}^k \geq ({\bm \alpha}_k  - 1) \1 + {\bm y}^k \\
          \overline{\bm \gamma}^k  \leq {\bm \gamma}^k + M (\1 - {\bm w}), \; \overline{\bm \gamma}^k  \leq M {\bm w} , \;
          \overline{\bm \gamma}^k \geq -M{\bm w}, \;  \overline{\bm \gamma}^k \geq {\bm \gamma}^k - M (\1 - {\bm w})
        \end{array} \quad \right\} \quad \forall k \in \sets K.
    \end{array}
\label{eq:Kadapt_obj_MBLP}
\end{equation}
Suppose that we are only in the presence of exogenous uncertainty, i.e., ${\bm w}=\1$, ${\bm D}={\bm 0}$, and ${\bm V}={\bm 0}$. Then, the last set of constraints in Problem~\eqref{eq:Kadapt_obj_MBLP} implies that $\overline{\bm \gamma}^k= {\bm \gamma}^k$ for all $k \in \sets K$. Since ${\bm \gamma}^k$ is free, the second and third constraints are equivalent to
$$
 {\bm A}^\top{\bm \beta}  =  \displaystyle \sum_{k \in \mathcal K}  {\bm C} \overline{\bm x}^k + {\bm Q} \overline{\bm y}^k  -   {\bm A}^\top{\bm \beta}^k.
$$
Exploiting the fact that ${\bm \alpha} \in \reals^K_+$, $\1^\top {\bm \alpha}=1$, and $\overline{\bm x}^k = {\bm \alpha}_k {\bm x}$, we can equivalently express this constraint as
$$
 {\bm A}^\top \left( {\bm \beta} + \sum_{k \in \mathcal K} {\bm \beta}^k \right)  = {\bm C} {\bm x} +  \displaystyle \sum_{k \in \mathcal K}  {\bm Q} \overline{\bm y}^k.
$$
We conclude that, in the presence of only exogenous uncertainty, Problem~\eqref{eq:Kadapt_obj_MBLP} is equivalent to
$$
    \begin{array}{cl}
         \minimize & \quad \displaystyle  {\bm b}^\top {\bm \beta}  \\
         \subjectto & \quad {\bm x} \in \sets X, \; {\bm w} \in \sets W, \; {\bm y}^k \in \sets Y, \; k \in \sets K \\
         & \quad {\bm \alpha} \in \mathbb R^K_+, \; {\bm \beta} \in \mathbb R^R_+, \; {\bm \beta}^k \in \mathbb R^R_+, \; \overline {\bm y}^k  \in \reals^{N_y}_+ , \;  k \in \mathcal K \\
        & \quad \1^\top {\bm \alpha}  = 1 , \; \displaystyle {\bm A}^\top  {\bm \beta} = {\bm C} {\bm x} +  \sum_{k \in \mathcal K}  {\bm Q} \overline{\bm y}^k \\
         & \quad {\bm T} {\bm x} + {\bm W}{\bm y}^k \leq {\bm h} \quad \forall k \in \sets K \\
         & \quad \overline{\bm y}^k \leq {\bm y}^k, \; \overline{\bm y}^k \leq {\bm \alpha}_k \1 ,\; \overline{\bm y}^k \geq ({\bm \alpha}_k  - 1) \1 + {\bm y}^k  \quad \forall k \in \sets K,
    \end{array}
$$
where we used the change of variables ${\bm \beta} \leftarrow {\bm \beta}  + \sum_{k\in \mathcal K}  {\bm \beta}^k$. We then recover the MBLP formulation of the $K$-adaptability counterpart of problems with exogenous uncertainty from~\cite{Hanasusanto2015}. Thus, our reformulation encompasses as a special case the one of~\cite{Hanasusanto2015}.~\Halmos
\endproof


\subsection{Proofs of Statements in Section~\ref{sec:kadaptability_constraint}}
\label{sec:EC_kadaptability_constraint}







\proof{Proof of Theorem~\ref{thm:hardness_endo_2}}
The proof is a direct consequence of Theorem~3 in~\cite{Hanasusanto2015}. Indeed, the authors show that evaluating the objective function of Problem~\eqref{eq:exo_2} is strongly NP-hard. Since Problem~\eqref{eq:exo_2} can be reduced in polynomial time to an instance of Problem~\eqref{eq:endo_3} by letting ${\bm D}={\bm 0}$, ${\bm V}={\bm 0}$, and ${\bm w}=\1$, this concludes the proof. 
\Halmos
\endproof


The proof below is a generalization of the proof of Proposition 1 in \cite{Hanasusanto2015} that operates in the \emph{lifted} uncertainty and decision spaces. Despite this key difference, the proof idea carries through.

\proof{Proof of Proposition~\ref{prop:Kadapt_ell}} Fix ${\bm x}$, ${\bm w}$, and $\{ {\bm y}^k \}_{k \in \sets K}$. We show that $\{ \Xi^K({\bm w},{\bm \ell}) \}_{ {\bm \ell} \in \sets L }$ is a cover of $\Xi^K({\bm w})$, i.e., that $\Xi^K({\bm w}) = \bigcup_{ {\bm \ell} \in \sets L }  \Xi^K({\bm w},{\bm \ell})$. Let $\{ {\bm \xi}^k \}_{k \in \sets K} \in \Xi^K({\bm w})$ and define
$$
{\bm \ell}_k = 
\begin{cases}
0 , \text{ if } {\bm T} {\bm x} + {\bm V} {\bm w} + {\bm W}{\bm y}^k \leq {\bm H}{\bm \xi}^k \\
\min\left\{ \ell \in \{1,\ldots,L\} \; : \; \left[{\bm T} {\bm x} + {\bm V} {\bm w} + {\bm W}{\bm y}^k \right]_{\ell} > [{\bm H}{\bm \xi}^k]_{\ell} \right\}  \text{, else.}
\end{cases} \forall k \in \sets K.
$$
Then, $\{ {\bm \xi}^k \}_{k \in \sets K} \in \Xi^K({\bm w},{\bm \ell})$. Moreover, by definition, we have $\Xi^K({\bm w},{\bm \ell}) \subseteq \Xi^K({\bm w})$ for all ${\bm \ell} \in \sets L$. Therefore $\{ \Xi^K({\bm w},{\bm \ell}) \}_{ {\bm \ell} \in \sets L }$ is a cover of $\Xi^K({\bm w})$. It then follows that
$$
\begin{array}{ccl}
     & \quad  \displaystyle \max_{ \{ {\bm \xi}^k \}_{k \in \sets K} \in \Xi^K({\bm w}) } & \;\;  \displaystyle \min_{ k \in \sets K } \;  \left\{  ({\bm \xi}^k)^\top {\bm C} \; {\bm x} + ({\bm \xi}^k)^\top {\bm D} \; {\bm w} + ({\bm \xi}^k)^\top {\bm Q} \; {\bm y}^k \; : \; {\bm T} {\bm x} + {\bm V} {\bm w} + {\bm W}{\bm y}^k \leq {\bm H}{\bm \xi}^k \right\}  \\
     = & \quad   \displaystyle \max_{ \{ {\bm \xi}^k \}_{k \in \sets K} \in \bigcup_{ {\bm \ell} \in \sets L  }\Xi^K({\bm w},{\bm \ell}) } & \;\;  \displaystyle \min_{ k \in \sets K } \;  \left\{  ({\bm \xi}^k)^\top {\bm C} \; {\bm x} + ({\bm \xi}^k)^\top {\bm D} \; {\bm w} + ({\bm \xi}^k)^\top {\bm Q} \; {\bm y}^k \; : \; {\bm T} {\bm x} + {\bm V} {\bm w} + {\bm W}{\bm y}^k \leq {\bm H}{\bm \xi}^k \right\}  \\
     = & \quad  \displaystyle \max_{{\bm \ell} \in \sets L} \;\; \max_{ \{ {\bm \xi}^k \}_{k \in \sets K} \in \Xi^K({\bm w},{\bm \ell}) } & \;\;  \displaystyle \min_{ k \in \sets K } \;  \left\{  ({\bm \xi}^k)^\top {\bm C} \; {\bm x} + ({\bm \xi}^k)^\top {\bm D} \; {\bm w} + ({\bm \xi}^k)^\top {\bm Q} \; {\bm y}^k \; : \; {\bm T} {\bm x} + {\bm V} {\bm w} + {\bm W}{\bm y}^k \leq {\bm H}{\bm \xi}^k \right\}.
\end{array}
$$
The definition of $\Xi^K( {\bm w} , {\bm \ell})$ implies that ${\bm \ell}_k =0$ if and only if $\; {\bm T} {\bm x} + {\bm V} {\bm w} + {\bm W}{\bm y}^k \leq {\bm H}{\bm \xi}^k$. This concludes the proof.\Halmos
\endproof


\proof{Proof of Theorem~\ref{thm:Kadapt_MBLP}}
The objective function of the approximate problem~\eqref{eq:Kadapt_ell_eps} is identical to
$$
 \displaystyle  \max_{{\bm \ell} \in \sets L} \;\; \max_{ \{ {\bm \xi}^k \}_{ k \in \sets K} \in \Xi^K_\epsilon({\bm w} , {\bm \ell}) } \;\;\min_{ {\bm \lambda} \in \Lambda_K({\bm \ell}) } \;\;  \left\{  \sum_{k \in \sets K} {\bm \lambda}_k \left[ ({\bm \xi}^k)^\top {\bm C} \; {\bm x} + ({\bm \xi}^k)^\top {\bm D} \; {\bm w} + ({\bm \xi}^k)^\top {\bm Q} \; {\bm y}^k \right] \right\},
$$
where $\Lambda_K({\bm \ell}) := \left\{ {\bm \lambda} \in \reals^K_+ \; : \; \1^\top {\bm \lambda} = 1, \; {\bm \lambda}_k = 0 \; \forall k \in \sets K : {\bm \ell}_k \neq 0 \right\}$. Note that $\Lambda_K({\bm \ell})=\emptyset$ if and only if ${\bm \ell} > {\bm 0}$. If $\Xi^K_\epsilon({\bm w},{\bm \ell}) = \emptyset$ for all ${\bm \ell} \in \sets L_+$, then the problem is equivalent to
$$
 \displaystyle  \max_{{\bm \ell} \in \partial \sets L} \;\; \max_{ \{ {\bm \xi}^k \}_{ k \in \sets K} \in \Xi^K_\epsilon({\bm w} , {\bm \ell}) } \;\;\min_{ {\bm \lambda} \in \Lambda_K({\bm \ell}) } \;\;  \left\{  \sum_{k \in \sets K} {\bm \lambda}_k \left[ ({\bm \xi}^k)^\top {\bm C} \; {\bm x} + ({\bm \xi}^k)^\top {\bm D} \; {\bm w} + ({\bm \xi}^k)^\top {\bm Q} \; {\bm y}^k \right] \right\},
$$
and we can apply the classical min-max theorem (since $\Lambda_K({\bm \ell})$ is nonempty for all ${\bm \ell} \in \partial \sets L$) to obtain the equivalent reformulation
$$
 \displaystyle  \max_{{\bm \ell} \in \partial \sets L} \;\;\min_{ {\bm \lambda} \in \Lambda_K({\bm \ell}) } \;\; \max_{ \{ {\bm \xi}^k \}_{ k \in \sets K} \in \Xi^K_\epsilon({\bm w} , {\bm \ell}) } \;\;  \left\{  \sum_{k \in \sets K} {\bm \lambda}_k \left[ ({\bm \xi}^k)^\top {\bm C} \; {\bm x} + ({\bm \xi}^k)^\top {\bm D} \; {\bm w} + ({\bm \xi}^k)^\top {\bm Q} \; {\bm y}^k \right] \right\},
$$
which in turn is equivalent to
$$
 \displaystyle \min_{ \begin{smallmatrix} {\bm \lambda}({\bm \ell}) \in \Lambda_K({\bm \ell}), \\ {\bm \ell} \in \partial \sets L  \end{smallmatrix} } \;\;  \max_{{\bm \ell} \in \partial \sets L} \;\;\max_{ \{ {\bm \xi}^k \}_{ k \in \sets K} \in \Xi^K_\epsilon({\bm w} , {\bm \ell}) } \;\;  \left\{  \sum_{k \in \sets K} {\bm \lambda}_k({\bm \ell}) \left[ ({\bm \xi}^k)^\top {\bm C} \; {\bm x} + ({\bm \xi}^k)^\top {\bm D} \; {\bm w} + ({\bm \xi}^k)^\top {\bm Q} \; {\bm y}^k \right] \right\}.
$$
If, on the other hand, $\Xi^K_\epsilon({\bm w},{\bm \ell}) \neq \emptyset$ for some ${\bm \ell} \in \sets L_+$, then the objective function in~\eqref{eq:Kadapt_ell_eps} evaluates to
$+\infty$. Using an epigraph reformulation, we thus conclude that~\eqref{eq:Kadapt_ell_eps} is equivalent to the problem
\begin{equation}
    \begin{array}{cll}
         \minimize & \;\;  \tau  \\
         \subjectto & \;\; {\bm x} \in \sets X, \; {\bm w} \in \sets W, \; {\bm y}^k \in \sets Y, \; k \in \sets K \\
         & \;\; \tau \in \mathbb R, \; {\bm \lambda}({\bm \ell}) \in \Lambda_K({\bm \ell}), \; {\bm \ell} \in \partial \sets L \\
         & \;\; \tau \geq \displaystyle \sum_{k \in \sets K} {\bm \lambda}_k({\bm \ell}) \left[ ({\bm \xi}^k)^\top {\bm C} \; {\bm x} + ({\bm \xi}^k)^\top {\bm D} \; {\bm w} + ({\bm \xi}^k)^\top {\bm Q} \; {\bm y}^k \right] & \quad \forall {\bm \ell} \in \partial \sets L, \; \{{\bm \xi}^k \}_{ k \in \sets K} \in \Xi^K_\epsilon({\bm w} , {\bm \ell}) \\
         & \;\; \Xi^K_\epsilon({\bm w},{\bm \ell}) = \emptyset & \quad \forall {\bm \ell}\in \sets L_+.
    \end{array}
\label{eq:Kadapt_ell_eps_equiv}
\end{equation}
The semi-infinite constraint associated with ${\bm \ell} \in \partial \sets L$ is satisfied if and only if the optimal value of
$$
\begin{array}{cll}
     \maximize & \quad \displaystyle \sum_{k \in \sets K} {\bm \lambda}_k({\bm \ell}) \left[ ({\bm \xi}^k)^\top {\bm C} \; {\bm x} + ({\bm \xi}^k)^\top {\bm D} \; {\bm w} + ({\bm \xi}^k)^\top {\bm Q} \; {\bm y}^k \right]  \\
     \subjectto & \quad \overline {\bm \xi} \in \reals^{N_\xi}, \; {\bm \xi}^k \in \reals^{N_\xi}, \; k \in \sets K \\
     & \quad {\bm A} \overline {\bm \xi} \leq {\bm b} \\
     & \quad {\bm A} {\bm \xi}^k \leq {\bm b} & \quad \forall k \in \sets K \\
     & \quad {\bm T} {\bm x} + {\bm V} {\bm w} + {\bm W}{\bm y}^k \leq {\bm H}{\bm \xi}^k & \quad \forall k \in \sets K : {\bm \ell}_k = 0 \\
     & \quad \left[{\bm T} {\bm x} + {\bm V} {\bm w} + {\bm W}{\bm y}^k \right]_{{\bm \ell}_k} \geq [{\bm H}{\bm \xi}^k]_{{\bm \ell}_k} + \epsilon & \quad \forall k \in \sets K : {\bm \ell}_k \neq 0 \\
     & \quad {\bm w} \circ {\bm \xi}^k = {\bm w} \circ \overline {\bm \xi}  & \quad \forall k \in \sets K
\end{array}
$$
does not exceed $\tau$. Strong linear programming duality implies that this problem attains the same optimal value as its dual problem which is given by
$$
\begin{array}{cll}
     \minimize & \multicolumn{2}{l}{ \quad \displaystyle {\bm b}^\top \left( {\bm \alpha} + \sum_{k \in \sets K} {\bm \alpha}^k \right) - \sum_{ \begin{smallmatrix} k \in \sets K  : \\ {\bm \ell}_k = 0 \end{smallmatrix}} ( {\bm T} {\bm x} + {\bm V} {\bm w} + {\bm W}{\bm y}^k )^\top {\bm \beta}^k + \sum_{ \begin{smallmatrix} k \in \sets K  : \\ {\bm \ell}_k \neq 0 \end{smallmatrix}} \left( \left[{\bm T} {\bm x} + {\bm V} {\bm w} + {\bm W}{\bm y}^k \right]_{{\bm \ell}_k} - \epsilon \right) {\bm \gamma}_k }   \\
     \subjectto & \quad  {\bm \alpha} \in \reals^R_+, \; {\bm \alpha}^k \in \reals^R_+, \; {\bm \beta}^k \in \reals^L_+, \; k \in \sets K, \; {\bm \gamma} \in \reals^K_+, \; {\bm \eta}^k \in \reals^{N_\xi}, \; k \in \sets K \\
     & \displaystyle \quad  {\bm A}^\top {\bm \alpha} = \sum_{k \in \sets K} {\bm w} \circ {\bm \eta}^k  \\
     & \quad {\bm A}^\top {\bm \alpha}^k - {\bm H}^\top {\bm \beta}^k  + {\bm w} \circ {\bm \eta}^k = {\bm \lambda}_k({\bm \ell}) \left[ {\bm C} \; {\bm x} + {\bm D} \; {\bm w} + {\bm Q} \; {\bm y}^k  \right]      & \quad \forall k \in \sets K : {\bm \ell}_k = 0 \\
     & \quad {\bm A}^\top {\bm \alpha}^k + [ {\bm H} ]_{ {\bm \ell}_k } {\bm \gamma}_k + {\bm w} \circ {\bm \eta}^k   = {\bm \lambda}_k({\bm \ell}) \left[ {\bm C} \; {\bm x} + {\bm D} \; {\bm w} + {\bm Q} \; {\bm y}^k  \right]  & \quad \forall k \in \sets K : {\bm \ell}_k \neq 0.
\end{array}
$$
Strong duality holds because the dual problem is always feasible. Indeed, one can show that the compactness of $\Xi$ implies that $\{ {\bm A}^\top {\bm \alpha} : {\bm \alpha} \geq {\bm 0}\} = \reals^{N_\xi}$. Note that the first constraint set in Problem \eqref{eq:Kadapt_MBLP} ensures that the optimal value of this dual problem does not exceed $\tau$ for all ${\bm \ell} \in \partial \sets L$.

The last constraint in \eqref{eq:Kadapt_ell_eps_equiv} is satisfied for ${\bm \ell} \in \sets L_+$ whenever the linear program
$$
\begin{array}{cll}
     \maximize & \quad 0  \\
     \subjectto & \quad \overline {\bm \xi} \in \reals^{N_\xi}, \; {\bm \xi}^k \in \reals^{N_\xi}, \; k \in \sets K \\
     & \quad {\bm A} \overline {\bm \xi} \leq {\bm b} \\
     & \quad {\bm A} {\bm \xi}^k \leq {\bm b} & \quad \forall k \in \sets K \\
     & \quad \left[{\bm T} {\bm x} + {\bm V} {\bm w} + {\bm W}{\bm y}^k \right]_{{\bm \ell}_k} \geq [{\bm H}{\bm \xi}^k]_{{\bm \ell}_k} + \epsilon & \quad \forall k \in \sets K \\
     & \quad {\bm w} \circ {\bm \xi}^k = {\bm w} \circ \overline {\bm \xi}  & \quad \forall k \in \sets K
\end{array}
$$
is infeasible. The dual to this problem reads
$$
\begin{array}{cll}
     \minimize & \multicolumn{2}{l}{ \quad \displaystyle {\bm b}^\top \left( {\bm \alpha} + \sum_{k \in \sets K} {\bm \alpha}^k \right) + \sum_{ k \in \sets K } \left( \left[{\bm T} {\bm x} + {\bm V} {\bm w} + {\bm W}{\bm y}^k \right]_{{\bm \ell}_k} - \epsilon \right) {\bm \gamma}_k }   \\
     \subjectto & \quad  {\bm \alpha} \in \reals^R_+, \; {\bm \alpha}^k \in \reals^R_+, \; k \in \sets K, \; {\bm \gamma} \in \reals^K_+, \; {\bm \eta}^k \in \reals^{N_\xi}, \; k \in \sets K \\
     & \quad \displaystyle {\bm A}^\top {\bm \alpha} = \sum_{k \in \sets K} {\bm w} \circ {\bm \eta}^k  \\
     & \quad {\bm A}^\top {\bm \alpha}^k + [ {\bm H} ]_{ {\bm \ell}_k } {\bm \gamma}_k + {\bm w} \circ {\bm \eta}^k   = {\bm 0} & \quad \forall k \in \sets K.
\end{array}
$$
The feasible set of this dual is a cone and thus feasible (set ${\bm \alpha}={\bm 0}$, ${\bm \eta}^k={\bm 0}$, ${\bm \gamma}_k=0$, $k\in \sets K$). Therefore, strong LP duality applies and the primal is infeasible if and only if the dual is unbounded. Since the feasible set of the dual is a cone, the dual is unbounded if and only if there exists a feasible solution attaining an objective value of $-1$. \Halmos
\endproof

\proof{Proof of Observation~\ref{obs:grani_equivalence_DVzero}}
Suppose that we are only in the presence of exogenous uncertainty, i.e., ${\bm w}=\1$, ${\bm D}={\bm 0}$, and ${\bm V}={\bm 0}$. Then, Problem~\eqref{eq:Kadapt_MBLP} reduces to
\begin{equation}
    \begin{array}{cll}
         \min & \;\;  \tau  \\
         \st & \;\; \tau \in \mathbb R, \; {\bm x} \in \sets X, \; {\bm y}^k \in \sets Y, \; k \in \sets K \\
         & \;\;  {\bm \alpha}({\bm \ell}) \in \reals^R_+, \; {\bm \alpha}^k({\bm \ell}) \in \reals^R_+, \; k \in \sets K, \; {\bm \gamma}({\bm \ell}) \in \reals^K_+, \; {\bm \eta}^k({\bm \ell}) \in \reals^{N_\xi}, \; k \in \sets K, {\bm \ell} \in \sets L \\
         &   \left. \begin{array}{l}
         \; {\bm \lambda}({\bm \ell}) \in \Lambda_K({\bm \ell}) ,\; {\bm \beta}^k({\bm \ell}) \in \reals^L_+, \; k \in \sets K,  \\
        \displaystyle {\bm A}^\top {\bm \alpha}({\bm \ell}) = \sum_{k \in \sets K} {\bm \eta}^k({\bm \ell})  \\
        {\bm A}^\top {\bm \alpha}^k({\bm \ell}) - {\bm H}^\top {\bm \beta}^k({\bm \ell})  + {\bm \eta}^k({\bm \ell}) = {\bm \lambda}_k({\bm \ell}) \left[ {\bm C} \; {\bm x} + {\bm Q} \; {\bm y}^k  \right]  \quad \forall k \in \sets K : {\bm \ell}_k = 0 \\
        {\bm A}^\top {\bm \alpha}^k({\bm \ell}) + [ {\bm H} ]_{ {\bm \ell}_k } {\bm \gamma}_k({\bm \ell}) +  {\bm \eta}^k({\bm \ell})   = {\bm \lambda}_k({\bm \ell}) \left[ {\bm C} \; {\bm x}  + {\bm Q} \; {\bm y}^k  \right]   \quad \forall k \in \sets K : {\bm \ell}_k \neq 0 \\
         \tau \geq \displaystyle {\bm b}^\top  \left( {\bm \alpha}({\bm \ell}) + \sum_{k \in \sets K}  {\bm \alpha}^k({\bm \ell}) \right) - \sum_{ \begin{smallmatrix} k \in \sets K  : \\ {\bm \ell}_k = 0 \end{smallmatrix}} ( {\bm T} {\bm x} + {\bm W}{\bm y}^k )^\top {\bm \beta}^k({\bm \ell}) \\
         \qquad \qquad \qquad \displaystyle + \sum_{ \begin{smallmatrix} k \in \sets K  : \\ {\bm \ell}_k \neq 0 \end{smallmatrix}} \left( \left[{\bm T} {\bm x}  + {\bm W}{\bm y}^k \right]_{{\bm \ell}_k} - \epsilon \right) {\bm \gamma}_k({\bm \ell})
         \end{array} \quad \right\} \quad \forall {\bm \ell} \in \partial \sets L  \\
         &   \left. \begin{array}{l}
          \displaystyle {\bm A}^\top {\bm \alpha}({\bm \ell}) =  \sum_{ k \in \sets K} {\bm \eta}^k({\bm \ell})   \\
          {\bm A}^\top {\bm \alpha}^k({\bm \ell}) + [ {\bm H} ]_{ {\bm \ell}_k({\bm \ell}) } {\bm \gamma}_k({\bm \ell}) +  {\bm \eta}^k ({\bm \ell})  = {\bm 0} \quad \forall k \in \sets K\\
          \displaystyle {\bm b}^\top  \left( {\bm \alpha}({\bm \ell}) + \sum_{k \in \sets K} {\bm \alpha}^k({\bm \ell}) \right) + \sum_{ k \in \sets K } \left( \left[{\bm T} {\bm x} + {\bm W}{\bm y}^k \right]_{{\bm \ell}_k} - \epsilon \right) {\bm \gamma}_k({\bm \ell}) \leq -1
         \end{array} \quad \right\}
         \quad \forall {\bm \ell}\in \sets L_+.
    \end{array}
\label{eq:Kadapt_MBLP_exo_ours}
\end{equation}
Since ${\bm \eta}^k({\bm \ell})$ is free for all $k \in \sets K$ and ${\bm \ell} \in \sets L$, the first set of constraints associated with ${\bm \ell}\in \partial \sets L$ in~\eqref{eq:Kadapt_MBLP_exo_ours} is equivalent to
$$
 \displaystyle {\bm A}^\top  \left( {\bm \alpha}({\bm \ell}) + \sum_{k\in \sets K} {\bm \alpha}^k({\bm \ell}) \right)  - \sum_{ \begin{smallmatrix} k \in \sets K : \\ {\bm \ell}_k = 0 \end{smallmatrix} } {\bm H}^\top {\bm \beta}^k({\bm \ell}) + \sum_{ \begin{smallmatrix} k \in \sets K :\\ {\bm \ell}_k \neq 0 \end{smallmatrix} } [ {\bm H} ]_{{\bm \ell}_k} {\bm \gamma}_k({\bm \ell}) = {\bm C}{\bm x} + \sum_{k \in \sets K} {\bm \lambda}_k({\bm \ell}) \cdot {\bm Q}{\bm y}^k,
$$
where we have exploited the fact that ${\bm \lambda}({\bm \ell})\geq {\bm 0}$ and $\1^\top {\bm \lambda}({\bm \ell})=1$.
Similarly, the first set of constraints associated with ${\bm \ell}\in \sets L_+$ in~\eqref{eq:Kadapt_MBLP_exo_ours} is equivalent to
$$
\displaystyle {\bm A}^\top  \left( {\bm \alpha}({\bm \ell}) + \sum_{k \in \sets K}  {\bm \alpha}^k({\bm \ell}) \right) + \sum_{ k \in \sets K} [ {\bm H} ]_{{\bm \ell}_k} {\bm \gamma}_k({\bm \ell}) = {\bm 0}.
$$
We conclude that, in the presence of only exogenous uncertainty, Problem~\eqref{eq:Kadapt_MBLP} reduces to
\begin{equation}
    \begin{array}{cll}
         \min & \;\;  \tau  \\
         \st & \;\; \tau \in \mathbb R, \; {\bm x} \in \sets X, \; {\bm y}^k \in \sets Y, \; k \in \sets K \\
         & \;\;  {\bm \alpha}({\bm \ell}) \in \reals^R_+, \; {\bm \gamma}({\bm \ell}) \in \reals^K_+,  \; {\bm \ell}\in \sets L \\
         &   \left. \begin{array}{l}
         \; {\bm \lambda}({\bm \ell}) \in \Lambda_K({\bm \ell}) ,\; {\bm \beta}^k({\bm \ell}) \in \reals^L_+, \; k \in \sets K,  \\
        \displaystyle {\bm A}^\top  {\bm \alpha}({\bm \ell})  - \sum_{ \begin{smallmatrix} k \in \sets K : \\ {\bm \ell}_k = 0 \end{smallmatrix} } {\bm H}^\top {\bm \beta}^k({\bm \ell}) + \sum_{ \begin{smallmatrix} k \in \sets K :\\ {\bm \ell}_k \neq 0 \end{smallmatrix} } [ {\bm H} ]_{{\bm \ell}_k} {\bm \gamma}_k({\bm \ell}) = {\bm C}{\bm x} + \sum_{k \in \sets K} {\bm \lambda}_k ({\bm \ell}) \cdot {\bm Q}{\bm y}^k \\
         \tau \geq \displaystyle {\bm b}^\top   {\bm \alpha}({\bm \ell})  - \sum_{ \begin{smallmatrix} k \in \sets K  : \\ {\bm \ell}_k = 0 \end{smallmatrix}} ( {\bm T} {\bm x}  + {\bm W}{\bm y}^k )^\top {\bm \beta}^k({\bm \ell}) + \sum_{ \begin{smallmatrix} k \in \sets K  : \\ {\bm \ell}_k \neq 0 \end{smallmatrix}} \left( \left[{\bm T} {\bm x}  + {\bm W}{\bm y}^k \right]_{{\bm \ell}_k} - \epsilon \right) {\bm \gamma}_k({\bm \ell})
         \end{array} \quad \right\} \quad \forall {\bm \ell} \in \partial \sets L  \\
         &   \left. \begin{array}{l}
          \displaystyle {\bm A}^\top  {\bm \alpha}({\bm \ell}) + \sum_{ k \in \sets K} [ {\bm H} ]_{{\bm \ell}_k} {\bm \gamma}_k({\bm \ell}) = {\bm 0} \\
          \displaystyle {\bm b}^\top  {\bm \alpha}({\bm \ell})  + \sum_{ k \in \sets K } \left( \left[{\bm T} {\bm x}  + {\bm W}{\bm y}^k \right]_{{\bm \ell}_k} - \epsilon \right) {\bm \gamma}_k({\bm \ell}) \leq -1
         \end{array} \quad \right\}
         \quad \forall {\bm \ell}\in \sets L_+,
    \end{array}
\label{eq:Kadapt_MBLP_exo}
\end{equation}
where we use the change of variables ${\bm \alpha}({\bm \ell}) \leftarrow \left( {\bm \alpha}({\bm \ell}) + \sum_{k \in \sets K} {\bm \alpha}^k({\bm \ell}) \right)$. We thus recover the MBLP formulation of the $K$-adaptability problem from~\cite{Hanasusanto2015}, which concludes the proof.~\Halmos
\endproof


\section{Proofs of Statements in Section~\ref{sec:pwl_convex_objective}}
\label{sec:EC_pwl_convex_objective}



\proof{Proof of Lemma~\ref{lem:Kadapt_min_max_regret_mo}} It suffices to show that, for any fixed ${\bm x} \in \sets X$, ${\bm w} \in \sets W$, ${\bm y}^k \in \sets Y$, $k \in \sets K$, and $\overline{\bm \xi}\in \Xi$,
\begin{equation}
    \begin{array}{cll}
        \displaystyle \min_{k\in\sets K}\;\;\max_{{\bm \xi} \in \Xi({\bm w} , \overline {\bm \xi})}\;\;  \left\{ \max_{i\in \sets I} \;\;{\bm \xi}^\top {\bm C^i} \; {\bm x} + {\bm \xi}^\top {\bm D^i} \; {\bm w} + {\bm \xi}^\top {\bm Q^i} \; {\bm y^k} \right\}
    \end{array}
\label{eq:PE_proof_1}
\end{equation}
and
\begin{equation}
    \begin{array}{cll}
        \displaystyle \max_{ \begin{smallmatrix} {\bm \xi}^k \in \Xi({\bm w} , \overline {\bm \xi}) , \\ k \in \sets K \end{smallmatrix}  }\;\; \min_{k\in\sets K}\;\;  \left\{ \max_{i\in \sets I} \;\; ({\bm \xi}^k)^\top {\bm C^i} \; {\bm x} + ({\bm \xi}^k)^\top {\bm D^i} \; {\bm w} + ({\bm \xi}^k)^\top {\bm Q^i} \; {\bm y^k} \right\}
    \end{array}
\label{eq:PE_proof_2}
\end{equation}
are equivalent.

First, note that Problem~\eqref{eq:PE_proof_1} is always feasible and has a finite objective by virtue of the compactness of $\Xi({\bm w} , \overline {\bm \xi})$ which is non-empty. Similarly, Problem~\eqref{eq:PE_proof_2} is always feasible and has a finite objective.

We now show that both problems have the same objective. Let $\tilde k$ and $\{\tilde{\bm \xi}^k\}_{k\in\sets K}$ be feasible in~\eqref{eq:PE_proof_1} and~\eqref{eq:PE_proof_2}, respectively. The objective value attained by $\tilde{k}$ in Problem~\eqref{eq:PE_proof_1} is
$$
\displaystyle \max_{{\bm \xi} \in \Xi({\bm w} , \overline {\bm \xi})}\;\;  \left\{ \max_{i\in \sets I} \;\;{\bm \xi}^\top {\bm C^i} \; {\bm x} + {\bm \xi}^\top {\bm D^i} \; {\bm w} + {\bm \xi}^\top {\bm Q^i} \; {\bm y^{\tilde{k}}} \right\}.
$$
Accordingly, the objective value attained by $\{\tilde{\bm \xi}^k\}_{k\in\sets K}$ in Problem~\eqref{eq:PE_proof_2} is
$$
\displaystyle \min_{k \in\sets K}\;\;  \left\{ \max_{i\in \sets I} \;\;{(\tilde{\bm \xi}^k)}^\top {\bm C^i} \; {\bm x} + {(\tilde{\bm \xi}^k)}^\top {\bm D^i} \; {\bm w} + {(\tilde{\bm \xi}^k)}^\top {\bm Q^i} \; {\bm y^{k}} \right\}.
$$
Note that
$$
\begin{array}{cl}
& \quad \displaystyle \min_{k \in\sets K}\;\;  \left\{ \max_{i\in \sets I} \;\;{(\tilde{\bm \xi}^k)}^\top {\bm C^i} \; {\bm x} + {(\tilde{\bm \xi}^k)}^\top {\bm D^i} \; {\bm w} + {(\tilde{\bm \xi}^k)}^\top {\bm Q^i} \; {\bm y^{k}} \right\}\\
\leq & \quad \displaystyle \max_{i\in \sets I} \;\;{(\tilde{\bm \xi}^{\tilde k})}^\top {\bm C^i} \; {\bm x} + {(\tilde{\bm \xi}^{\tilde k})}^\top {\bm D^i} \; {\bm w} + {(\tilde{\bm \xi}^{\tilde k})}^\top {\bm Q^i} \; {\bm y}^{\tilde k} \\
\leq & \quad \displaystyle \max_{ {\bm \xi} \in \Xi({\bm w},\overline{\bm \xi}) } \; \; \max_{i\in \sets I}\;\;{\bm \xi}^\top {\bm C} \; {\bm x} + {\bm \xi}^\top {\bm D} \; {\bm w} + {\bm \xi}^\top {\bm Q} \; {\bm y}^{\tilde k}.
\end{array}
$$
Since the choice of $\tilde k \in \sets K$ and $\tilde{\bm \xi}^k \in \Xi({\bm w},\overline{\bm \xi})$ was arbitrary, it follows that Problem~\eqref{eq:PE_proof_1} upper bounds Problem~\eqref{eq:PE_proof_2}.

Next, we show that the converse also holds. Let
$$
\bm \xi^{k,*}\in \displaystyle \argmax_{{\bm \xi} \in \Xi({\bm w},\overline{\bm \xi})}\;\;\left\{ \max_{i\in \sets I} \;\;{\bm \xi}^\top {\bm C^i} \; {\bm x} + {\bm \xi}^\top {\bm D^i} \; {\bm w} + {\bm \xi}^\top {\bm Q^i} \; {\bm y^k} \right\}.
$$
Then, the optimal objective value of Problem~\eqref{eq:PE_proof_1} is expressible as
$$
\displaystyle \min_{k \in\sets K}\;\;  \left\{ \max_{i\in \sets I} \;\;{({\bm \xi}^{k,*})}^\top {\bm C^i} \; {\bm x} + {({\bm \xi}^{k,*})}^\top {\bm D^i} \; {\bm w} + {({\bm \xi}^{k,*})}^\top {\bm Q^i} \; {\bm y^{k}} \right\}.
$$
The solution $\{{\bm \xi}^{k,*}\}_{k\in\sets K}$ is feasible in~\eqref{eq:PE_proof_2} with objective
$$
\displaystyle \min_{k \in\sets K}\;\;  \left\{ \max_{i\in \sets I} \;\;{({\bm \xi}^{k,*})}^\top {\bm C^i} \; {\bm x} + {({\bm \xi}^{k,*})}^\top {\bm D^i} \; {\bm w} + {({\bm \xi}^{k,*})}^\top {\bm Q^i} \; {\bm y^{k}} \right\}.
$$
Thus, the optimal objective value of Problem~\eqref{eq:PE_proof_2} upper bounds that of Problem~\eqref{eq:PE_proof_1}.

Combining the two parts of the proof, we conclude that Problems~\eqref{eq:PE_proof_1} and~\eqref{eq:PE_proof_2} are equivalent.\Halmos
\endproof


\proof{Proof of Theorem~\ref{thm:Kadapt_mo_MBLP}} The objective function of Problem~\eqref{eq:Kadapt_min_max_regret_mo} is expressible as
$$
\max_{  \overline {\bm \xi} \in \Xi } \;\; \max_{ {\bm \xi^k} \in \Xi({\bm w} , \overline {\bm \xi}) , k \in \sets K}  \;\;\min_{ {k} \in {\sets K} } \;\; \left\{ \max_{i\in \sets I} \;\;{\bm \xi}^\top {\bm C^i} \; {\bm x} + {\bm \xi}^\top {\bm D^i} \; {\bm w} + {\bm \xi}^\top {\bm Q^i} \; {\bm y^k} \right\}.
$$
Using an epigraph reformulation, we can write it equivalently as
\begin{equation}
\begin{array}{cl}
    \maximize & \quad \tau \\
    \subjectto & \quad \tau\in \reals,\; \overline{\bm \xi}\in\Xi,\; {\bm \xi^k} \in \Xi({\bm w}, \overline {\bm \xi}), \; k \in \sets K \\
    & \quad \tau \; \leq \; \displaystyle{\max_{i\in\sets I}}\;\; (\bm\xi^k)^\top \bm C^{i}\;{\bm x} + {(\bm \xi^k)}^\top {\bm D^i} \; {\bm w} + {(\bm \xi^k)}^\top {\bm Q^i} \; {\bm y^k} \quad\forall k \in \sets K.
\end{array}
\label{eq:obj_Kadapt_mo}
\end{equation}
Noting that, for each $k \in \sets K$, the choice of $i \in \sets K$ can be made, in conjuction with the choice in $\tau$, $\overline{\bm \xi}$, and ${\bm \xi}^k$, $k\in \sets K$, Problem~\eqref{eq:obj_Kadapt_mo} can be written equivalently as
\begin{equation}
\begin{array}{ccl}
    \displaystyle \mathop{\maximize_{ {\bm i}_k\in\sets I , \;  k \in \sets K  }} \quad \quad & \max & \quad \tau \\
    & \st & \quad \tau\in\mathbb{R},\; \overline{\bm \xi}\in\Xi,\; {\bm \xi^k} \in \Xi({\bm w} , \overline {\bm \xi}),\; \forall k\in\sets K \\
    && \quad \tau \; \leq \; (\bm\xi^k)^\top {\bm C}^{{\bm i}_k}\;{\bm x} + ({\bm \xi}^k)^\top {\bm D}^{{\bm i}_k} \; {\bm w} + ({\bm \xi}^k)^\top {\bm Q}^{{\bm i}_k} \; {\bm y}^k \quad\forall k\in\sets K.
\end{array}
\label{eq:obj_Kadapt_mo_2}
\end{equation}

Dualizing the inner maximization problem yields
\begin{equation}
\begin{array}{ccl}
     \displaystyle \mathop{\maximize}_{ {\bm i}_k\in\sets I , \;  k \in \sets K  } \quad \quad & \min &\quad \bm b^\top\bm \beta + \sum_{k\in\sets K}\bm b^\top\bm\beta^k \\
     & \st & \quad {\bm \alpha} \in \mathbb R^K_+, \; {\bm \beta} \in \mathbb R^R_+, \; {\bm \beta}^k \in \mathbb R^R_+, \; {\bm \gamma}^k \in \mathbb R^{N_\xi}, \;\forall k \in \mathcal K \\
&& \quad \1^\top {\bm \alpha}  = 1  \\
&& \quad     {\bm A}^\top{\bm \beta}^k +  {\bm w} \circ {\bm \gamma}^k = {\bm \alpha}_k \left( {\bm C^{{\bm i}_k}} {\bm x} + {\bm D^{{\bm i}_k}} {\bm w} + {\bm Q^{{\bm i}_k}} {\bm y}^k \right) \quad \forall k \in \mathcal K \\
&& \quad  {\bm A}^\top{\bm \beta}  =  \displaystyle \sum_{k \in \mathcal K}  {\bm w} \circ {\bm \gamma}^k.
\end{array}
\label{eq:obj_Kadapt_mo_3}
\end{equation}
Equivalence of Problems~\eqref{eq:obj_Kadapt_mo_2} and~\eqref{eq:obj_Kadapt_mo_3} follows by strong LP duality which applies since the inner maximization problem in~\eqref{eq:obj_Kadapt_mo_2} is feasible and bounded. We next interchange the max and min operators, indexing each of the decision variables by ${\bm i}: = ({\bm i}_1,\ldots,{\bm i}_k) \in \sets I^K$. We obtain 
$$
\begin{array}{cl}
     \displaystyle \minimize  & \quad \displaystyle \max_{ {\bm i} \in\sets I^K }  \quad \bm b^\top\bm \beta^{\bm i} + \sum_{k\in\sets K}\bm b^\top\bm\beta^{{\bm i},k} \\
      \subjectto & \quad {\bm \alpha}^{\bm i} \in \mathbb R^K_+, \; {\bm \beta}^{\bm i} \in \mathbb R^R_+, \; {\bm \beta}^{{\bm i},k} \in \mathbb R^R_+, \; {\bm \gamma}^{{\bm i},k} \in \mathbb R^{N_\xi}, \; \forall k \in \mathcal K, \; {\bm i} \in \sets I^K \\
 & \!\! \quad  \left.  \begin{array}{l}
 \1^\top {\bm \alpha}^{\bm i}  = 1  \\

 {\bm A}^\top{\bm \beta}^{{\bm i}, k} +  {\bm w} \circ {\bm \gamma}^{{\bm i}, k} = {\bm \alpha}^{\bm i}_k \left( {\bm C^{{\bm i}_k}} {\bm x} + {\bm D^{{\bm i}_k}} {\bm w} + {\bm Q^{{\bm i}_k}} {\bm y}^k \right) \quad \forall k \in \mathcal K \\
 {\bm A}^\top{\bm \beta}^{\bm i}  =  \displaystyle \sum_{k \in \mathcal K}  {\bm w} \circ {\bm \gamma}^{{\bm i}, k}
 \end{array} \quad \right\} \quad  \forall {\bm i} \in \sets I^K.
\end{array}
$$
Finally, we write the above problem as a single minimization using an epigraph formulation, as follows
\begin{equation}
\begin{array}{cl}
     \displaystyle \minimize  & \quad \tau \\
      \subjectto & \quad \tau \in \reals, \; {\bm \alpha}^{\bm i} \in \mathbb R^K_+, \; {\bm \beta}^{\bm i} \in \mathbb R^R_+, \; {\bm \beta}^{{\bm i},k} \in \mathbb R^R_+, \; {\bm \gamma}^{{\bm i},k} \in \mathbb R^{N_\xi}, \; \forall k \in \mathcal K, \; {\bm i} \in \sets I^K \\
 & \!\! \quad  \left.  \begin{array}{l}
  \tau \; \geq \; \bm b^\top\bm \beta^{\bm i} + \sum_{k\in\sets K}\bm b^\top\bm\beta^{{\bm i},k} \\
  \1^\top {\bm \alpha}^{\bm i}  = 1  \\
 {\bm A}^\top{\bm \beta}^{{\bm i}, k} +  {\bm w} \circ {\bm \gamma}^{{\bm i}, k} = {\bm \alpha}^{\bm i}_k \left( {\bm C^{{\bm i}_k}} {\bm x} + {\bm D^{{\bm i}_k}} {\bm w} + {\bm Q^{{\bm i}_k}} {\bm y}^k \right) \quad \forall k \in \mathcal K \\
 {\bm A}^\top{\bm \beta}^{\bm i}  =  \displaystyle \sum_{k \in \mathcal K}  {\bm w} \circ {\bm \gamma}^{{\bm i}, k}
 \end{array} \quad \right\} \quad  \forall {\bm i} \in \sets I^K.
\end{array}
\label{eq:obj_Kadapt_mo_4}
\end{equation}
The claim then follows by grouping the outer minimization problem in~\eqref{eq:Kadapt_min_max_regret_mo} with the minimization problem in~\eqref{eq:obj_Kadapt_mo_4}.
\Halmos
\endproof


\section{Proofs of Statements in Section~\ref{sec:multistage}}
\label{sec:EC_multistage}


While Problem~\eqref{eq:endo_multistage_2_Kadapt} appears significantly more complicated than its two-stage counterpart, it can be brought to a min-max-min form at the cost of lifting the dimension of the uncertainty, as shown in the following lemma.

\begin{lemma}
Problem~\eqref{eq:endo_multistage_2_Kadapt} is equivalent to the two-stage robust problem
\begin{equation}\renewcommand{\arraystretch}{1.5}
    \begin{array}{cl}
         \displaystyle \min & \quad \displaystyle \max_{ \begin{smallmatrix} {\bm \xi}^{t,k_1 \cdots k_t} \in \Xi(  {\bm w}^{t-1,k_1 \ldots k_{t-1}}, {\bm \xi}^{t-1,k_1 \cdots k_{t-1}} )  \\ \forall k_1, \ldots, k_t \in \sets K, \; t \in \sets T \end{smallmatrix} } \cdots \;\; \\
         & \qquad \qquad \qquad \cdots \displaystyle \min_{k_1,\ldots,k_T \in \sets K} \;\; \sum_{t \in \sets T}  ({\bm \xi}^{T,k_1\cdots k_T})^\top {\bm D}^t  \; {\bm w}^{t,k_1 \ldots k_t} + ({\bm \xi}^{T,k_1\cdots k_T})^\top {\bm Q}^t  \; {\bm y}^{t,k_1 \ldots k_t} \\
         \st &  \quad {\bm y}^{t,k_1 \ldots k_t} \in \sets Y_t, \; {\bm w}^{t,k_1 \ldots k_t} \in \sets W_t \quad \forall t \in \sets T, \; k \in \sets K \\
 & \quad {\bm w}^{t,k_1 \ldots k_t} \geq {\bm w}^{t-1,k_1 \ldots k_{t-1}} \quad \forall t \in \sets T, \; k_t\in \sets K, \; k_{t-1} \in \sets K \\
 & \quad \displaystyle \sum_{t\in \sets T}{\bm V}^{t} {\bm w}^{t,k_1 \ldots k_t} + {\bm W}^{t} {\bm y}^{t,k_1 \ldots k_t} \leq {\bm h} \quad \forall k_1 , \ldots, k_T \in \sets K \\
 & \quad  {\bm w}^{1,k_1} = {\bm w}^{1,k_1'} \quad \forall k_1, k_1' \in \sets K, 
    \end{array}
\label{eq:endo_multistage_3_Kadapt}
\end{equation}
\label{lem:Kadapt_min_max_min_multistage}
\end{lemma}
The proof of Lemma~\ref{lem:Kadapt_min_max_min_multistage} follows directly by applying the proof of Lemma~\ref{lem:Kadapt_cstr_min_max_min} iteratively, starting at the last period.


\proof{Proof of Theorem~\ref{thm:Kadapt_obj_MBLP_multistage}} For any fixed $\{ {\bm w}^{t,k_1 \ldots k_t} \}_{t \in \sets T, k_1,\ldots k_t \in \sets K}$ and $\{ {\bm y}^{t,k_1 \ldots k_t} \}_{t \in \sets T, k_1,\ldots k_t \in \sets K}$, the inner problem in the objective of Problem~\eqref{eq:endo_multistage_3_Kadapt} can be written in epigraph form as
\begin{equation*}
    \begin{array}{cl}
         \displaystyle \maximize & \quad \displaystyle \tau \\
         \subjectto & \quad \tau \in \reals, \; {\bm \xi}^{T,k_1 \cdots k_T} \in \Xi^T(  {\bm w}^{1,k_1},\ldots, {\bm w}^{T-1,k_1 \ldots k_{T-1}} ) \quad \forall  k_1, \ldots, k_T \in \sets K \\
         & \quad \displaystyle \tau \; \leq \;  \sum_{t \in \sets T}  ({\bm \xi}^{T,k_1\cdots k_T})^\top {\bm D}^t  \; {\bm w}^{t,k_1 \ldots k_t} + ({\bm \xi}^{T,k_1\cdots k_T})^\top {\bm Q}^t  \; {\bm y}^{t,k_1 \ldots k_t} \quad \forall k_1,\ldots, k_T \in \sets K.
    \end{array}
\end{equation*}
From the definition of $\Xi^T(\cdot)$ in Lemma~\ref{lem:Kadapt_min_max_min_multistage}, the above problem can be equivalently written as
\begin{equation*}
    \begin{array}{cl}
         \displaystyle \maximize & \quad \displaystyle \tau \\
         \subjectto & \quad \tau \in \reals, \; {\bm \xi}^{t,k_1 \cdots k_t} \in \Xi \quad \forall t \in \sets T, \; k_1,\ldots, k_t \in \sets K\\
         & \quad \displaystyle \tau \; \leq \;  \sum_{t \in \sets T}  ({\bm \xi}^{T,k_1\cdots k_T})^\top {\bm D}^t  \; {\bm w}^{t,k_1 \ldots k_t} + ({\bm \xi}^{T,k_1\cdots k_T})^\top {\bm Q}^t  \; {\bm y}^{t,k_1 \ldots k_t} \quad \forall k_1,\ldots, k_T \in \sets K \\
         & \quad {\bm w}^{t-1,k_1 \ldots k_{t-1}} \circ {\bm \xi}^{t,k_1 \cdots k_t} = {\bm w}^{t-1,k_1 \ldots k_{t-1}} \circ {\bm \xi}^{t-1,k_1\cdots k_{t-1}}  \quad \forall t \in \sets T  \backslash \{1\}, \; k_1,\ldots, k_t \in \sets K.
    \end{array}
\end{equation*}
Writing the set $\Xi$ explicitly yields
\begin{equation*}
    \begin{array}{cl}
         \displaystyle \maximize & \quad \displaystyle \tau \\
         \subjectto & \quad \tau \in \reals, \; {\bm \xi}^{t,k_1 \ldots k_t} \in \reals^{N_\xi} \quad \forall t \in \sets T, \; k_1,\ldots, k_t \in \sets K \\
         & \quad \displaystyle \tau \; \leq \;   \sum_{t \in \sets T} \left( {\bm D}^t  \; {\bm w}^{t,k_1 \ldots k_t} + {\bm Q}^t  \; {\bm y}^{t,k_1 \ldots k_t} \right)^\top {\bm \xi}^{T,k_1\cdots k_T} \quad \forall k_1,\ldots, k_T \in \sets K \\
         & \quad {\bm A} {\bm \xi}^{t,k_1 \ldots k_t} \leq {\bm b} \quad \forall t \in \sets T, \; k_1,\ldots, k_t \in \sets K  \\
         & \quad {\bm w}^{t-1,k_1 \ldots k_{t-1}} \circ {\bm \xi}^{t,k_1 \ldots k_t} = {\bm w}^{t-1,k_1 \ldots k_{t-1}} \circ {\bm \xi}^{t-1,k_1\cdots k_{t-1}}  \quad \forall t \in \sets T  \backslash \{1\}, \; k_1,\ldots, k_t \in \sets K.
    \end{array}
\end{equation*}
The dual of this problem reads
\begin{equation*}
    \begin{array}{cl}
         \displaystyle \minimize & \quad \displaystyle \sum_{t\in \sets T} \sum_{k_1 \in \sets K} \cdots \sum_{k_t \in \sets K} {\bm b}^\top {\bm \beta}^{t,k_1\cdots k_t}   \\
         \subjectto & \quad {\bm \alpha} \in \reals_+^{K^T} , \; {\bm \beta}^{t, k_1 \cdots k_t} \in \reals_+^{R}, \; {\bm \gamma}^{t, k_1 \cdots k_t} \in \reals^{N_\xi}, \; t \in \sets T, \; k_1, \ldots, k_t \in \sets K \\
         & \quad \1^\top {\bm \alpha} = 1 \\
         & \quad \displaystyle {\bm A}^\top {\bm \beta}^{1,k_1} = \sum_{k_2 \in \sets K} {\bm w}^{1,k_1} \circ {\bm \gamma}^{2,k_1 k_2} \quad \forall k_1 \in \sets K  \\
         & \quad {\bm A}^\top {\bm \beta}^{t,k_1\cdots k_t} + {\bm w}^{t-1,k_1 \ldots k_{t-1}} \circ {\bm \gamma}^{t,k_1 \cdots k_t}  = \\
         & \qquad \qquad \qquad \displaystyle \sum_{k_{t+1} \in \sets K} {\bm w}^{t,k_1 \cdots k_t} \circ {\bm \gamma}^{t+1,k_1\cdots k_{t+1}} \displaystyle  \quad \forall t\in \sets T \backslash \{ 1, T\} , \; k_1,\ldots,k_t \in \sets K \\
         & \quad {\bm A}^\top {\bm \beta}^{T,k_1\cdots k_T} + {\bm w}^{T-1,k_1\cdots k_{T-1}} \circ {\bm \gamma}^{T,k_1 \cdots k_T} \displaystyle  \; = \; \\
         & \qquad \qquad \qquad {\bm \alpha}_{k_1 \cdots k_T}  \sum_{t \in \sets T} \left( {\bm D}^t  \; {\bm w}^{t,k_1\cdots k_t} + {\bm Q}^t  \; {\bm y}^{t,k_1\cdots k_t} \right) \quad \forall k_1, \ldots k_T .
    \end{array}
\end{equation*}
Moreover, strong duality applies by virtue of the compactness of $\Xi$. Merging the problem above with the outer minimization problem in~\eqref{eq:endo_multistage_3_Kadapt} yields
\begin{equation*}
    \begin{array}{cl}
         \displaystyle \minimize & \quad \displaystyle \sum_{t\in \sets T} \sum_{k_1 \in \sets K} \cdots \sum_{k_t \in \sets K} {\bm b}^\top {\bm \beta}^{t,k_1\cdots k_t}   \\
         \subjectto & \quad {\bm \alpha} \in \reals_+^{K^T} , \; {\bm \beta}^{t, k_1 \cdots k_t} \in \reals_+^{R}, \; {\bm \gamma}^{t, k_1 \cdots k_t} \in \reals^{N_\xi}, \; t \in \sets T, \; k_1, \ldots, k_t \in \sets K \\
         &\quad {\bm y}^{t,k_1 \ldots k_t} \in \sets Y_t, \; {\bm w}^{t,k_1 \ldots k_t} \in \sets W_t \quad \forall t \in \sets T, \; k \in \sets K \\
         & \quad \1^\top {\bm \alpha} = 1 \\
         & \quad \displaystyle {\bm A}^\top {\bm \beta}^{1,k_1} = \sum_{k_2 \in \sets K} {\bm w}^{1,k_1} \circ {\bm \gamma}^{2,k_1 k_2} \quad \forall k_1 \in \sets K  \\
         & \quad {\bm A}^\top {\bm \beta}^{t,k_1\cdots k_t} + {\bm w}^{t-1,k_1\cdots k_{t-1}} \circ {\bm \gamma}^{t,k_1 \cdots k_t}  =  \\
         & \qquad \qquad \qquad \displaystyle \sum_{k_{t+1} \in \sets K} {\bm w}^{t,k_1 \cdots k_t} \circ {\bm \gamma}^{t+1,k_1\cdots k_{t+1}} \displaystyle  \quad \forall t\in \sets T \backslash \{ 1, T\} , \; k_1,\ldots,k_t \in \sets K \\
         & \quad {\bm A}^\top {\bm \beta}^{T,k_1\cdots k_T} + {\bm w}^{T-1,k_1 \cdots k_{T-1}} \circ {\bm \gamma}^{T,k_1 \cdots k_T} \displaystyle  \; = \; \\
         & \qquad \qquad \qquad {\bm \alpha}_{k_1 \cdots k_T}  \sum_{t \in \sets T} \left( {\bm D}^t  \; {\bm w}^{t,k_1 \cdots k_t} + {\bm Q}^t  \; {\bm y}^{t,k_1 \cdots k_t} \right) \quad \forall k_1, \ldots k_T \in \sets K  \\
        & \quad {\bm w}^{t,k_1 \cdots k_t} \geq {\bm w}^{t-1,k_1 \cdots k_{t-1}} \quad \forall t \in \sets T, \; k_t\in \sets K, \; k_{t-1} \in \sets K \\
        & \quad \displaystyle \sum_{t\in \sets T}{\bm V}^{t} {\bm w}^{t,k_1 \cdots k_t} + {\bm W}^{t} {\bm y}^{t,k_1 \cdots k_t} \leq {\bm h} \quad \forall k_1, \ldots, k_T \in \sets K,
    \end{array}
\end{equation*}
and our proof is complete. \Halmos
\endproof


\section{Proofs of Statements in Section~\ref{sec:EC_pwl}}

\proof{Proof of Proposition~\ref{prop:algo_correct}}
Since $({\bm x},{\bm w},\{{\bm y}^k\}_{k\in \sets K})$ is feasible in the relaxed master problem~\eqref{eq:ccg_master}, it follows that ${\bm x}\in \sets X$, ${\bm w} \in \sets W$, and ${\bm y}^k \in \sets Y$, $k\in \sets K$. Thus, $({\bm x},{\bm w},\{{\bm y}^k\}_{k\in \sets K})$ is feasible in Problem~\eqref{eq:min_max_regret_general_mo_Kadapt}. An inspection of the Proof of Theorem~\ref{thm:Kadapt_mo_MBLP} reveals that the objective value of $({\bm x},{\bm w},\{{\bm y}^k\}_{k\in \sets K})$ in Problem~\eqref{eq:min_max_regret_general_mo_Kadapt} is given by the optimal value of Problem~\eqref{eq:obj_Kadapt_mo}. The proof then follows by noting that Problems~\eqref{eq:obj_Kadapt_mo} and~\eqref{eq:ccg_feas} are equivalent.
\Halmos
\endproof


\proof{Proof of Lemma~\ref{lem:algo_correct2}}
\begin{enumerate}[label=\emph{(\roman*)}]
\item By virtue of Proposition~\ref{prop:algo_correct}, it follows that $\theta\geq \tau$.
\item Suppose that $\theta = \tau$ and that there exists ${\bm i} \in \sets I^K$ such that Problem~\eqref{eq:ccg_subproblem} is infeasible. This implies that there exists ${\bm i} \in \sets I^K$ such that $\tau$ is strictly smaller than the optimal objective value of
\begin{equation}
\begin{array}{cl}
     \displaystyle \minimize  & \quad \bm b^\top\bm \beta^{\bm i} + \sum_{k\in\sets K}\bm b^\top\bm\beta^{{\bm i},k} \\
      \subjectto & \quad {\bm \alpha}^{\bm i} \in \mathbb R^K_+, \; {\bm \beta}^{\bm i} \in \mathbb R^R_+, \; {\bm \beta}^{{\bm i},k} \in \mathbb R^R_+, \; {\bm \gamma}^{{\bm i},k} \in \mathbb R^{N_\xi}, \; \forall k \in \mathcal K \\
 & \!\! \quad   \begin{array}{l}
  \1^\top {\bm \alpha}^{\bm i}  = 1  \\
 {\bm A}^\top{\bm \beta}^{{\bm i}, k} +  {\bm w} \circ {\bm \gamma}^{{\bm i}, k} = {\bm \alpha}^{\bm i}_k \left( {\bm C^{{\bm i}_k}} {\bm x} + {\bm D^{{\bm i}_k}} {\bm w} + {\bm Q^{{\bm i}_k}} {\bm y}^k \right) \quad \forall k \in \mathcal K \\
 {\bm A}^\top{\bm \beta}^{\bm i}  =  \displaystyle \sum_{k \in \mathcal K}  {\bm w} \circ {\bm \gamma}^{{\bm i}, k}.
 \end{array}
\end{array}
\label{eq:ccg_subproblem_2}
\end{equation}
Equivalently, by dualizing this problem, we conclude that there exists ${\bm i} \in \sets I^K$ such that $\tau$ is strictly smaller than the optimal objective value of
\begin{equation}
\begin{array}{cl}
    \maximize & \quad \theta' \\
    \subjectto & \quad \theta' \in\mathbb{R},\; \overline{\bm \xi}\in\Xi,\; {\bm \xi^k} \in \Xi({\bm w} , \overline {\bm \xi}),\; \forall k\in\sets K \\
    & \quad \theta' \; \leq \; (\bm\xi^k)^\top {\bm C}^{{\bm i}_k}\;{\bm x} + ({\bm \xi}^k)^\top {\bm D}^{{\bm i}_k} \; {\bm w} + ({\bm \xi}^k)^\top {\bm Q}^{{\bm i}_k} \; {\bm y}^k \quad\forall k\in\sets K.
\end{array}
\label{eq:obj_Kadapt_mo_tmp}
\end{equation}
Since Problem~\eqref{eq:obj_Kadapt_mo_tmp} lower bounds Problem~\eqref{eq:ccg_feas} with optimal objective value $\theta$, we conclude that $\tau<\theta$, a contradiction.
\item Suppose that $\theta>\tau$ and let ${\bm i}$ be defined as in the premise of the lemma. Then, ${\bm i}$ is optimal in \eqref{eq:obj_Kadapt_mo_3} with associated optimal objective value $\theta$. This implies that the optimal objective value of Problem~\eqref{eq:ccg_subproblem_2} is $\theta$. Since $\theta > \tau$, this implies that subproblem~\eqref{eq:ccg_subproblem} is infeasible, which concludes the proof.
\end{enumerate}
We have thus proved all claims. \Halmos
\endproof

\proof{Proof of Theorem~\ref{thm:algo_converges}} First, note that finite termination is guaranteed since at each iteration, either ${\rm{UB}}-{\rm{LB}} \leq \delta$ (in which case the algorithm terminates) or a new set of constraints (indexed by the infeasible index~${\bm i}$) is added to the master problem~\eqref{eq:ccg_master}, see Lemma~\ref{lem:algo_correct2}. Since the set of all indices, $\mathcal I^K$, is finite, the algorithm will terminate in a finite number of steps. Second, by construction, at any iteration of the algorithm, $\tau$ (i.e., ${\rm{LB}}$) provides a lower bound on the optimal objective value of the problem. On the other hand, the returned (feasible) solution has as objective value $\theta$ (i.e., ${\rm{UB}}$). Since the algorithm only terminates if ${\rm{UB}}-{\rm{LB}} \leq \delta$, we are guaranteed that, at termination, the returned solution will have an objective value that is within~$\delta$ of the optimal objective value of the problem. This concludes the proof.
\Halmos
\endproof

\proof{Proof of Observation~\ref{obs:min_max_regret_as_mo}} Suppose that $\sets X := \{ {\bm x}: \1^\top {\bm x}=1\}$, $\sets W := \{ {\bm w}: \1^\top {\bm w}=1\}$, and $\sets Y := \{ {\bm y} : \1^\top {\bm y}=1\}$. Then,
$$
\begin{array}{clccl}
 \max    & \quad {\bm \xi}^\top {\bm C} \; {\bm x}' + {\bm \xi}^\top {\bm D} \; {\bm w}' + {\bm \xi}^\top {\bm Q} \; {\bm y}'
 \quad \qquad  & = & \quad \displaystyle \max_{i,j,k}    & \;\;  \left\{ {\bm \xi}^\top {\bm C} \; \1_i + {\bm \xi}^\top {\bm D} \; \1_j + {\bm \xi}^\top {\bm Q} \; \1_k \right\}. \\
 \st   & \quad {\bm x}' \in \sets X, \; {\bm w}' \in \sets W,\; {\bm y}' \in \sets Y
\end{array}
$$
Thus, in this case, the objective function is expressible in the form~\eqref{eq:regret_max_objective} and the claim follows. An analogous argument can be made if ${\bm C}={\bm 0}$, ${\bm D}={\bm 0}$, or ${\bm Q}={\bm 0}$. 
\Halmos
\endproof

\section{$K$-Adaptable Binary Linear Decision Rule and Lifting}

Throughout Sections~\ref{sec:kadaptability_objective},~\ref{sec:kadaptability_constraint}, and~\ref{sec:pwl_convex_objective}, we showed that if the original two-stage robust optimization problem with decision-dependent information discovery presents binary first- and second-stage decisions, then it can be written equivalently as an MBLP, while it is a bilinear problem if any of the decision variables are real valued. We also noted that there exists new off-the-shelf solvers for tackling such bilinear problems. In this section, we propose an alternative conservative solution approach applicable to problems that have real-valued wait-and-see decisions.


Consider the following variant of Problem~\eqref{eq:endo_3} where the wait-and-see decisions ${\bm y}$ are real-valued ($\sets Y \subseteq \reals^{N_y})$ and its coefficients in the objective function are deterministic.
\begin{equation}
    \begin{array}{cl}
         \min & \;\; \displaystyle \max_{\overline {\bm \xi} \in \Xi} \;\;\min_{ {\bm y} \in \sets Y } \;  \left\{ \max_{ {\bm \xi} \in \Xi({\bm w},\overline {\bm \xi}) } \; \; {\bm \xi}^\top {\bm C} \; {\bm x} + {\bm \xi}^\top {\bm D} \; {\bm w} + {\bm q}^\top {\bm y} \; : \; {\bm T} {\bm x} + {\bm V} {\bm w} + {\bm W}{\bm y} \leq {\bm H}{\bm \xi}  \; \; \; \; \forall {\bm \xi} \in \Xi({\bm w},\overline {\bm \xi}) \right\}  \\
         \st & \;\; {\bm x} \in \sets X, \; {\bm w} \in \sets W, 
    \end{array}
\label{eq:endo_cont_rec}
\end{equation}
where ${\bm q} \in \reals^{N_y}$. In the spirit of the linear decision rule approximation approach proposed in the stochastic and robust optimization literature, see e.g.,~\cite{AdjustableRSols_uncertainLP,kuhn_primal_dual_rules,Bodur2018}, we propose to restrict the recourse decisions ${\bm y}$ to those that are expressible as
$$
{\bm y}({\bm \xi}) = {\bm Y}{\bm \xi},
$$
for some matrix ${\bm Y} \in \{0,1\}^{N_y \times N_\xi}$. We refer to this approximation as the \emph{binary linear decision rule}.

Under this approximation, Problem~\eqref{eq:endo_cont_rec} is equivalent to  
$$
    \begin{array}{cl}
         \minimize & \;\; \displaystyle \max_{\overline {\bm \xi} \in \Xi({\bm w})} \;\;\min_{ {\bm Y} \in \{0,1\}^{N_y \times N_\xi} } \;  \left\{
         \begin{array}{cl}
         \displaystyle \max_{ {\bm \xi} \in \Xi({\bm w},\overline {\bm \xi}) } & \; \; {\bm \xi}^\top {\bm C} \; {\bm x} + {\bm \xi}^\top {\bm D} \; {\bm w} + {\bm q}^\top {\bm Y}{( {\bm w} \circ {\bm \xi} )} \\ \st & {\bm T} {\bm x} + {\bm V} {\bm w} + {\bm W}{\bm Y}{( {\bm w} \circ {\bm \xi} )} \leq {\bm H}{\bm \xi}  \; \; \; \; \forall {\bm \xi} \in \Xi({\bm w},\overline {\bm \xi}) \\
         & {\bm Y} {( {\bm w} \circ {\bm \xi} )} \in \sets Y  \; \; \; \; \forall {\bm \xi} \in \Xi({\bm w},\overline {\bm \xi}) 
         \end{array}  \right\}  \\
         \subjectto & \;\; {\bm x} \in \sets X, \; {\bm w} \in \sets W.
    \end{array}
$$
This problem can be written in the form~\eqref{eq:endo_3} with the matrix ${\bm W}$ being affected by uncertainty (left-handside uncertainty), after linearizing the product of ${\bm Y}$ and ${\bm w}$. From Remark~\ref{rmk:random_rec}, our $K$-adaptability approximation framework applies in this case too. It results in a number~$K$ of binary linear contingency plans or \emph{operating regimes}. For this reason, we refer to it as the \emph{$K$-adaptable binary linear decision rule}. This approximation is very natural since it enables us to choose between several modes of operation for the wait-and-see decisions. 

\end{document}